\documentclass[11pt]{article}
\usepackage[top=1.1in, left=1in, right=1in, bottom=1.1in]{geometry}
\usepackage{amsmath,amssymb,amsthm,amsfonts,verbatim}
\usepackage{enumerate}
\usepackage{url}
\usepackage{amsmath,amssymb}
\usepackage{graphicx}
\usepackage{color}
\usepackage{setspace}
\usepackage{xcolor}

\newtheorem{theorem}{Theorem}

\def\mr{\mathbb{R}}
\def\mx{\mathbf{x}}
\def\mb{\mathbf{b}}

\def\mv{\mathbf{v}}

\title{Ordered Line Integral Methods \\
for Computing the Quasi-potential}

\author{Daisy Dahiya$^1$  and  Maria Cameron$^1$}

\begin{document}

\maketitle
\abstract{
The quasi-potential is a key function in the Large Deviation Theory. It characterizes the difficulty of the escape 
from the neighborhood of an attractor of a stochastic non-gradient dynamical system due to the influence of small white noise.
It also gives an estimate of the invariant probability distribution in the neighborhood of the attractor up { to} the exponential order.
We present a new family of methods for computing the quasi-potential on a regular mesh 
named the Ordered Line Integral Methods (OLIMs).  In comparison with the first proposed quasi-potential finder 
based on the Ordered Upwind Method (OUM) (Cameron, 2012), the new methods are 1.5 to { 4 } times faster, 
can produce error two to three orders of magnitude smaller, and may exhibit faster convergence.
Similar to the OUM, OLIMs employ the dynamical programming principle. Contrary to it, 
they $(i)$ have an optimized strategy for the use of computationally expensive { triangle} updates leading to a notable speed-up, and $(ii)$
directly solve local minimization problems using quadrature rules instead of solving the corresponding 
Hamilton-Jacobi-type equation by the first order finite difference upwind scheme. 
The OLIM with the right-hand quadrature rule is equivalent to OUM.
The use of higher order quadrature rules in local minimization problems dramatically boosts up the accuracy of OLIMs. 
We offer a detailed discussion  on the origin of numerical errors in OLIMs and propose rules-of-thumb for the choice of the 
important parameter, the update factor, in the OUM and OLIMs.
Our results are supported by extensive numerical tests on two challenging 2D examples. 


\footnotetext[1]{ Daisy Dahiya: {\tt ddahiya@math.umd.edu}; 
Maria Cameron {\tt cameron@math.umd.edu};\\
Department of Mathematics, University of Maryland, College Park, MD 20742}

}

%%%%%%%%%%%

\section{Introduction}
\label{sec:intro}
The quasi-potential is a key concept of Large Deviation Theory \cite{FW}. 
Let us consider a nongradient stochastic differential equation (SDE) of the form
\begin{equation}
\label{sde1}
d\mx=\mb(\mx)dt+\sqrt{\epsilon} dW, \quad \mx\in\mr^d
\end{equation}
where $\mb(\mx)$ is a twice continuously differentiable vector field, 
$W$ is the standard $d$-dimensional Brownian motion, and $\epsilon$ is a small parameter.
We assume that the corresponding deterministic system $\dot{\mx}=\mb(\mx)$ 
has a finite number of attractors.
Let $A\subset \mr^d$ be an attractor and $\mx$ be a point. 
{ The quasi-potential at $\mx$ with respect to the
attractor $A$ can be defined as (details are provided in Appendix A)}
 \begin{equation}
 \label{Qpot_def}
 U_A(\mx):=\inf_{\psi:~\psi(0)\in A,~\psi(L) = \mx}\left\{S(\psi): = \int_0^L\left[\|\psi'\|\|\mb(\psi)\| - \psi'\cdot\mb(\psi)\right]ds\right\}.
 \end{equation}
In Eq. \eqref{Qpot_def}, $\|\cdot\|$ is the 2-norm, the functional $S(\psi)$ is the geometric action \cite{hey1,hey2}, 
 $\psi(s)$ is an absolutely continuous path from $A$ to $\mx$ parametrized by its arc length (i.e., { $\|\psi'\|=1$}), 
 $L$ is its length which can be infinite. 
 According to the Large Deviation Theory, the expected first passage time  from a small neighborhood of the attractor $A$
 to a small neighborhood of the point $\mx$ 
 lying in its basin of attraction
 is logarithmically equivalent to $\exp(U_A(\mx)/\epsilon)$ \cite{FW}. 
 The curve corresponding to the maximum likelihood path to $\mx$ from $A$  is the global minimizer of the  geometric action $S(\psi)$ \cite{FW,hey1,hey2}. 
The invariant probability distribution in the basin of attraction of $A$ is logarithmically equivalent to $\exp(-U_A(\mx)/\epsilon)$ \cite{FW}.
% Some details regarding the properties of the quasi-potential as a function of $\mx$ are given in \cite{quasi}.
 A nice and visual account on the significance of the quasi-potential as a tool for quantifying the stability of attractors  is 
 found in \cite{NA}. 
 
 Nongradient SDEs of the form \eqref{sde1} often arise in modeling of ecological and biological systems. 
 A number of models of
 population dynamics were analyzed in \cite{NA} using the quasi-potential.
 The quasi-potential for a genetic switching model was constructed in \cite{tiejun1}.

Therefore, the quasi-potential is of crucial  importance for the quantification of the dynamics of systems evolving 
according to SDE \eqref{sde1} with a small noise.
Unfortunately, it is not readily available for nongradient SDEs and can be found analytically only in special cases. For example, 
analytic formulas for the quasi-potential are available for linear SDEs \cite{quasi,chen}. 
For general nonlinear SDEs, finding the quasi-potential is a difficult task. 
In any bounded neighborhood of an attractor $A$ of $\dot{\mx}=\mb(\mx)$,
 { the quasi-potential} is a Lipschitz-continuous but not necessarily differentiable function as it is for the Maier-Stein 
\cite{maier-stein,quasi} model. One can show that the quasi-potential $U_A(\mx)$, i.e.,  the 
solution of the functional minimization problem \eqref{Qpot_def}, 
is also a viscosity solution of a Hamilton-Jacobi-type PDE \cite{quasi}
\begin{equation}
\label{HJ}
\|\nabla U_A(\mx)\|^2 + 2\mb(\mx)\cdot\nabla U_A(\mx) = 0,\quad U_A(A) = 0.
\end{equation}
Eq. \eqref{HJ} has at least two viscosity solutions due to the fact that $A$ is an attractor and  $\nabla U(\mx) = 0$ for all $\mx\in A$:  
$U_A(\mx)\equiv 0$ and $U_A(\mx)$ is the quasi-potential (i.e, the solution of  Eq. \eqref{Qpot_def}). 
We refer an interested reader to Refs. \cite{crandall} and \cite{ishii} 
where the concept of viscosity solution was introduced and the questions of existence and uniqueness
of solution of Eq. \eqref{HJ} were investigated respectively.

%%%%%
\subsection{Background}
Minimum Action Paths (MAPs), i.e.,  minimizers of the geometric action $S(\psi)$ ((Eq. \eqref{Qpot_def}) can be found numerically by path-based methods.
The Geometric Minimum Action Method (GMAM) \cite{hey1,hey2} and the Adaptive Minimum Action Method (AMAM) \cite{xiang1,xiang2} 
iteratively update paths
connecting given initial and final points starting from user-provided initial guesses
so that the paths approach a local minimizer of the geometric action $S(\psi)$ (Eq. \eqref{Qpot_def}) 
and a local minimizer of the original Freidlin-Wentzell action (Eq. \eqref{FWA}) respectively.
The advantage of the path-based methods is that they are suitable for high- and even infinite-dimensional systems, i.e., 
stochastic partial differential equations (SPDEs) (e.g., \cite{hey1,grafke}). 
However, the found action minimizer is biased by the initial { path} and hence might not be the global minimizer, and
it can be inexact due to slow convergence of the iterative process and numerical effects 
in the case if the actual MAP 
exhibits complex behavior. 
If  the global action minimizer is found, the quasi-potential can be calculated along it.

{ Finding the quasi-potential in the whole  region surrounding the attractor has important advantages over the search for MAPs. 
First, the quasi-potential allows us to estimate the invariant probability density near the attractor up to exponential order \cite{FW,quasi}.
Second, suppose a system has more than two attractors. Once the quasi-potential is computed  with respect to an attractor $A$ 
in a large enough region, the most likely 
escape set $E$ from the basin of attraction  of $A$ is automatically detected as the quasi-potential is constant 
along any trajectory running from $E$ to another attractor.
The escape set $E$ can be a saddle point, a hyperbolic periodic orbit, or a more complex set of points in the phase space.
Once the quasi-potential is computed, 
the global minimizer of the geometric action connecting the 
attractor $A$ with the escape set $E$ is guaranteed to be found by straightforward numerical  integration 
\cite{quasi} (see Section \ref{sec:tests},  Fig. \ref{fig:MAP} and its caption).
On the contrary, path-based methods might have wrong endpoints, e.g., the most likely escape from $A$ is to the attractor $B$, while
a path-based method is set up to seek action minimizers connecting attractors $A$ and $C$.
Furthermore, even if the endpoints are identified correctly,  path-based methods might converge to local minimizers.

On the other hand, it is clear that the computation of the quasi-potential in whole 
regions of the phase space is limited to low-dimensional systems.}

To the best of our knowledge, the first attempt to compute the quasi-potential on a 
regular mesh  was undertaken in \cite{quasi}. 
It was done for 2D systems { for both possible types of attractors, asymptotically stable equilibrium and stable limit cycle.
The proposed numerical technique was} an adjustment of the Ordered Upwind Method  (OUM) \cite{oum1,oum2} 
for the case where the  { anisotropy} coefficient $\Upsilon$ was unbounded. (What is the anisotropy coefficient will be explained right below.)
The original OUM was designed for 
solving Hamilton-Jacobi equations of the form $\frac{1}{F\left(\mx, \hat{a}\right)}\|\nabla u(\mx)\| = 1$, 
where $\hat{a} \equiv \frac{\nabla u}{\|\nabla u\|}$, with the boundary condition $u(A) = 0$ where $A$ can be a point or a curve.
The function $F(\mx,\hat{a})$, the speed of front propagation in the normal direction, 
was assumed to be bounded: $0<F_{\min} \le F(\mx,\hat{a}) \le F_{\max}<\infty$.
The Hamiltonian $\|\nabla u\| /F(\mx, \hat{a})$ was assumed to be convex, Lipschitz-continuous, and homogeneous of degree one in $\nabla u$. 
The anisotropy coefficient $\Upsilon$ is the ratio of the maximal and the minimal values of the speed function: $\Upsilon:=F_{\max}/F_{\min}$.
One can think of $u(\mx)$ being the minimal possible travel time from $A$ 
to the point $\mx$ in the following associated 
control problem. Suppose an astronaut on Mars (a 2D manifold with no roads) needs to minimize  his traveltime from $A$ to $\mx$.
He can pick the direction of motion $\hat{q}$ of his rover at every moment of time.
The speed of the rover is the function $f(\mathbf{y},\hat{q})$ of the position $\mathbf{y}$ and the direction of motion $\hat{q}$.
The speed functions $f(\mathbf{y},\hat{q})$ and $F(\mx,\hat{a})$ 
relate via a Legendre transform. We refer an interested reader to Ref. \cite{oum2}.

The OUM \cite{oum1,oum2} employs the dynamical programming principle, 
i.e., at each step it solves the driver's optimal control program by means of 
a finite difference upwind scheme for $\frac{1}{F\left(\mx, \hat{a}\right)}\|\nabla u(\mx)\| = 1$.
The role of the anisotropy coefficient $\Upsilon$ is very important.  The value of $u$ at a mesh point $\mx$ can be updated { only from} the mesh points lying within the distance of $\Upsilon h$ from it ($h = \max\{h_1,h_2\}$, $h_1$, $h_2$ are the mesh steps in $x_1$ and $x_2$
respectively). 

 Eq. \eqref{HJ}, the Hamilton-Jacobi equation for the quasi-potential, can be put in the form  $\frac{1}{F\left(\mx, \hat{a}\right)}\|\nabla u(\mx)\| = 1$:
\begin{equation}
\label{HJ1}
\frac{1}{-2\mb(\mx)\cdot\hat{a}}\|\nabla U_A(\mx)\| = 1,\quad \hat{a}\equiv \frac{\nabla U_A(\mx)}{\|\nabla U_A(\mx)\|}, \quad U_A(A) = 0,
\end{equation}
Hence, the front speed function is $F\left(\mx, \hat{a} \right) = \frac{1}{-2\mb(\mx)\cdot\hat{a} }$.
The corresponding speed function $f$ is the integrand in Eq. \eqref{Qpot_def}.
It { is} clear that the speed function $F\left(\mx, \hat{a} \right) = \frac{1}{-2\mb(\mx)\cdot\hat{a} }$
is unbounded at any direction $\hat{a}$ at the points $\mx$ where $\mb(\mx) = \mathbf{0}$, 
and at the directions $\hat{a}$ normal to $\mb(\mx)$ at any point $\mx$ where $\mb(\mx)\neq \mathbf{0}$.
This difficulty was overcome in \cite{quasi} by the introduction of the update factor $K$  ($K$ is a positive integer)
such that the value of $U_A$ at a mesh point $\mx$ could be  { updated}
only from the mesh points lying within the distance of $Kh$ {(i.e., within the update radius $Kh$)} from it.
The additional numerical error due to this adjustment was investigated \cite{quasi} and 
shown to decay quadratically for a fixed update radius with mesh refinement.
The OUM for finding the quasi-potential was incorporated into an R-package {\tt QPot} available at cran.r-project.org \cite{qpotR,QPot_package}.

\subsection{A brief summary of main results}
In this paper, we present a new family of methods for computing the quasi-potential on a regular rectangular mesh.
Like the OUM, these methods are based on the dynamical programming principle. Contrary to OUM, they
do not use the upwind scheme for solving local minimization problems. Moreover, they completely abandon Hamilton-Jacobi Eq. \eqref{HJ}
and refer only to the minimization problem \eqref{Qpot_def}. 
At each step, the local minimization problem is solved on the set of local straight line paths 
of length at most $Kh$\footnotemark[1]. 
\footnotetext[1]{Actually, in our codes, the maximal update length for the one-point update  is $Kh$,
while it is  $Kh + \sqrt{h_1^2+h_2^2}$ for the triangle update.}
The Lagrangian $L(\psi,\psi'):= \|\psi'\|\|\mb(\psi)\| - \psi'\cdot\mb(\psi)$ 
is integrated along them using quadrature rules.
 Therefore, we name this family of methods the Ordered Line Integral Methods. 
The names of the methods in this family reflect which quadrature rule is 
employed: 
\begin{center}
\begin{tabular}{ll}
OLIM-R: & Right-hand rectangle rule;\\
OLIM-MID: & Midpoint rule;\\
OLIM-TR: & Trapezoid rule;\\
OLIM-SIM: & Simpson's rule.\\
\end{tabular}
\end{center}
We prove a theorem showing that the solution of the minimization problem in OLIM-R method is  
equivalent to the solution of the upwind scheme in the OUM. 
Our extensive numerical experiments with the other OLIMs
show that they are significantly more accurate than the OUM.  Our least squares fits to the error formula $E=Ch^q$ indicate 
that their error constants are up to 100 times smaller and the exponents $q$ are larger than those for OUM and OLIM-R. 
To make sure, OLIMs like OUM are at most first order due to the use of linear interpolation, 
but the major portion of the numerical error, the quadrature rule error, decays 
at least quadratically with the mesh refinement, which affects the { exponents} $q$ (makes  { them} greater than 1)
obtained by the least squares fits. 
Furthermore, we propose a CPU-time-saving implementation for OLIMs that reduces the number of calls of computationally expensive triangle update.
%This measure is justified by the fact for any fixed point $\mx$, 
%the Lagrangian $L(\mx,\hat{q}):= \|\hat{q}\|\|\mb(\mx)\| - \hat{q}'\cdot\mb(\mx)$ has a single minimum.
As a result, the OLIM-R is about { four} times faster { than} the OUM. 
The other OLIMs are also faster than the OUM by some more modest factors. The graphs of the CPU times versus errors eloquently display that the OLIMs
with second and higher order quadrature rules produce at least as accurate solutions as the OLIM-R (and hence the OUM) 
at smaller by several orders of magnitude CPU times. 

So far, we have implemented the OLIMs in 2D. Our C codes {\tt OLIM\_righthand.c}, {\tt OLIM\_midpoint.c}, {\tt OLIM\_trapezoid.c} and {\tt OLIM\_simpson.c}
are posted on M. Cameron's web site \cite{mariakc}.
{ A promotion of OLIMs to 3D in underway and will be reported elsewhere in the future.}

The rest of the paper is organized as follows.
The OLIMs are described in Section \ref{sec:OLIM}.
The results of our numerical tests are presented in Section \ref{sec:tests}.
The origin of numerical errors in the OLIMs are investigated in Section \ref{sec:errors}.
The results of this work are summarized in Section \ref{sec:conclusion}.

%%%%%%%%%%%%%%%%%%%%%%%%%%%
%%%%
%%%%  O L I M
%%%%
%%%%%%%%%%%%%%%%%%%%%%%%%%%

\section{Ordered Line Integral Methods}
\label{sec:OLIM}
\begin{figure}
\begin{center}
  \includegraphics[width=0.4\textwidth]{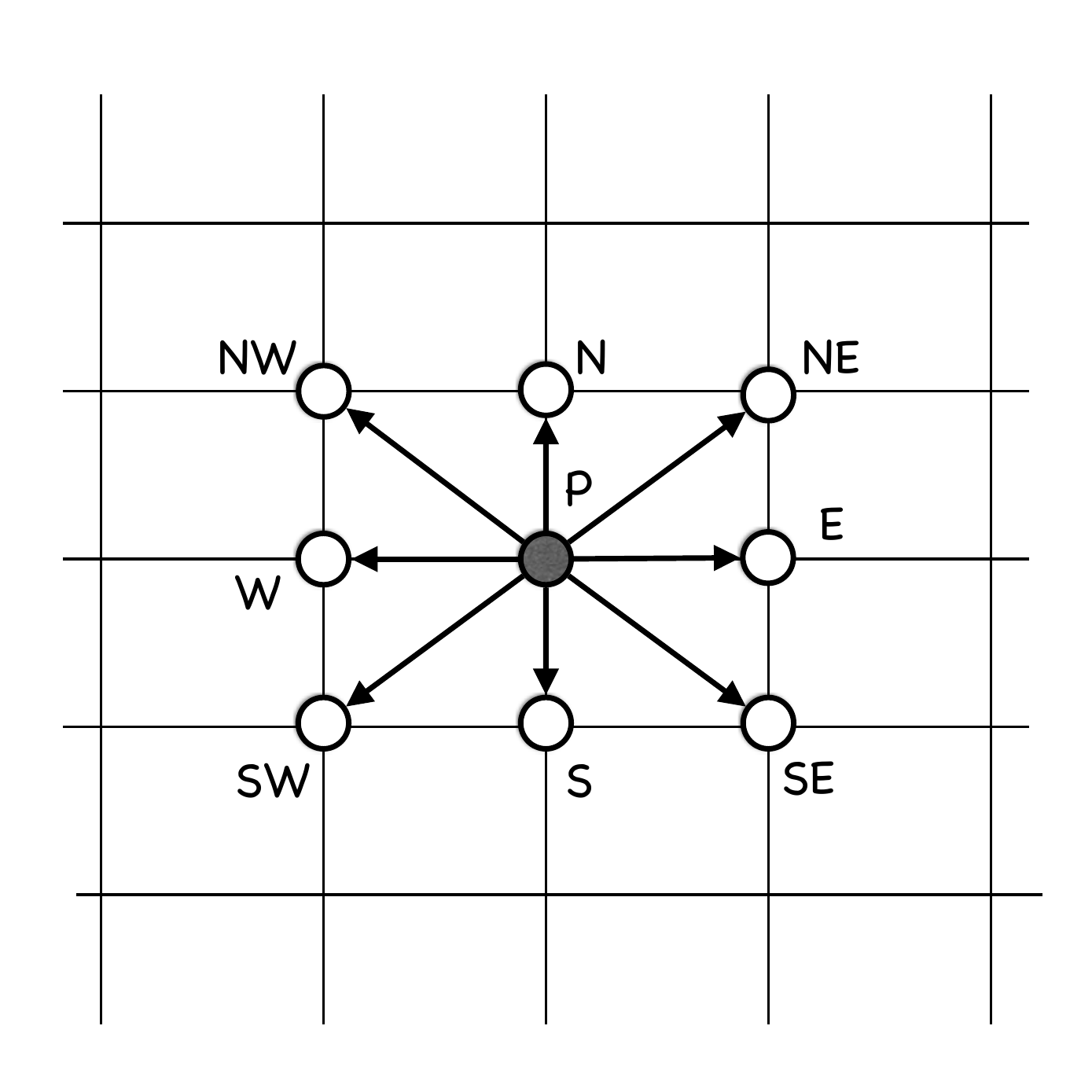}

\caption{The eight nearest neighbors of the mesh point $P$: $N$, $NE$, $E$, $SE$, $S$, $SW$, $W$, and $NW$.}
\end{center}
\label{fig:NN}       % Give a unique label
\end{figure}

Throughout the rest of the paper, we will assume that the attractor $A$ 
of the system under consideration is fixed and omit the subscript $A$
in the notation for the quasi-potential with respect to $A$: $U_A(\mx)\equiv U(\mx)$. 
We consider a regular rectangular mesh with mesh steps $h_1$ and $h_2$ in $x_1$ and $x_2$
directions respectively, and set $h:=\max\{h_1,h_2\}$. 
The mesh defines the sets of nearest neighbors for each mesh points. 
Every inner mesh point $P$ has eight nearest neighbors
surrounding it as shown in Fig. \ref{fig:NN}

Similar to as it is done in the OUM \cite{oum1,oum2,quasi}, the mesh points are divided into the following four categories.
\begin{itemize}
\item[0.] {\sf Unknown points:} the points where the solution $U$ has not been computed yet, and none of its nearest neighbors is {\sf Accepted} or {\sf Accepted Front}.
\item[1.] {\sf Considered:} the points { that} have {\sf Accepted Front} nearest neighbors. 
Tentative values of $U$ that might change as the algorithm proceeds, are available at them.
\item[2.] {\sf Accepted Front:} the points at which $U$ has been computed and no longer can be updated, and they have at least one {\sf Considered} nearest neighbor.
\item[3.] {\sf Accepted:} the points at which $U$ has been computed and no longer can be updated, 
and they have only {\sf Accepted} and/or {\sf Accepted Front} nearest neighbors.
\end{itemize}

The general outline of the OLIMs coincides with the one of the OUM.

{\bf Initialization}\\ 
Start with all mesh points being {\sf Unknown}. 
Compute tentative values { of $U$} at mesh points near the attractor $A$ and make them {\sf Considered}.

{\bf The main body}\\
{\bf while} $\{$ the boundary of the mesh is not reached $\}$ {\sf and} $\{$ the set of {\sf Considered} points is not empty $\}$ \\
\hspace*{5ex}{\bf 1:} Make the {\sf Considered} point $\mx$ with the smallest tentative value of $U$ {\sf Accepted Front}. \\
\hspace*{5ex}{\bf 2:} Make all  {\sf Accepted Front} nearest neighbors of $\mx$ that no longer have {\sf Considered} nearest neighbors  {\sf Accepted}.\\
\hspace*{5ex}{\bf 3:} Update all  {\sf Considered} points within the distance $Kh$ from $\mx$ using $\mx$ and maybe its nearest {\sf Accepted Front} neighbors.\\
\hspace*{5ex}{\bf 4:} Make all {\sf Unknown} nearest neighbors of $\mx$ {\sf Considered} and compute tentative values of $U$ at them using the 
{\sf Accepted Front} points lying within the distance $Kh$ from them.\\
{\bf end while}

In the rest of this section, we will elaborate the initialization and steps 3 and 4 of the while-cycle.

%%%%%%%%%%%%%%%%%%%%%%%%%%%
\subsection{Initialization}
\label{sec:init}
The way the initialization is performed is important for the accuracy of the computation of the  quasi-potential.
The two types of attractors in 2D, the asymptotically stable equilibrium point and the asymptotically stable limit cycle,
correspond to the boundary condition $U(A)=0$ given at an initial point or curve respectively.
The initialization procedure depends on whether the computation starts from the initial point  or  the initial curve.

\begin{figure*}
\begin{center}
  \includegraphics[width=0.7\textwidth]{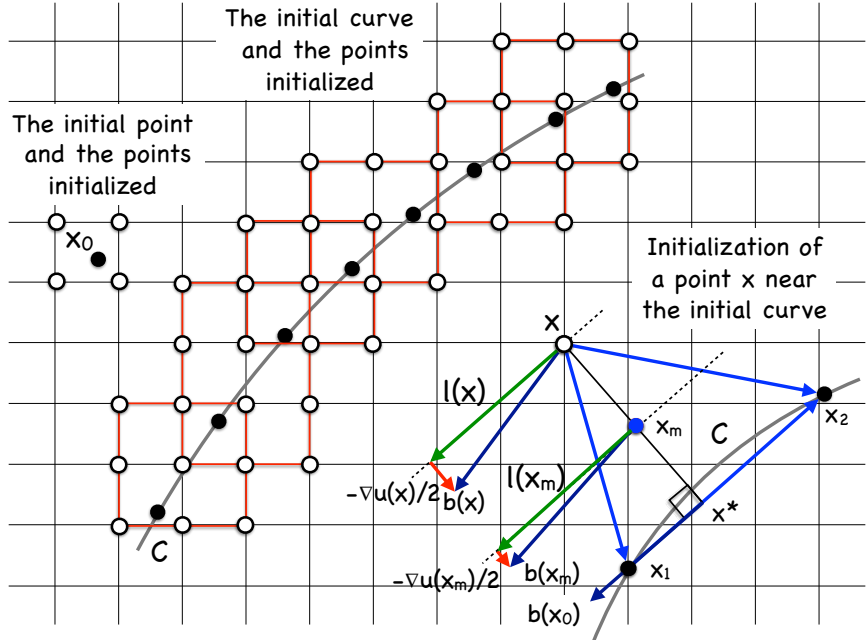}
\caption{Initialization of the mesh points near the initial point corresponding to an asymptotically 
stable equilibrium and a stable limit cycle of  $\dot{\mx}=\mb(\mx)$.}
\end{center}
\label{init_fig}       % Give a unique label
\end{figure*}
%

%%%
\subsubsection{Initialization from the initial point}
\label{sec:ipoint}
Let $\mx_0$ be an asymptotically stable equilibrium of $\dot{\mx}=\mb(\mx)$. 
Since  $\mb(\mx) = [b_1(\mx),b_2(\mx)]^T$ is twice continuously differentiable and $\mb(\mx_0) = \mathbf{0}$, we get
\begin{equation}
\label{bappr}
\mb(x) = \left[\begin{array}{cc}\partial_{x_1}b_1&\partial_{x_2}b_1\\\partial_{x_1}b_2&\partial_{x_2}b_2\end{array}\right](\mx-\mx_0) + O(\|\mx-\mx_0\|^2)=:
\mathbf{A}(\mx-\mx_0) + O(\|\mx-\mx_0\|^2)
\end{equation}
 in the neighborhood of $\mx_0$. 
We calculate the quasi-potential for the linear approximation $\mb(\mx)\approx \mathbf{A}(\mx-\mx_0)$
at the four mesh points
 surrounding the point $\mx_0$  (Fig. \ref{init_fig}, Left) using the analytic formula \cite{quasi}
 for the quasi-potential for a linear SDE with an asymptotically stable equilibrium. Let $A_{11}$, $A_{12}$, $A_{21}$, and $A_{22}$ 
 be the entries of the matrix $\mathbf{A}$ in Eq. \eqref{bappr}. Then the quasi-potential in the neighborhood of $\mx_0$ is approximated by
 \begin{align}
 U(\mx)& \approx(\mx-\mx_0)^T\left[\begin{array}{cc} \mathcal{A}&\mathcal{B}\\\mathcal{B}&\mathcal{C}\end{array}\right](\mx-\mx_0),
 \quad{\rm where}\notag \\
\mathcal{A} &= -(\alpha A_{11} + \beta A_{21}),\notag \\
\mathcal{B} &= -(\alpha A_{12} + \beta A_{22}),\label{lin_approx}\\
\quad \mathcal{C} &= -(\alpha A_{22} - \beta A_{12}),
\notag \\
\alpha& = \frac{(A_{11} + A_{22})^2}{(A_{11} + A_{22})^2 + (A_{21} - A_{12})^2},\quad 
\beta = \frac{ (A_{21} - A_{12})(A_{11} + A_{22})}{(A_{11} + A_{22})^2 + (A_{21} - A_{12})^2}.\notag
\end{align}
These four mesh points surrounding the point $\mx_0$ become {\sf Considered}.

%%%%%
\subsubsection{Initialization from the initial curve}
\label{sec:icurve}
The initialization of the mesh points in the neighborhood of the initial curve $C$
is based on the observation that the gradient of the quasi-potential vanishes at $C$ (see Fig. \ref{init_fig}, Center and Right). 
Assuming that $C$ is a smooth curve, there is a neighborhood $\mathcal{N}(C)$  of $C$ 
where the gradient of the quasi-potential  $\nabla U(\mx)$ is a smooth vector field. 
 Hence,  in $\mathcal{N}(C)$, the following decomposition into two smooth vector fields takes place: 
\begin{equation}
\label{bcurve1}
\mb(\mx) = -\frac{\nabla U(\mx)}{2} + \mathbf{l}(\mx),\quad  \nabla U(\mx)\cdot \mathbf{l}(\mx) = 0.
\end{equation}
One can easily deduce this decomposition using Eq. \eqref{HJ}. 
Since the rotational component $\mathbf{l}(\mx)$ is smooth, its direction in  $\mathcal{N}(C)$ changes continuously and 
 can be approximated by  the direction of $\mathbf{l}$ at the orthogonal projection ${\rm Proj}_C\mx$ of $\mx$ onto the curve $C$. 
Taking into account that $\nabla U(\mx)$ vanishes at $C$ and therefore
\begin{equation}
\label{bcurve2}
 \left.\mb(\mx)\right\vert_{C} = \left.\mathbf{l}(\mx)\right\vert_{C},
 \end{equation}
we approximate the gradient of the quasi-potential in $\mathcal{N}(C)$  by 
\begin{equation}
\label{bcurve3}
-\frac{\nabla U(\mx)}{2} = \mb(\mx) - \frac{\mb(\mx)\cdot\mb({\rm Proj}_C\mx)}{\mb({\rm Proj}_C\mx)\cdot \mb({\rm Proj}_C\mx)}\mb({\rm Proj}_C\mx).
\end{equation}
The second term in the right-hand side of Eq. \eqref{bcurve3} is the projection of $\mb(\mx)$ onto $\mb({\rm Proj}_C\mx)$.
The curve $C$ is provided by the user by a set of points $\{\mx_i\}_{i=0}^{M-1}$. It is assumed that $\{\mx_{i},\mx_{i+1}\}$, $i = 0,\ldots,M-2$,  
and $\{\mx_{M-1},\mx_0\}$ are consecutive  points along the curve.  
We define 
${\sf next}(i) = i + 1$ for $0\le i\le M-2$, and ${\sf next}(M-1) = 0$. 
For each pair of the consecutive points along the curve $(\mx_i,\mx_{{\sf next}(i)})$, 
we define the smallest rectangle  containing $\mx_i$ and $\mx_{{\sf next}(i)}$
with sides 
lying on the mesh lines 
and denote it by ${\sf Rectangle}(\mx_i,\mx_{{\sf next}(i)})$. 
A collection of such rectangles is shown  in Fig. \ref{init_fig}, Center,
by red contours.
We will initialize all mesh points lying in the union of the rectangles
\begin{equation}
\label{icurve1}
{\sf Neib}(C): = \bigcup_{0\le i < M}{\sf Rectangle}(\mx_i,\mx_{{\sf next}(i)}).
\end{equation}

Each mesh point $\mx$ located within ${\sf Neib}(C)$ is initialized using Simpson's quadrature rule as follows. 
{ First, observe that for any pair of points $\mx_i$ and $\mx_f$ we have
\begin{equation}
\label{icur_aux}
U(\mx_f) = U(\mx_i) + \int_0^1 \nabla U\left(\mx_i + t(\mx_f - \mx_i)\right)\cdot(\mx_f -\mx_i)dt.
\end{equation}
Let us take $\mx_f\equiv \mx$, the point to be initialized, and $\mx_i = {\rm Proj}_C\mx$, the projection of $\mx$ onto the curve $C$.
Note that $U({\rm Proj}_C\mx) = 0$ and $\nabla U({\rm Proj}_C\mx)=0$. An approximation $\mx^{\ast}$
to ${\rm Proj}_C\mx$ is constructed as follows.}
Let $\mx_1$ be the closest point of the curve $C$ to the point $\mx$, and the point $\mx_2$ be the neighboring point of $\mx_1$ along the curve $C$
such that the angle  between the vectors $\mx-\mx_1$ and $\mx_2-\mx_1$ does not exceed $\pi/2$ (Fig. \ref{init_fig}, Right).
The projection ${\rm Proj}_C(\mx)$ is approximated by the orthogonal projection of $\mx$ onto the line 
passing through $\mx_1$ and $\mx_2$: 
{
$$
\mx^{\ast} =\mx_1 +((\mx - \mx_1)\cdot \hat{\tau})\hat{\tau},\quad{\rm where}\quad \hat{\tau}:=\frac{\mx_2 - \mx _1}{ \|\mx_2 - \mx _1\|}.
$$
To use Simpson's rule for computing the integral in Eq. \eqref{icur_aux},
we need to evaluate $U(\mx^{\ast})$, $\nabla U(\mx^{\ast})$,  $\nabla U(\mx)$, and $\nabla U(\mx_m)$, 
where $\mx_m$  is the midpoint between $\mx$ and $\mx^{\ast}$. 
Since $U = 0$  and $\nabla U=0$ on $C$, $U(\mx^{\ast})\approx 0$ and $\nabla U(\mx^{\ast})\approx 0$.
Next, we note that  on the line segment $[\mx^{\ast},\mx]$, $\nabla U$ is nearly parallel to $(\mx - \mx^{\ast})$, 
hence $\nabla U \cdot  (\mx - \mx^{\ast}) \approx \|\nabla U \|\| \mx - \mx^{\ast}\|$.
Finally,  we approximate $\|\nabla U(\mx_m) \|$ and $\| \nabla U(\mx)\| $ using Eq. \eqref{bcurve3} and get}
\begin{equation}
\label{icurve2}
U(\mx) \approx \|\mx - \mx^{\ast}\|\left(4\|\nabla U(\mx_m)/2\| + \|\nabla U(\mx)/2\|\right)/3.
\end{equation}
All initialized mesh points become {\sf Considered}.

%%%%%%%%%%%%%%%%%%%%%%%%%%%%%%%%%%%%%%%%%%%%
\subsection{The one-point updates and the triangle updates} 
The OUM \cite{oum1,oum2,quasi} involves two types of updates: the {one-point update} and the {triangle update}.
The OUM {one-point update} applied to  \eqref{HJ} reads
\begin{align}
{\sf R}_{1pt}(\mx_0,\mx) &= U(\mx_0)  +\|\mx - \mx_0\|\|\mb(\mx)\| - (\mx - \mx_0)\cdot \mb(\mx),\label{1ptR}\\
U(\mx) &= \min\{ {\sf R}_{1pt}(\mx_0,\mx) ,U(\mx)\},\notag
\end{align}
where $\mx$ is the mesh point being updated and $\mx_0$ is an {\sf Accepted Front} point within the distance $Kh$ from $\mx_0$.
Eq. \eqref{1ptR} is merely  the integrand in Eq. \eqref{Qpot_def} integrated along the line segment from $\mx_0$ to $\mx$
using the right-hand rectangle quadrature rule. 
The OUM {triangle update} is done by the upwind finite difference scheme applied to PDE \eqref{HJ}. 
Its details are worked out in Section \ref{sec:equiv} below.

In the proposed OLIMs, the finite difference scheme as well as PDE \eqref{HJ} are completely abandoned. 
The {one-point update} and  the {triangle update} are done as follows. Let $\mathcal{Q}$ 
be a basic quadrature rule, i.e., right-hand rectangle, midpoint, trapezoid, or Simpson's.
The  OLIM {one-point update} results from the application of $\mathcal{Q}$ for the
integration of the integrand of Eq. \eqref{Qpot_def}
along the line segment $[\mx_0,\mx]$:
\begin{equation}
\label{1ptOLIM}
{\sf Q}_{1pt}(\mx_0,\mx) = U(\mx_0) + \mathcal{Q}(\mx_0,\mx), \quad 
U(\mx) = \min\{{\sf Q}_{1pt}(\mx_0,\mx),U(\mx)\}.
\end{equation}
{I.e., we compute the value ${\sf Q}_{1pt}(\mx_0,\mx)$ 
and replace the current value of $U$ at $\mx$ with ${\sf Q}_{1pt}(\mx_0,\mx)$ if and only if it is smaller
than the current value.}
The OLIM {triangle update} at the mesh point $\mx$ from the triangle $(\mx_1,\mx_0,\mx)$, where the points $\mx_1$ and $\mx_0$ are assumed to be nearest neighbors,
results from solving the minimization problem
\begin{align}
{\sf Q}_{\Delta}(\mx_1,\mx_0,\mx) &= \min_{s\in[0,1]} \left\{ sU(\mx_0)  + (1-s)U(\mx_1) + \mathcal{Q}(s\mx_0+(1-s)\mx_1,\mx)\right\}, \label{2ptOLIM}\\
U(\mx) &= \min\{{\sf Q}_{\Delta}(\mx_1,\mx_0,\mx),U(\mx)\},\notag
\end{align}
where the vector field $\mb$ is approximated with a linear field within the triangle $(\mx_1,\mx_0,\mx)$ 
for the purpose of the application of the quadrature rule. 
The minimization problem \eqref{2ptOLIM} is solved by taking the derivative of the function to be minimized, setting it to zero, i.e.,
$$
 g(s):= \frac{d}{ds}\left(sU(\mx_0)  + (1-s)U(\mx_1) + \mathcal{Q}(s\mx_0+(1-s)\mx_1,\mx)\right) = 0,
$$
and then trying to find a root of $g$ in the interval $[0,1]$. 
For this purpose, we use the hybrid secant/bisection method \cite{wilkinson,stewart}.

Depending on the quadrature rule used, we denote the methods OLIM-R, OLIM-MID, OLIM-TR, and OLIM-SIM (see Section \ref{sec:intro}).
The details of the {triangle update} for each OLIM are worked out in Appendix B.

%%%%%%%%%%%%%%%%%%%%%%%%%%%%%%%%%%%%%%%%%

\subsection{Equivalence of update rules in the OUM and OLIM-R}
\label{sec:equiv}
The {one-point update} in  OLIM-R and the OUM are done according to Eq. \eqref{1ptR}.
In this section, we show that the {triangle updates} in OLIM-R and the OUM, if successful, give identical results.

The {triangle update} in the OUM \cite{quasi} is done as it is proposed in \cite{oum2}.
The quasi-potential $U(\mx)$ is assumed to be linear within the triangle $(\mx_1,\mx_0,\mx)$
and have {\sf Accepted} values $u_1$ and $u_0$ at $\mx_1$ and $\mx_0$  respectively (see Fig. \ref{fig:equiv}). 
Its value at $\mx$ is to be found.
The finite-difference scheme proposed in \cite{oum1,oum2} applied to Eq. \eqref{HJ}
is derived from the observation that 
\begin{equation}
\label{eq1}
\left[\begin{array}{c}u - u_0\\u-u_1\end{array}\right] = \left[\begin{array}{c}(\mx - \mx_0)^T\\(\mx-\mx_1)^T\end{array}\right]\nabla U,
\end{equation}
where $u = U(\mx)$, $u_0=U(\mx_0)$, and $u_1 = U(\mx_1)$ are the values of $U$ at the points $\mx$, $\mx_0$, and $\mx_1$ respectively.
The vector field $\mb$ within the triangle $(\mx_1,\mx_0,\mx)$ is approximated by the constant field $\mb\equiv\mb(\mx)$ ($\mb$ evaluated at $\mx$).
Throughout the rest of this section, we will omit the argument of $\mb$. 
Denoting the matrix { on} the right-hand side of Eq. \eqref{eq1} by $P$ and plugging $\nabla U$ into Eq. \eqref{HJ}, we obtain 
the following quadratic equation for $u$:
\begin{equation}
\label{dfeq}
\left[u - u_0, u-u_1\right] P^{-T}P^{-1}\left[\begin{array}{c}u - u_0\\u-u_1\end{array}\right] + 2\mb^TP^{-1}\left[\begin{array}{c}u - u_0\\u-u_1\end{array}\right] = 0.
\end{equation}
If Eq. \eqref{dfeq} has a solution $u$ satisfying $u\ge\min\{u_0,u_1\}$, it is subjected to the consistency test that the characteristic
passing through the point $\mx$ crosses the segment $[\mx_1,\mx_0]$. The direction of the characteristic is given by $\mb(\mx)+\nabla U$ \cite{quasi}.
The condition that it crosses the segment $[\mx_1,\mx_0]$ is equivalent to the fact that  
$\mb(\mx)+\nabla U = c_0(\mx-\mx_0) + c_1(\mx-\mx_1)$ is a linear combination with nonnegative coefficients $c_0$, $c_1$, i.e., the vector 
\begin{equation}
\label{eq2}
\left[\begin{array}{c}c_0\\c_1\end{array}\right]=P^{-T}\left(\mb(\mx) + P^{-1} \left[\begin{array}{c}u - u_0\\u-u_1\end{array}\right]\right)
\end{equation}
has nonnegative entries. If $u$ has  passed the consistency test then $U(\mx) = \min\{u,U(\mx)\}$.

OLIM-R performs the {triangle update} by solving the following minimization problem
\begin{align}
&u = \min_{s\in[0,1]}\left[su_0 + (1-s)u_1 + \|\mb\|\|\mx - \mx_s\| - \mb\cdot(\mx - \mx_s)\right],
\label{r_olim}\\
&\text{where}~~\mx_s = s\mx_0 + (1-s)\mx_1.\notag
\end{align}
Like Eq. \eqref{dfeq}, Eq. \eqref{r_olim} is set up under the assumption that  the function $U(\mx)$ is linear within the triangle $(\mx_1,\mx_0,\mx)$.

%%%
The solution by the finite difference scheme \eqref{dfeq} and the solution of the minimization \eqref{r_olim} are equivalent in the sense 
specified by the following theorem.
\begin{theorem}
\label{thm_equiv}
Let $u$ be a solution of Eq. \eqref{dfeq} such that the coefficients $c_0$ and $c_1$ in Eq. \eqref{eq2} are positive, i.e., 
the characteristic passing through the point $\mx$ crosses the open interval $(\mx_1,\mx_0)$. Then $u$ is also the solution of the minimization problem 
\eqref{r_olim} and the corresponding minimizer $s^{\ast}$ satisfies $0<s^{\ast}<1$. 

Conversely, let $u$ be the solution of the minimization problem 
\eqref{r_olim} with the corresponding minimizer $s^{\ast}$ satisfying $0<s^{\ast}<1$. Then $u$ is also a solution of Eq. \eqref{dfeq}
such that the characteristic passing through the point $\mx$ crosses the open interval $(\mx_1,\mx_0)$.
\end{theorem}

The proof of Theorem \ref{thm_equiv} is given in Appendix C.

%%%%%%%%%%%%%%%%%%%%%%%%%%%

\subsection{Reducing CPU time: { a hierarchical update strategy}}
\label{sec:redCPU}
\begin{figure*}
\begin{center}
(a)\includegraphics[width=0.6\textwidth]{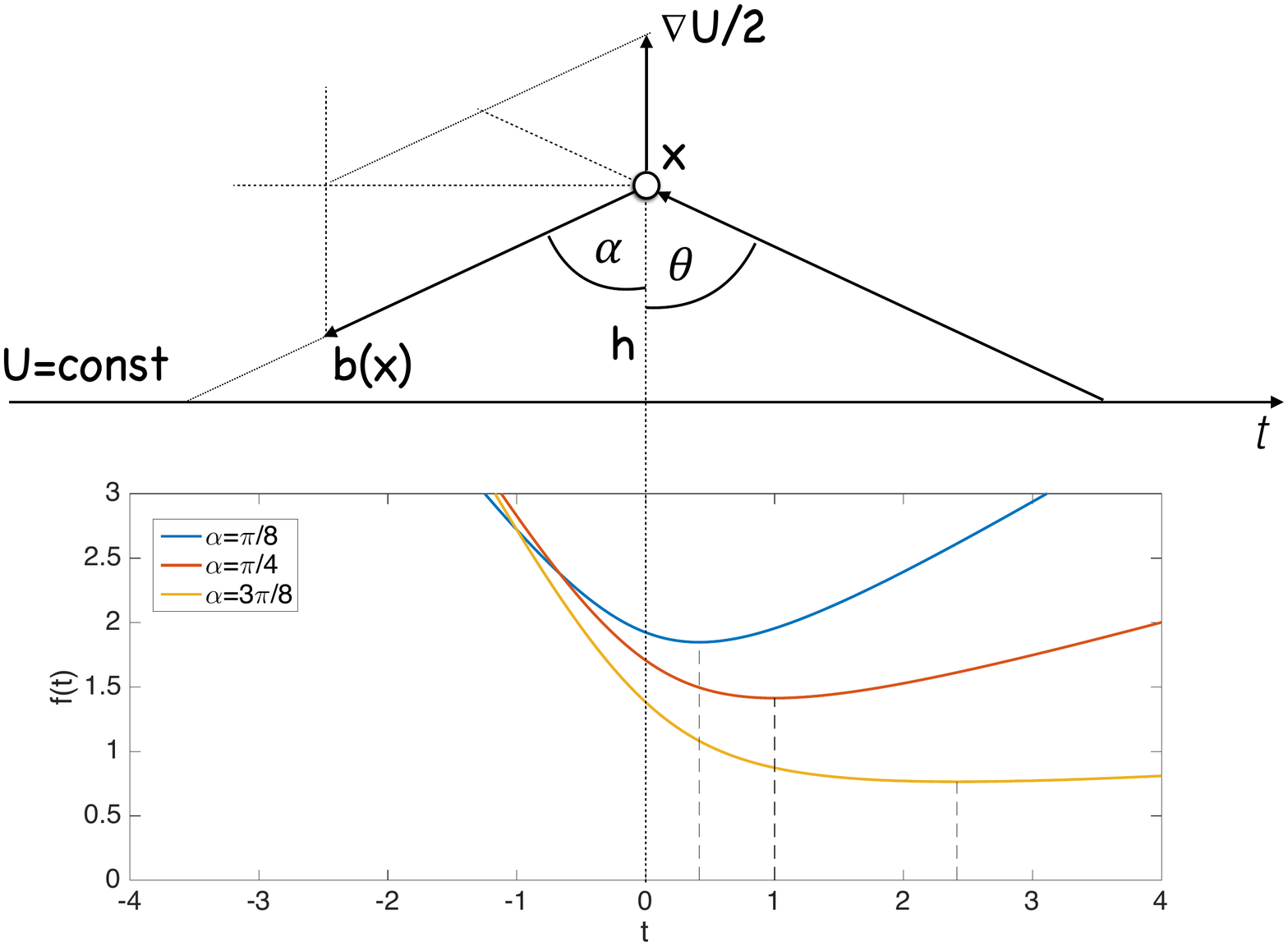}\\
(b)\includegraphics[width=0.6\textwidth]{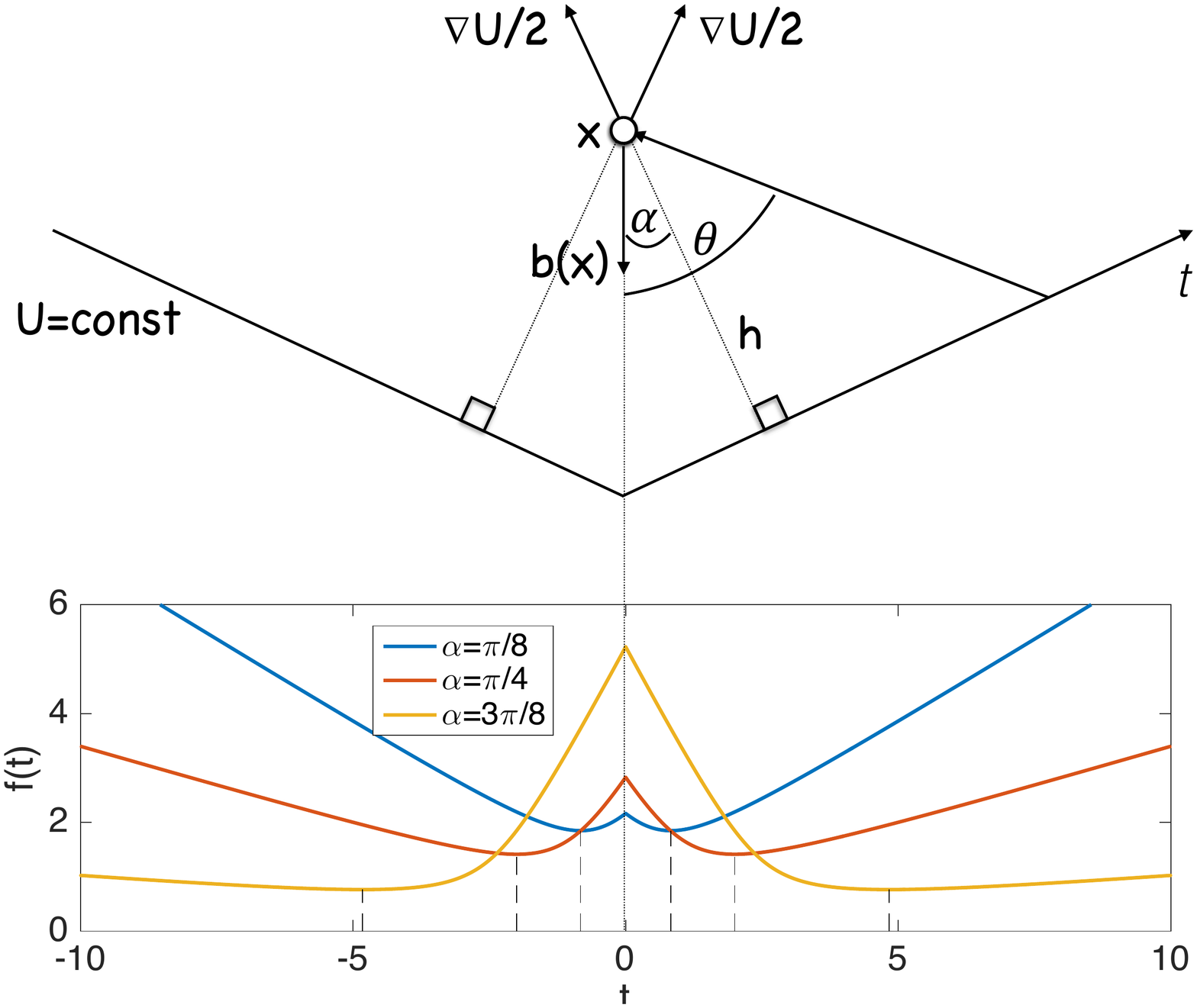}
\caption{ An illustration for Section \ref{sec:redCPU}. (a): The case with no kink. 
Top: The $t$-axis coincides with a level set of $U$. $\mx$ is a point to be updated.
Bottom: The graphs of the update function $f(t)$  (Eq. \eqref{update_fun}) for three different values between $\mb$ and $-\nabla U$ for $|\mb| = h = 1$. 
These functions have unique minima at $t = h\tan(\alpha)$ marked by the dashed lines.
(b): The case with a kink. Top: $-\infty<t<\infty$ is the arclength parameter along a level set of $U$ with a kink. $\mx$ is a point on the kink curve to be updated.
$\mb(x)$ for $x$ on the kink curve is directed toward the kink.
Bottom: The graphs of the update function $f(t)$ for three different values between $\mb$ and $-\nabla U$ for $|\mb| = h = 1$. 
These functions have two local minima at $t = \pm 2h\tan(\alpha)$ marked by the dashed lines.
}

\end{center}
\label{fig:rat}
\end{figure*}

Clearly, the {one-point updates} require significantly fewer floating point operations  { than} the {triangle updates} in both the OUM
and the OLIMs. Therefore, we want to develop of a strategy for reducing the number of calls for the {triangle update}. 
Suppose that the mesh is fine enough so that $\mb$ is approximately constant within the update radius $Kh$, 
and level sets of $U$ can be approximated by line segments  within the distance $Kh$ from the mesh point $\mx$ to be updated.

Motivated by this, we consider the following idealized situation. Let $\mb$ be constant
and the level set $U = U_0$ be a straight line.
Let $\mx$ be a point at the distance $h$ from this level set as shown in Fig. \ref{fig:rat} (a,Top).
Then 
$$
U(\mx) = \min_{t\in \mr}[U_0 +  f(t)],
$$
where $t$ is the axis along the level set, and $f(t)$ is the update function
\begin{equation}
\label{update_fun}
f(t): = |\mb|\frac{h}{\cos(\theta(t))} - |\mb|\frac{h}{\cos(\theta(t))}\cos(\pi - (\alpha + \theta)).
\end{equation}
The angles $\alpha\in[0,\pi/2)$ and $\theta\in(-\pi/2,\pi/2)$ are defined as shown in Fig. \ref{fig:rat} (a, Top).
Differentiating $f$ with respect to $\theta$,
setting its derivative to zero, and then returning to the variable $t$, we find that $f$ has a single minimum on $-\infty < t < \infty$ 
achieved at $t=h\tan(\alpha)$ (see Fig. \ref{fig:rat} (a, Bottom)).

Now we address the case where the quasi-potential $U$ is not differentiable along 
some curve, and hence its level sets have kinks.
Typically, the vector field at the kink curve is directed approximately toward the kink. 
For example, see the quasi-potential computed  in Ref. \cite{quasi} for the Maier-Stein model \cite{maier-stein}. 
To account for such a case, we consider the following idealized situation.
Let $\mb$ be constant, and the level set $U = U_0$ consists of two line segments.
Let $t$ be the arclength parameter along the level set, $t=0$ at the kink, 
and the angles $\alpha $ and $\theta$ be defined as shown in Fig. \ref{fig:rat} (b, Top).
Let the point $\mx$  lie on the kink line. 
The update function $f(t)$ for this case is obtained from the one in Eq. \eqref{update_fun} using an appropriate reflection. 
Its graphs for three values of $\alpha$ are shown  in Fig.  \ref{fig:rat} (b, Bottom).
It has two local minima  at $t = \pm 2h\tan(\alpha)$.

These considerations suggest that the update values of $U$ at a {\sf Considered} point $\mx$ lying near the {\sf Accepted Front} 
as  functions of an arclength parameter along the level set of $U$ approximating the {\sf Accepted Front}
will behave similar to the graphs in Fig. \ref{fig:rat}(a, Bottom, or b, Bottom), i.e., have a single local minimum,
or two equal  local minima.
Therefore, we have rationalized the update procedures as follows.
When step 4  of the main body of the algorithm is performed, i.e., an {\sf Unknown} nearest neighbor $\mx$ 
of the new {\sf Accepted Front} point $\mx^{\ast}$ becomes {\sf Considered},
we compute tentative values at $\mx$ from every {\sf Accepted Front} point at distance at most $Kh$ from $\mx$ using only  the cheap {one-point update}.
Suppose the minimal tentative value at $\mx$ by  the {one-point update} has been computed from the {\sf Accepted Front} point $\mx_0$,
i.e., 
$$
\mx_0 =  \arg\min_{\mathbf{y}\in{\sf Accepted\medspace Front},\medspace\|\mx-\mathbf{y}\|\le Kh} {\sf Q}_{1pt}(\mathbf{y},\mx).
$$
Only then the {triangle updates} ${\sf Q}_{\Delta}(\mx_1,\mx_0,\mx)$ are called for all {\sf Accepted Front} nearest neighbors $\mx_1$ of $\mx_0$ (if any).
Then 
$$
U(\mx) = \min\{ \min_{\mathbf{\mx_1}\in{\sf Accepted\medspace Front}, \medspace\mathbf{\mx_1} \in \mathcal{N}(\mx_0)} {\sf Q}_{\Delta}(\mx_1,\mx_0,\mx), U(\mx)\}.
$$
The symbol $\mathcal{N}(\mx_0)$ denotes the set of the nearest neighbors of $\mx_0$ (eight nearest neighbors for every inner point as shown in Fig. \ref{fig:NN}).
{We will refer to this update strategy as the \emph{hierarchical update strategy}.}

Our numerical experiments in  Section \ref{sec:tests} below indicate that this strategy 
reduces CPU time in OLIM-R in comparison with the OUM by the factor of approximately 4, 
{ while the numerical errors in them either coincide, or the ones in OLIM-R exceed 
the ones in the OUM by less than 1\%.}

%%%%%%%%%%%%%%%%%%%%%%%%%%%%%%%%%%%
%
%    T E S T S
%
%%%%%%%%%%%%%%%%%%%%%%%%%%%%%%%%%%%

\section{Numerical tests}
\label{sec:tests}
We compare performances of the OLIMs (OLIM-R, OLIM-MID, OLIM-TR, and OLIM-SIM) and the original OUM-based \cite{quasi} 
quasi-potential solver on two examples for which analytic formulas for the quasi-potential are available. These two examples are quite challenging 
from the computational point of view due to the large rotational components of $\mb$ in comparison with $\nabla U$ and 
large curvatures of their MAPs (Minimum Action Paths). 
Furthermore, in the second example, the quasi-potential grows as a fourth degree polynomial.

%For brevity, we will refer to the original OUM-based method \cite{quasi} as the OUM.

\begin{figure*}
\begin{center}

(a)\includegraphics[width = 0.4\textwidth]{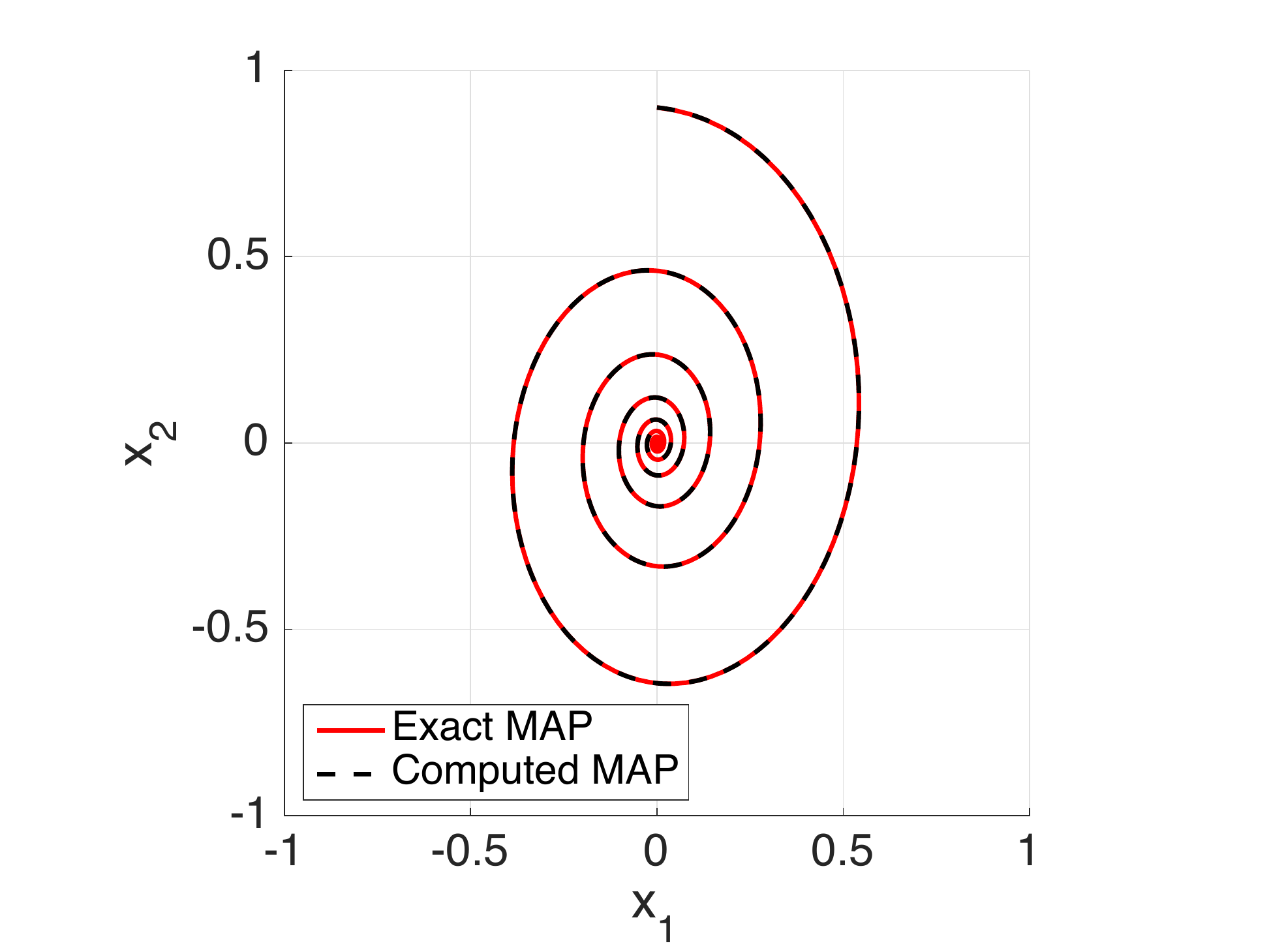}
(b)\includegraphics[width = 0.4\textwidth]{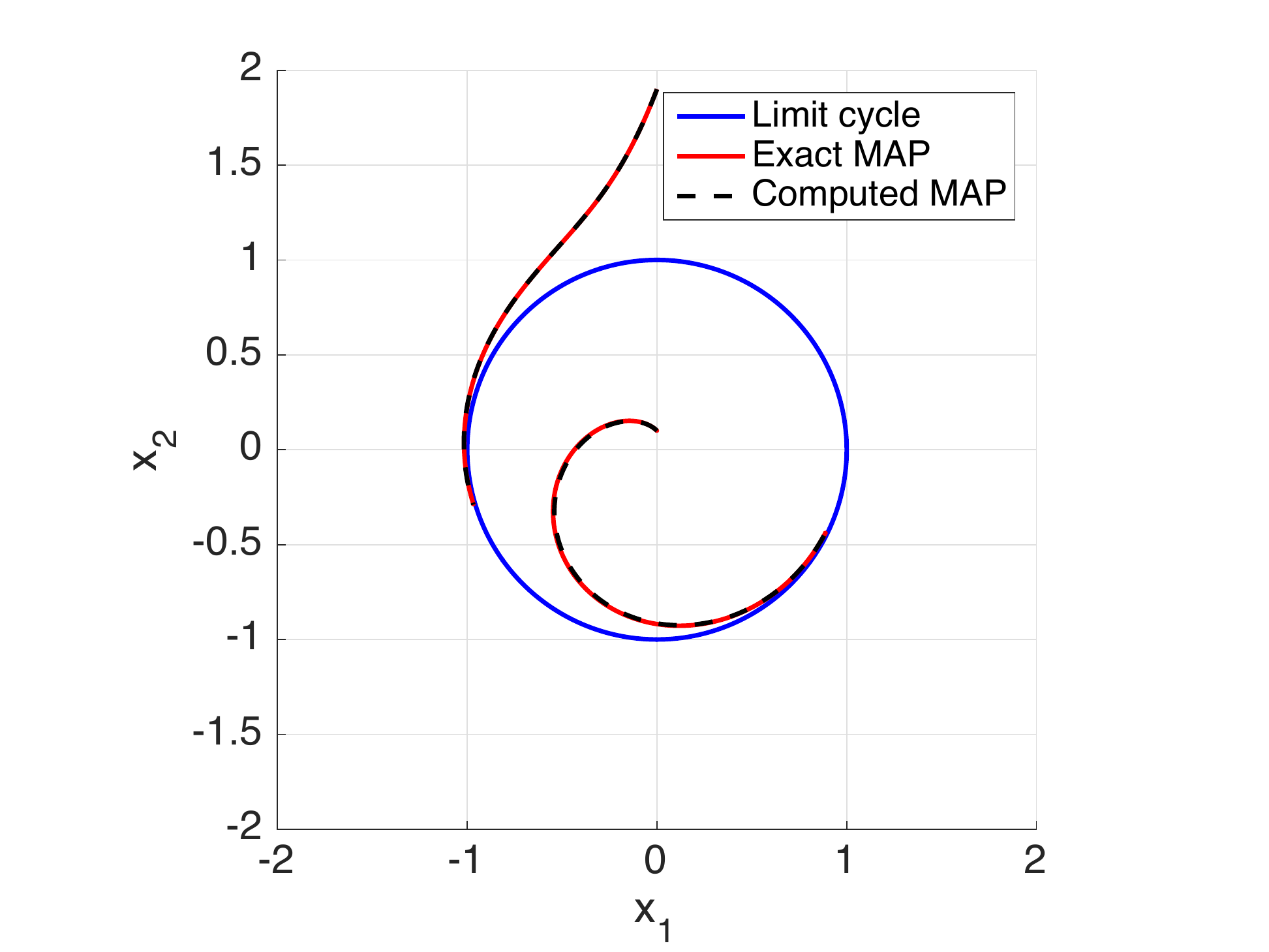}

\caption{
(a) The exact and the computed Minimum Action Paths (MAPs) for SDE \eqref{linSDE} arriving at the point $(0,0.9)$. 
(b) The exact and the computed Minimum Action Paths (MAPs) for SDE \eqref{linSDE} arriving at the points $(0,1.9)$ and $(0,0.1)$. 
In both cases, 
the computed MAPs were obtained by integrating the 
path $\dot{\phi} = -\left(\mb(\phi) + \nabla U(\phi)\right)$ \cite{quasi} using the 4-stage 4-th order Runge-Kutta method. $\nabla U(\mx)$ 
was found  by finite differences from the computed quasi-potential $U$ on the $1024\times 1024$ mesh with $K=20$. 
}
\end{center}
\label{fig:MAP}
\end{figure*}

%%%%%%   Linear Example
%\begin{example}
%\label{ex:lin}
\textbf{A linear SDE.}
For the linear SDE \cite{quasi}
\begin{align}
dx_1 &= (-2x_1 - ax_2)dt + \sqrt{\epsilon}dw_1,\notag\\
dx_2 &= (2ax_1 - x_2)dt +\sqrt{\epsilon}dw_2, \label{linSDE}
\end{align}
the quasi-potential with respect to the origin, the asymptotically stable equilibrium of the corresponding deterministic system, 
is the quadratic function 
\begin{equation}
\label{qlin}
U(x_1,x_2) = 2x_1^2 + x_2^2.
\end{equation}
The parameter $a$ is the quotient of the magnitudes of the rotational and the potential components ( $\mathbf{l}$ and $-\tfrac{1}{2}\nabla U$ respectively) 
of the vector field $\mb$:
\begin{equation}
\label{lin_bdec}
\mb(\mx) = - \tfrac{1}{2}\nabla U(\mx) + \mathbf{l}(\mx) = -\left[\begin{array}{c}2x_1\\x_2\end{array}\right] +
a\left[\begin{array}{c}-x_2\\2x_1\end{array}\right].
\end{equation}
We set $a=10$, i.e., 
the rotational component exceeds the potential component in magnitude by the factor of 10 
that makes the computation of the quasi-potential challenging.
The computational domain for this example is the square $[-1\le x_1\le 1]\times[-1\le x_2\le 1]$.
The Minimum Action Path for SDE \eqref{linSDE} arriving at the point $(0,0.9)$ is shown in Fig. \ref{fig:MAP}(a).
%\end{example}

%%%%% Circle Example
%\begin{example}
%\label{ex:cir}
\textbf{{ An} SDE with a limit cycle.}
The unit circle $C:=\{(x_1,x_2)\in\mr^2~|~ x_1^2 + x_2^2 = 1\}$ is the asymptotically stable limit cycle
of the deterministic system corresponding to the SDE \cite{hurewicz,quasi}
\begin{align}
dx_1 &= (x_2 + x_1(1-x_1^2 - x_2^2))dt + \sqrt{\epsilon}dw_1,\notag\\
dx_2 &= (-x_1+x_2(1-x_1^2 - x_2^2))dt +\sqrt{\epsilon}dw_2. \label{cirSDE}
\end{align}
The quasi-potential with respect to the unit circle is given by the quartic function\footnotemark[1]
\begin{equation}
\label{qlin}
U(x_1,x_2) = \tfrac{1}{2}(x_1^2 +x_2^2 - 1)^2.
\end{equation}
\footnotetext[1]{There is an error in Eq. (89) in \cite{quasi}. It should be $U=\tfrac{1}{2}(r^2-1)^2$.}
The decomposition of the vector field into the potential and the rotational components is
\begin{equation}
\label{lin_bdec}
\mb(\mx) = - \tfrac{1}{2}\nabla u(\mx) + \mathbf{l}(\mx) = \left[\begin{array}{c}x_1(1-x_1^2 - x_2^2)\\x_2(1-x_1^2 - x_2^2)\end{array}\right] +
\left[\begin{array}{c}x_2\\-x_1\end{array}\right].
\end{equation}
We set the computational domain for this example to be the square $[-2\le x_1\le 2]\times[-2\le x_2\le 2]$.
The Minimum Action Paths for SDE \eqref{linSDE} arriving at the points $(0,1.9)$ and $(0,0.1)$ are shown in Fig. \ref{fig:MAP}(b).
%\end{example}

We test the OLIMs and the OUM on SDEs \eqref{linSDE} and \eqref{cirSDE}.
The mesh sizes are $N\times N$ where $N=2^p$, $p = 7,8,\ldots,12$.
The update factor $K$ varies from $K=1$ to $K = 50$.
We have measured the maximum absolute error, the RMS error, and the CPU time 
for each mesh size $N$ and for each value of $K$.

%%%%

\subsection{Dependence of the accuracy on the update factor $K$}
\label{sec:Ktest}
The choice of the optimal value of the update factor $K$ is a subtle issue. 
On one hand, if the Minimum Action Paths (MAPs) for a considered SDE would be straight lines, 
large $K$ would enable all mesh points to be updated from the  correct triangles and hence enhance the accuracy.
However, typically, MAPs are not straight lines. Furthermore, 
large $K$ leads to numerical integration by simple quadrature rules along long linear segments and hence increases the integration error. 
Finally, large $K$ increases the CPU time. 
As a result, the optimal $K$ should be not too large and not too small.

\begin{figure*}
\begin{center}

(a)\includegraphics[width = 0.3\textwidth]{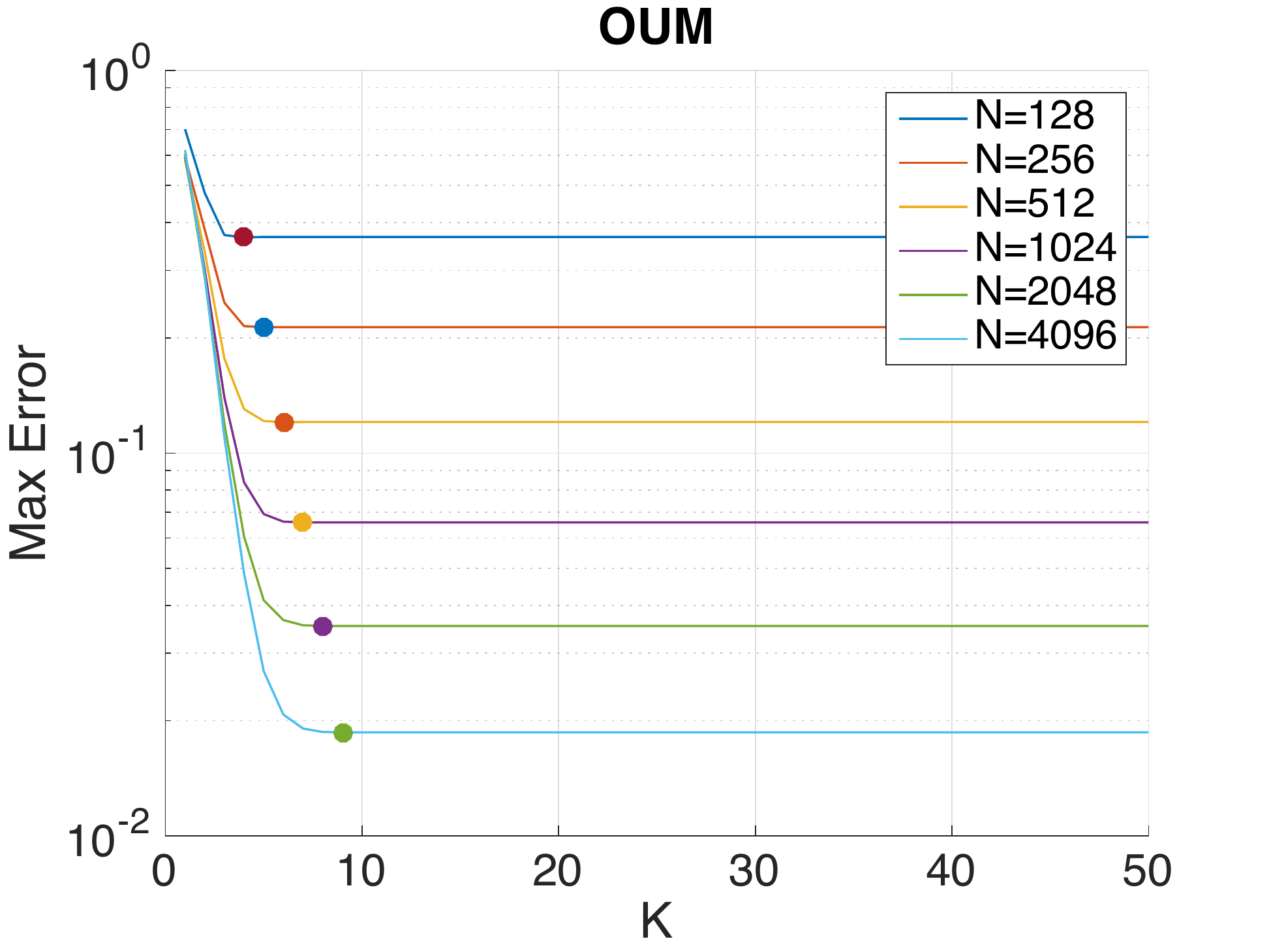}
(b)\includegraphics[width = 0.3\textwidth]{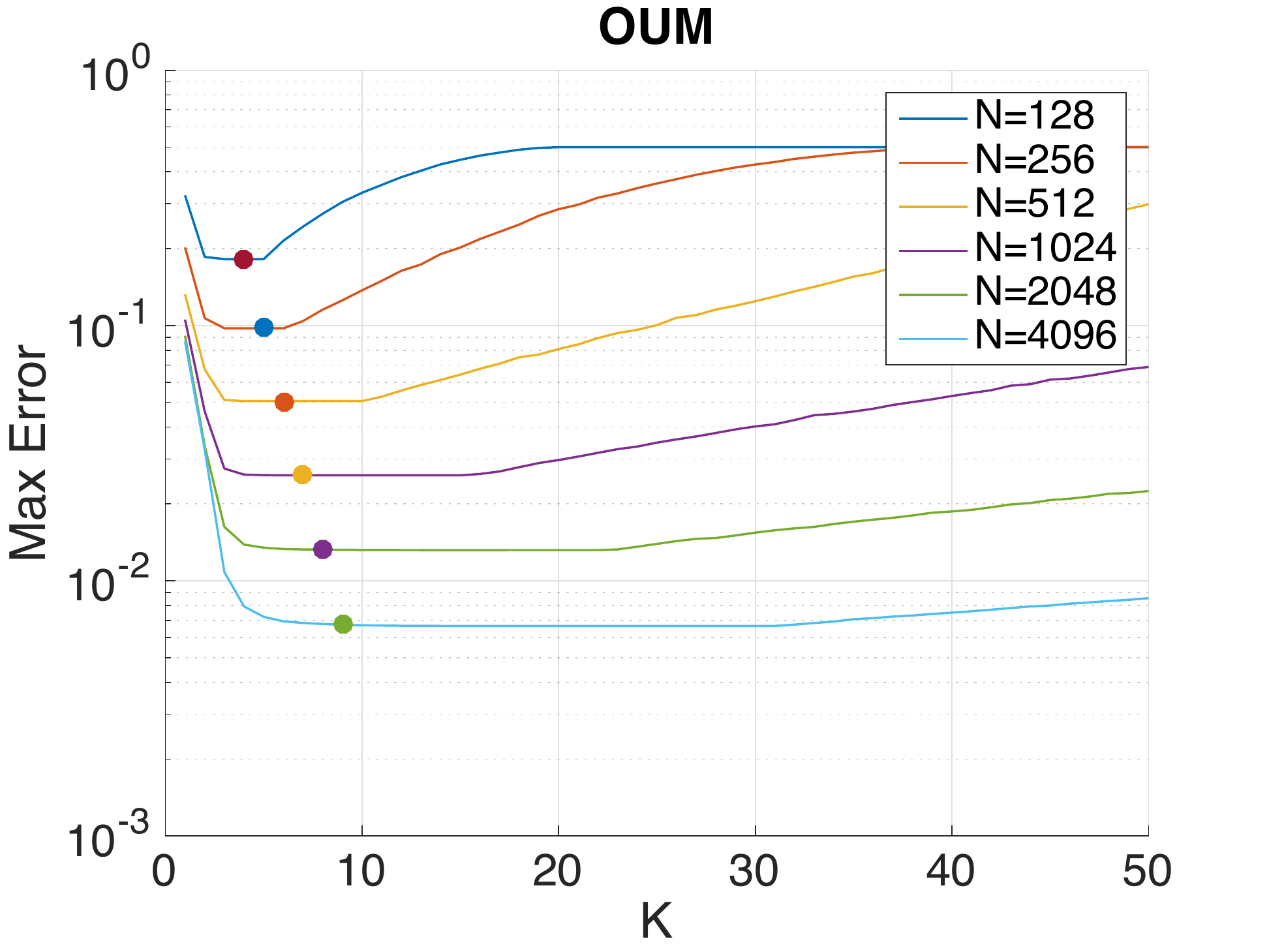}

(c)\includegraphics[width = 0.3\textwidth]{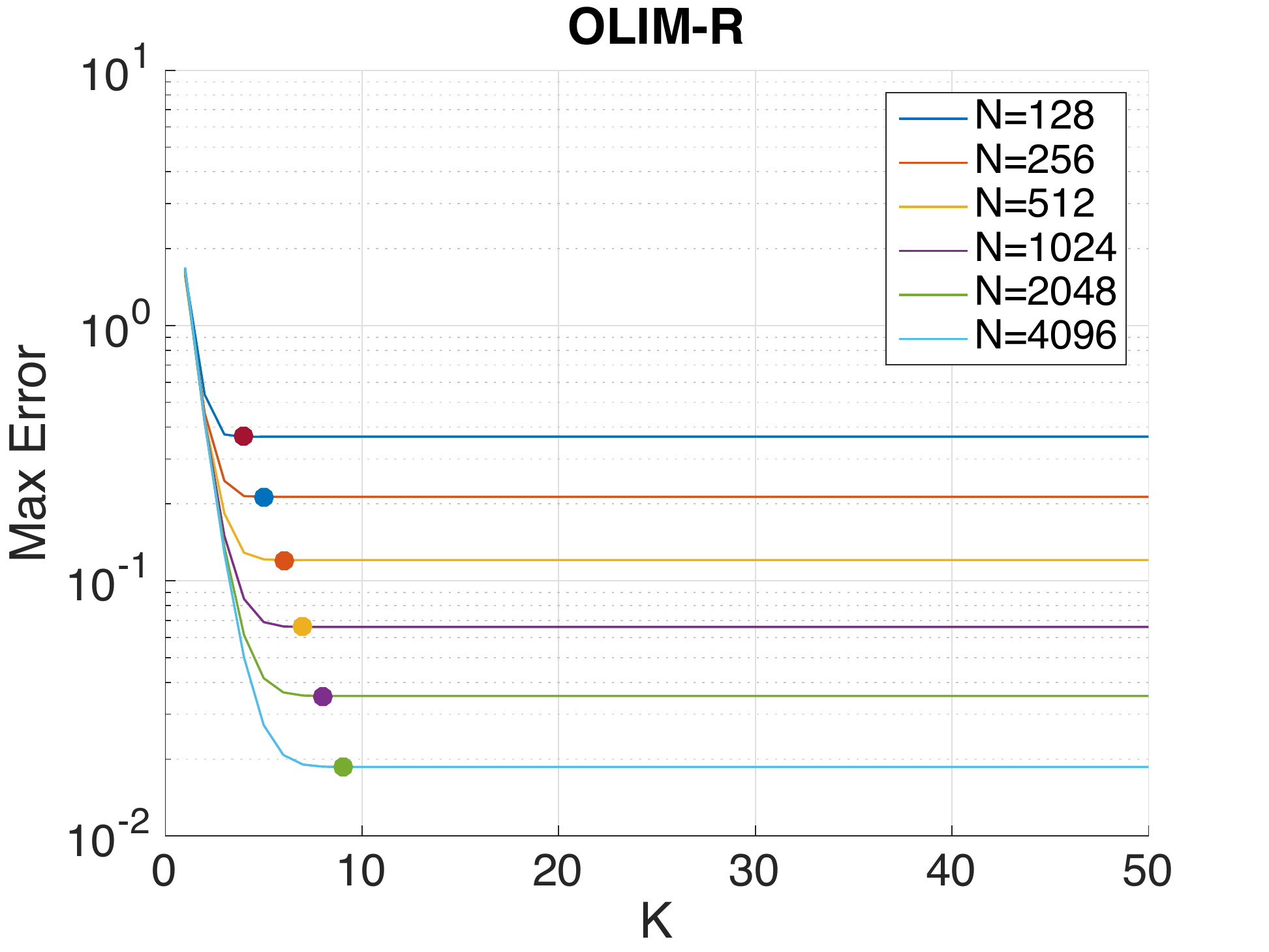}
(d)\includegraphics[width = 0.3\textwidth]{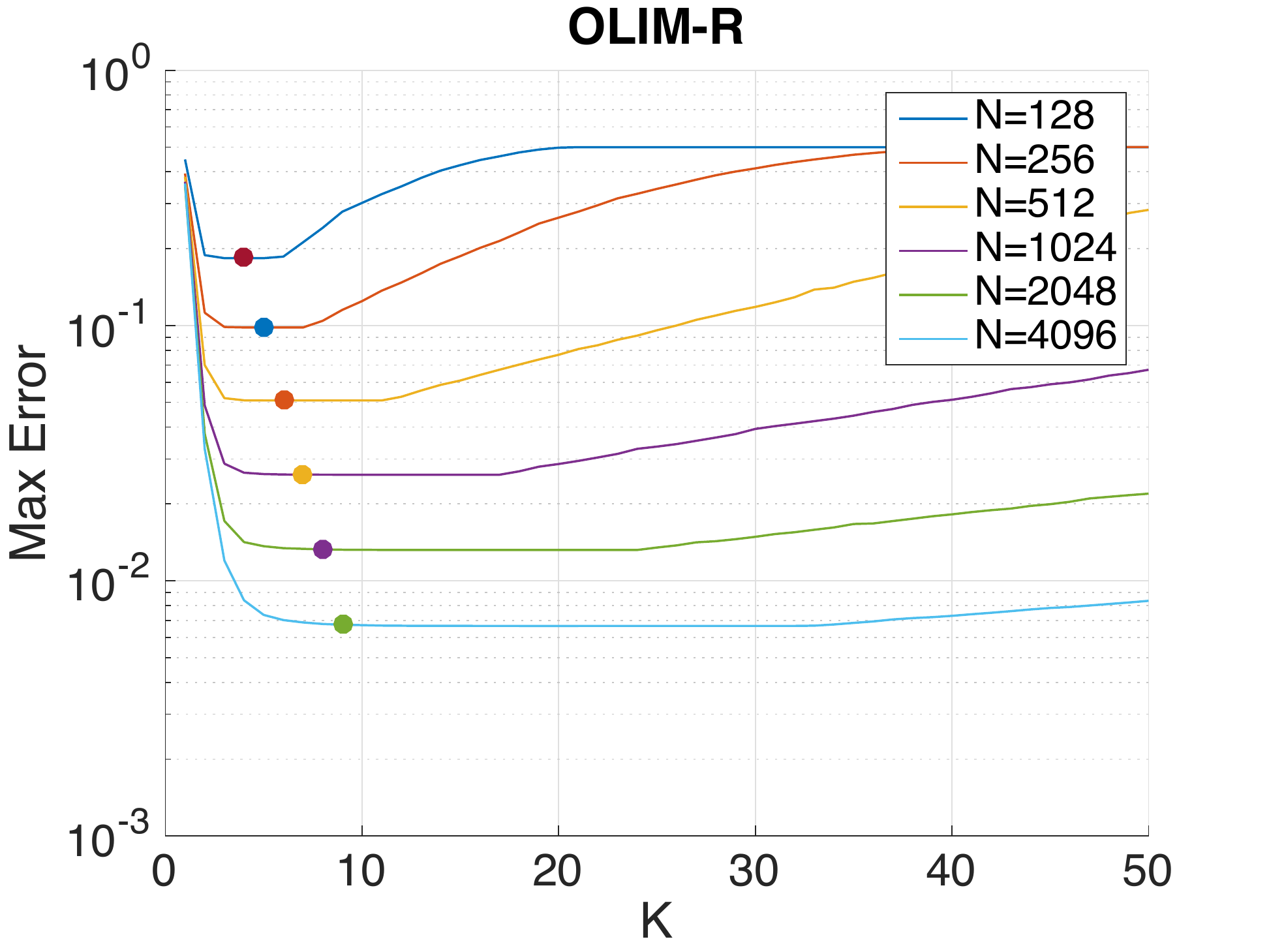}

(e)\includegraphics[width = 0.3\textwidth]{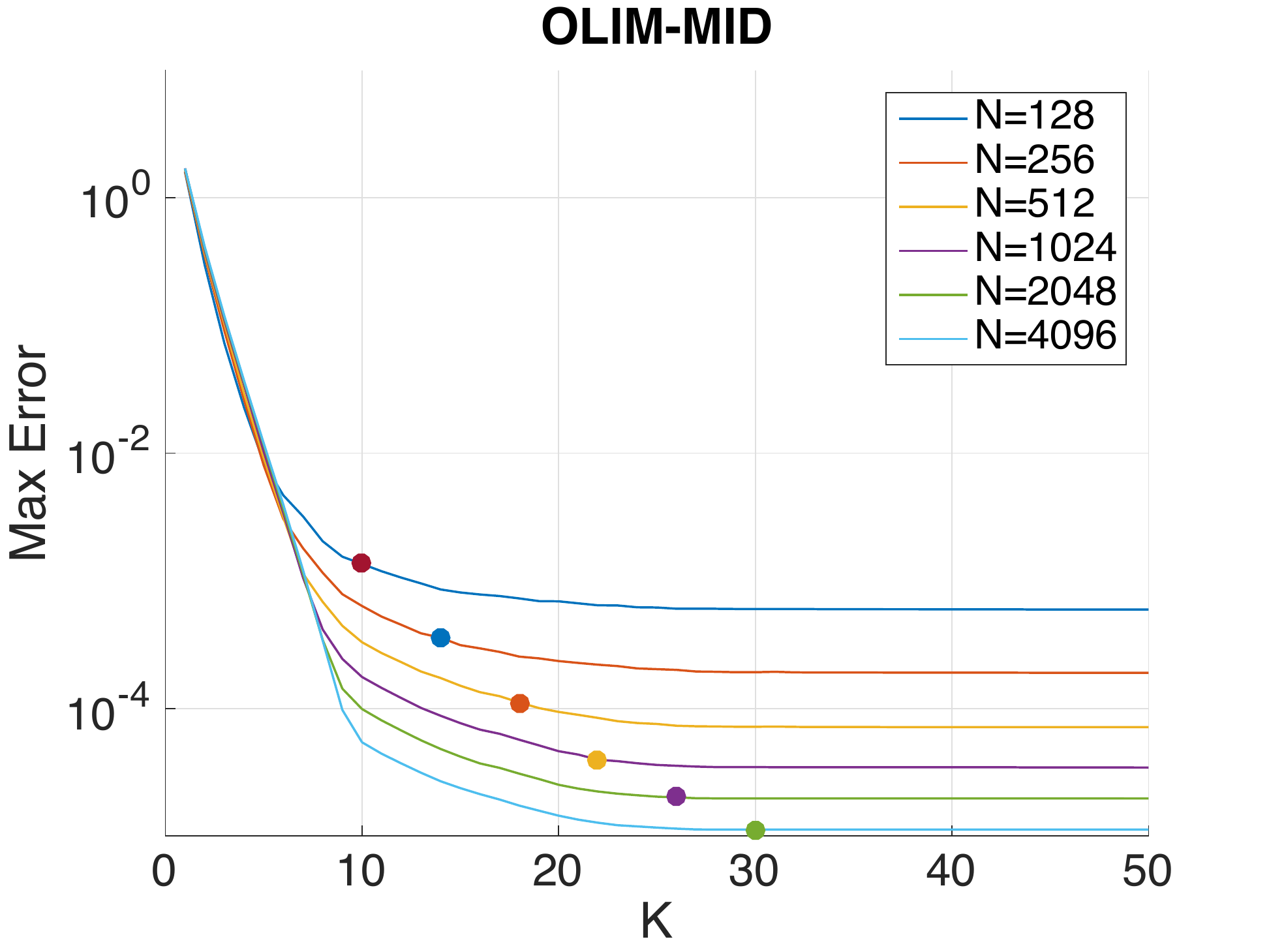}
(f)\includegraphics[width = 0.3\textwidth]{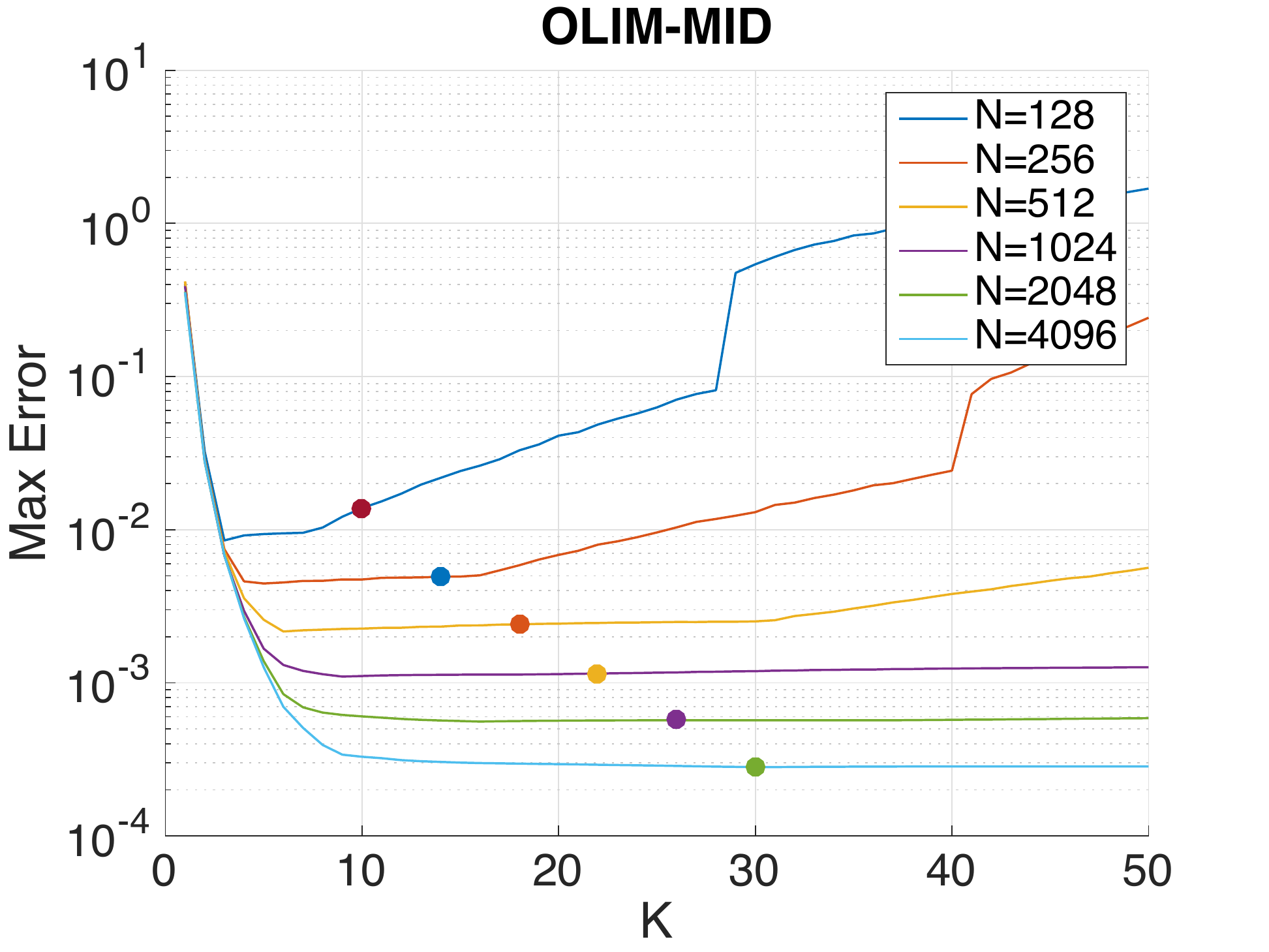}

(g)\includegraphics[width = 0.3\textwidth]{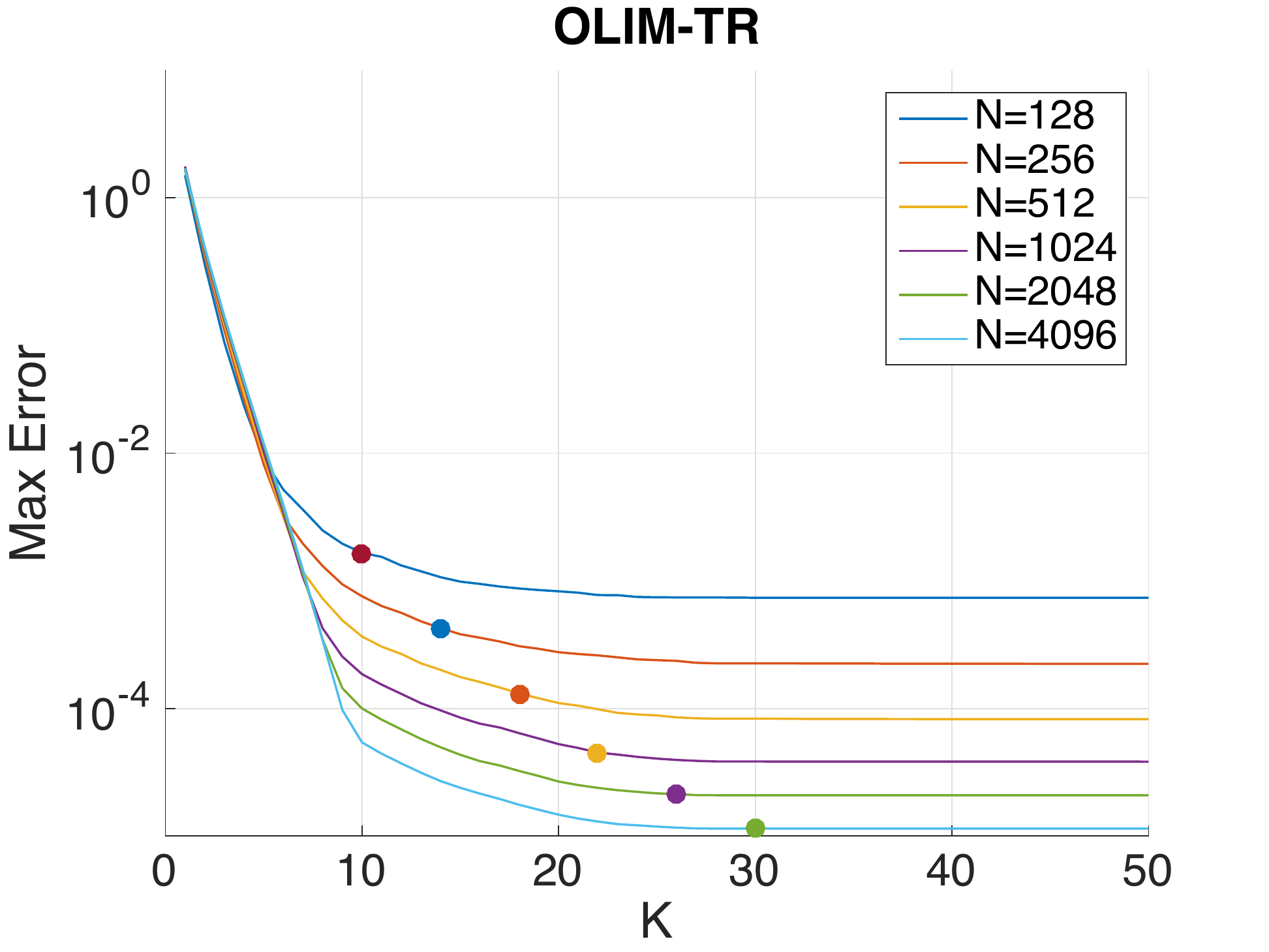}
(h)\includegraphics[width = 0.3\textwidth]{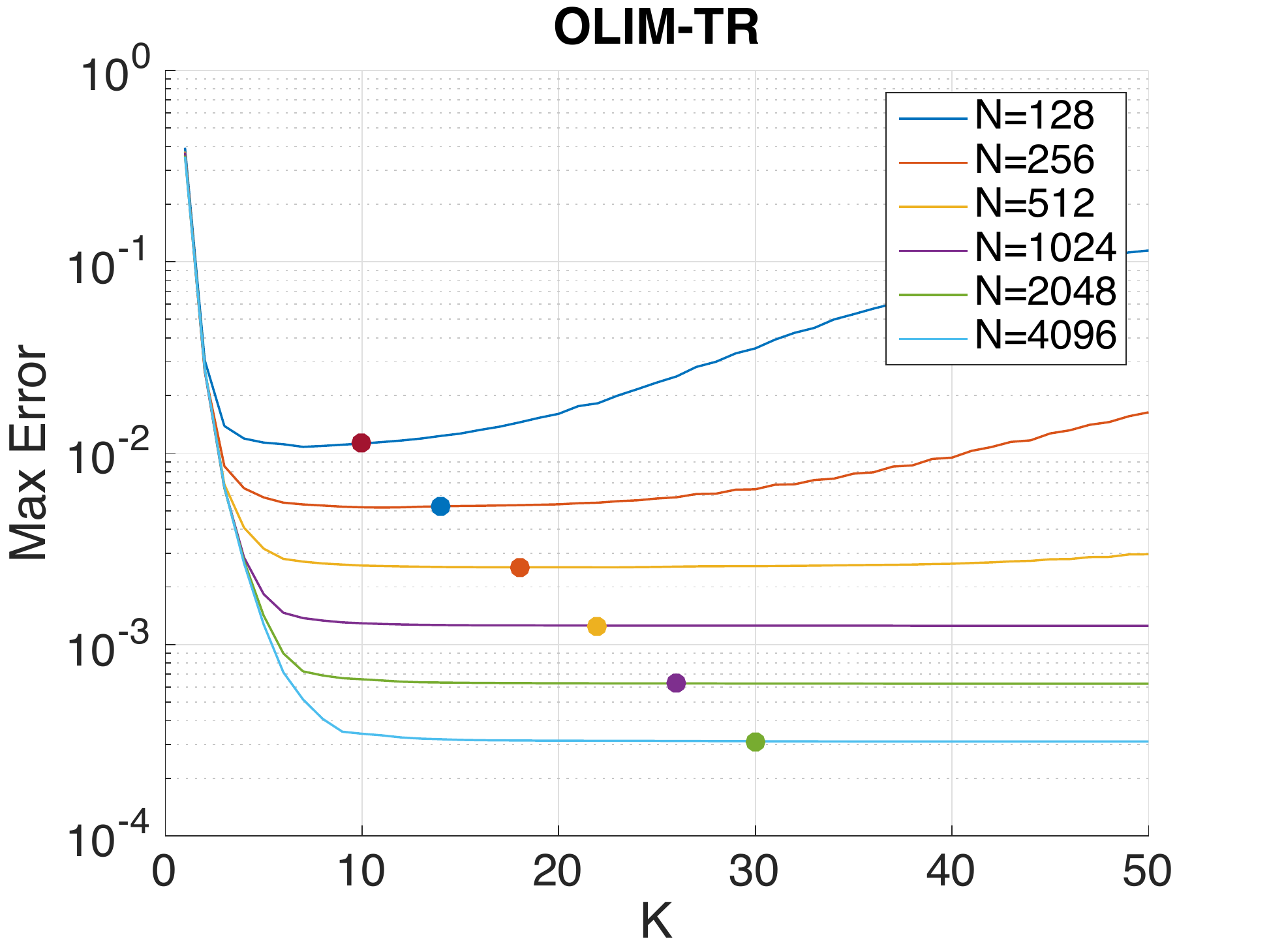}

(i)\includegraphics[width = 0.3\textwidth]{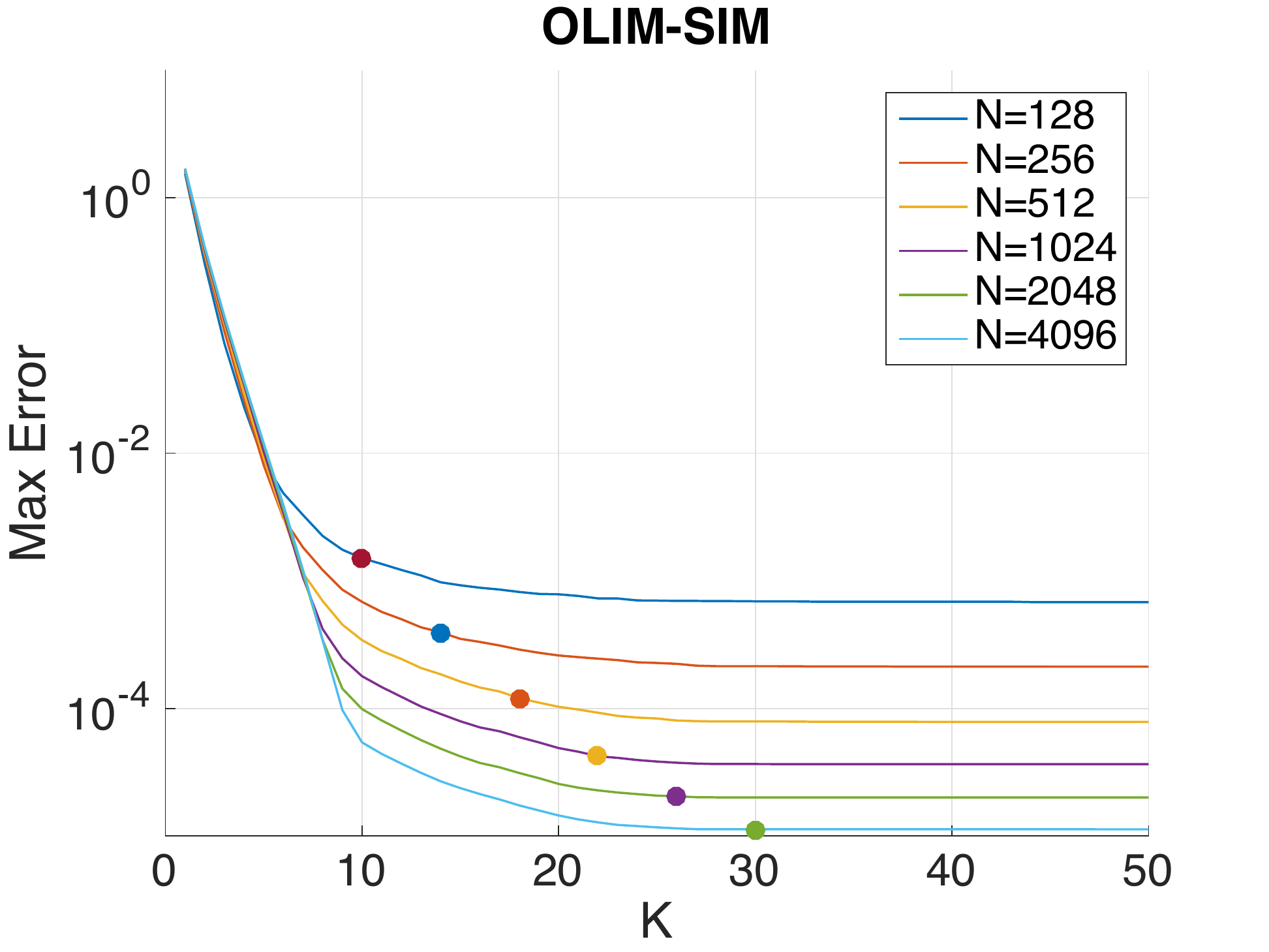}
(j)\includegraphics[width = 0.3\textwidth]{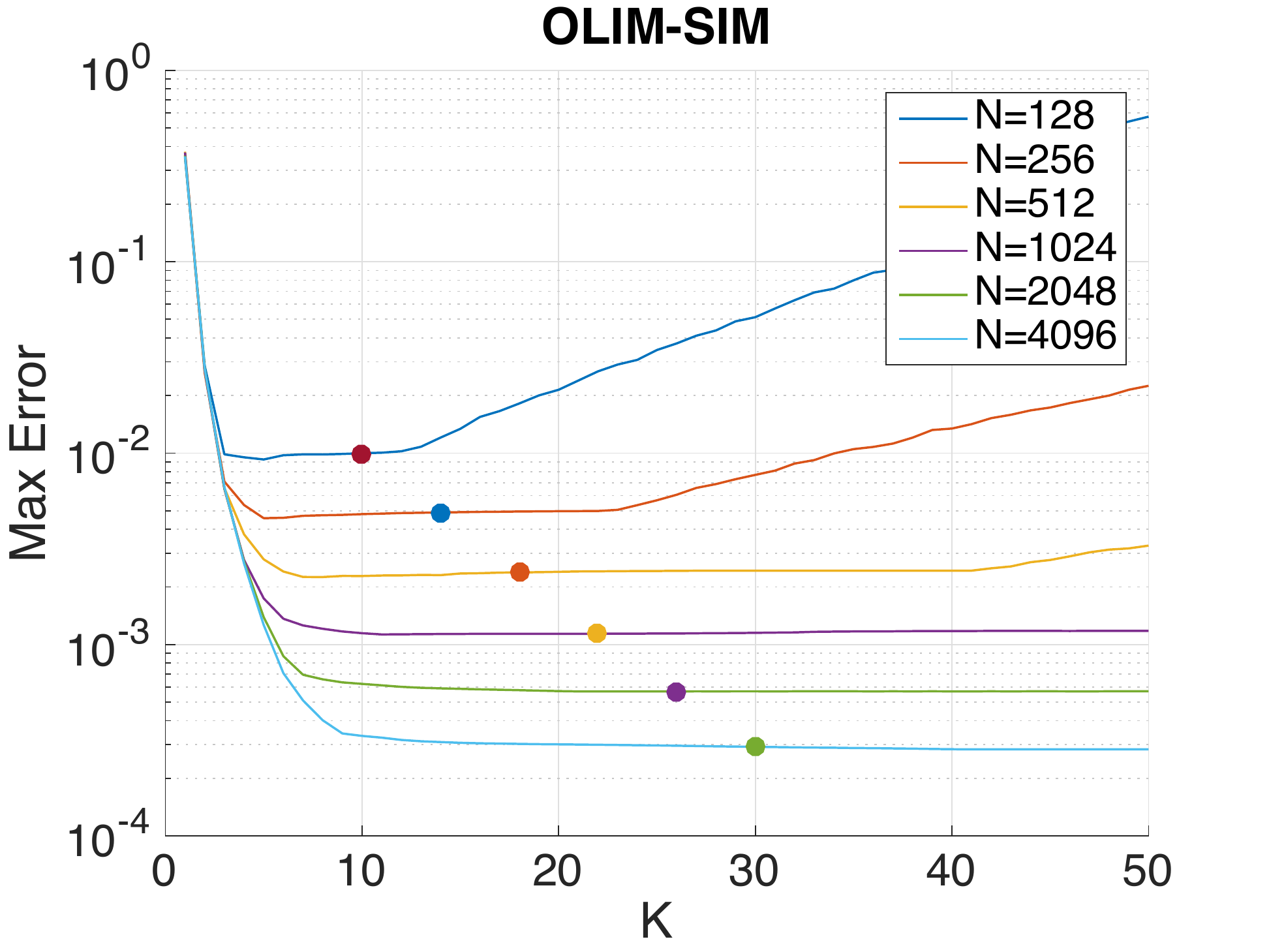}
\caption{The dependence of the maximum absolute error of the update factor $K$ 
for the OUM (a,b), OLIM-R (c,d), OLIM-MID (e,f), OLIM-TR (g,h), and OLIM-SIM (i,j).
Left column (a,c,e,g,i): SDE \eqref{linSDE}. Right column (b,d,f,h,j): SDE \eqref{cirSDE}.
The values of $K$ corresponding to Rules-of-Thumb \eqref{rule1} and \eqref{rule2} are marked on each graph.
}
\end{center}
\label{fig:K}
\end{figure*}

The results of our measurements are shown in Fig. \ref{fig:K}. For all methods, and all mesh sizes $2^p\times 2^p$, 
we plot the maximal absolute error versus $K$ for  SDEs  \eqref{linSDE} and \eqref{cirSDE}.
Note that the optimal values of $K$ are larger for SDE \eqref{linSDE} than those for SDE \eqref{cirSDE}. 

Based on our plots, we propose the following Rules-of-Thumb for choosing $K$. 
The points of the graphs corresponding to the proposed Rules-of-Thumb are marked on the graphs by large dots.

{\bf The Rule-of-Thumb for the OUM and OLIM-R.} For an $N\times N$ mesh where $2^7\le N\le 2^{12}$, pick
\begin{equation}
\label{rule1}
K(N) = {\sf round}[ \log_2N] - 3.
\end{equation}

{\bf The Rule-of-Thumb for OLIM-MID, OLIM-TR, and OLIM-SIM.} For an $N\times N$ mesh where $2^7\le N\le 2^{12}$, pick
\begin{equation}
\label{rule2}
K(N) =  10  + 4({\sf round}[\log_2N] - 7).
\end{equation}

%%%%

\subsection{Comparison of OLIM-MID, OLIM-TR, and OLIM-SIM to OLIM-R and the OUM}
Our comparison of OLIM-MID, OLIM-TR, and OLIM-SIM to OLIM-R and the OUM
is unambiguously in favor of the former ones.
The comparison is conducted using the update factor $K(p) = p - 3$ according to the Rule-of-Thumb for the OUM and OLIM-R, 
i.e., in the way benefiting the OUM and OLIM-R rather than the OLIMs with higher order quadrature rules. 

\begin{figure*}
\begin{center}

(a)\includegraphics[width = 0.4\textwidth]{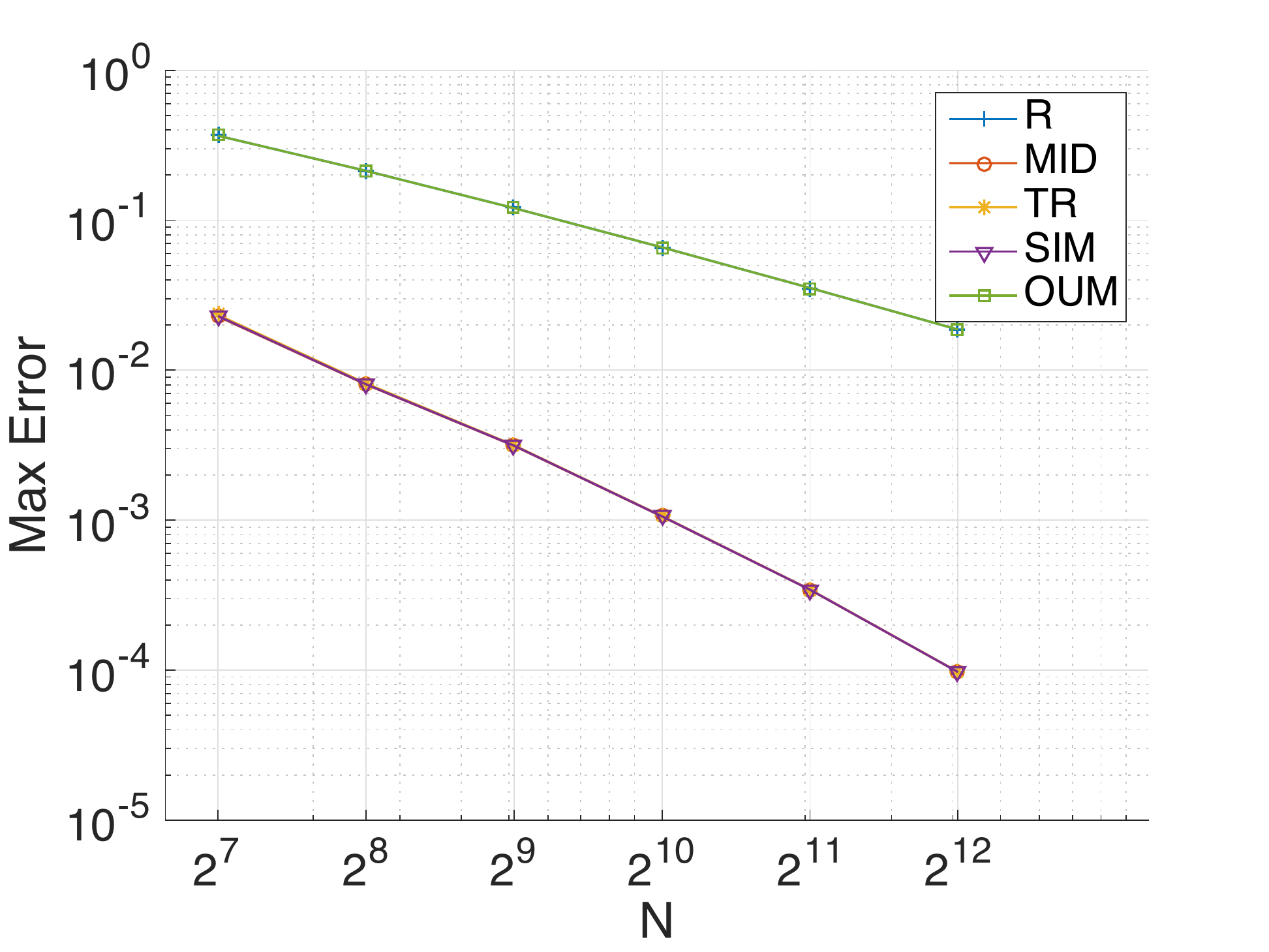}
(b)\includegraphics[width = 0.4\textwidth]{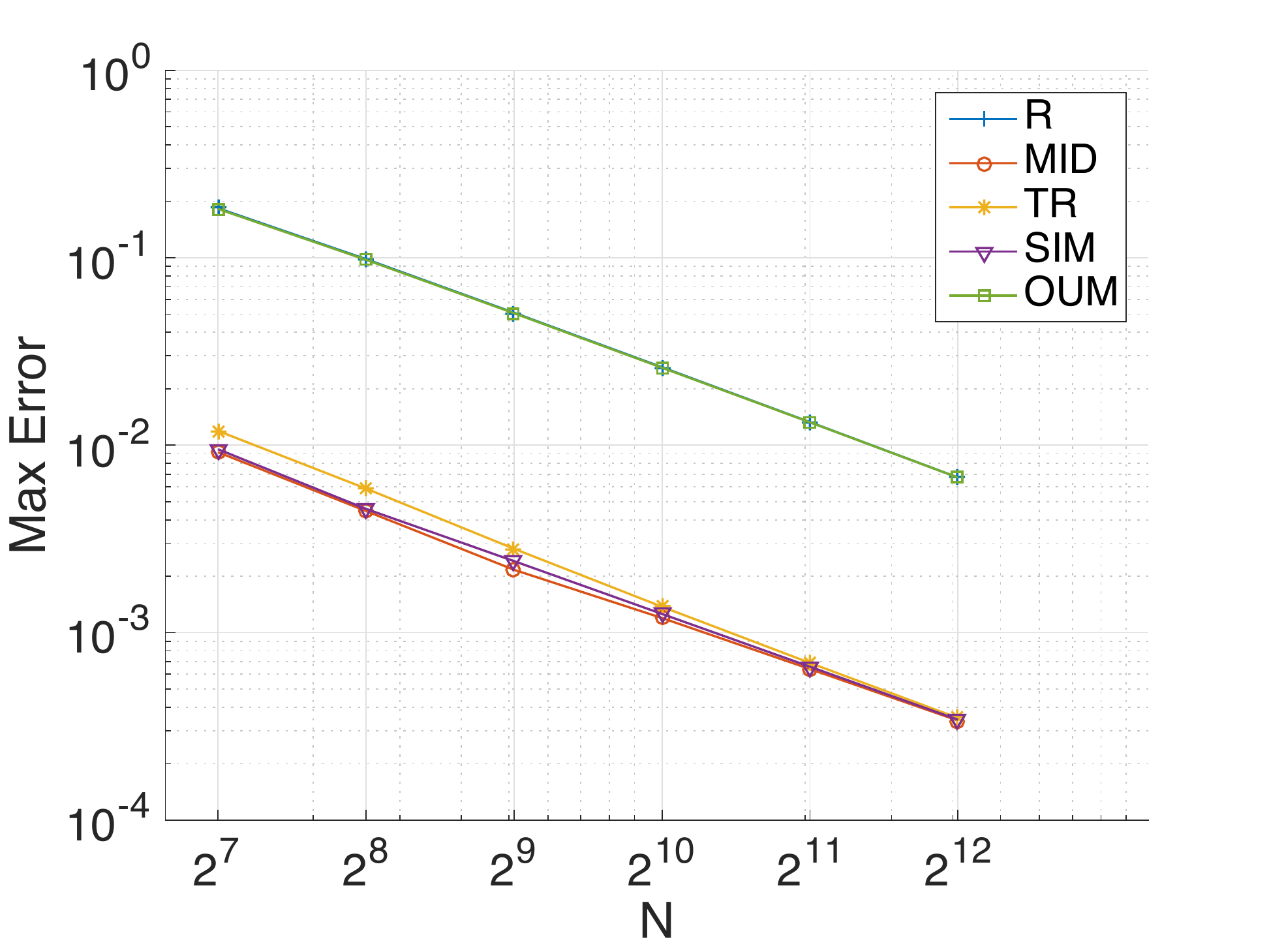}

(c)\includegraphics[width = 0.4\textwidth]{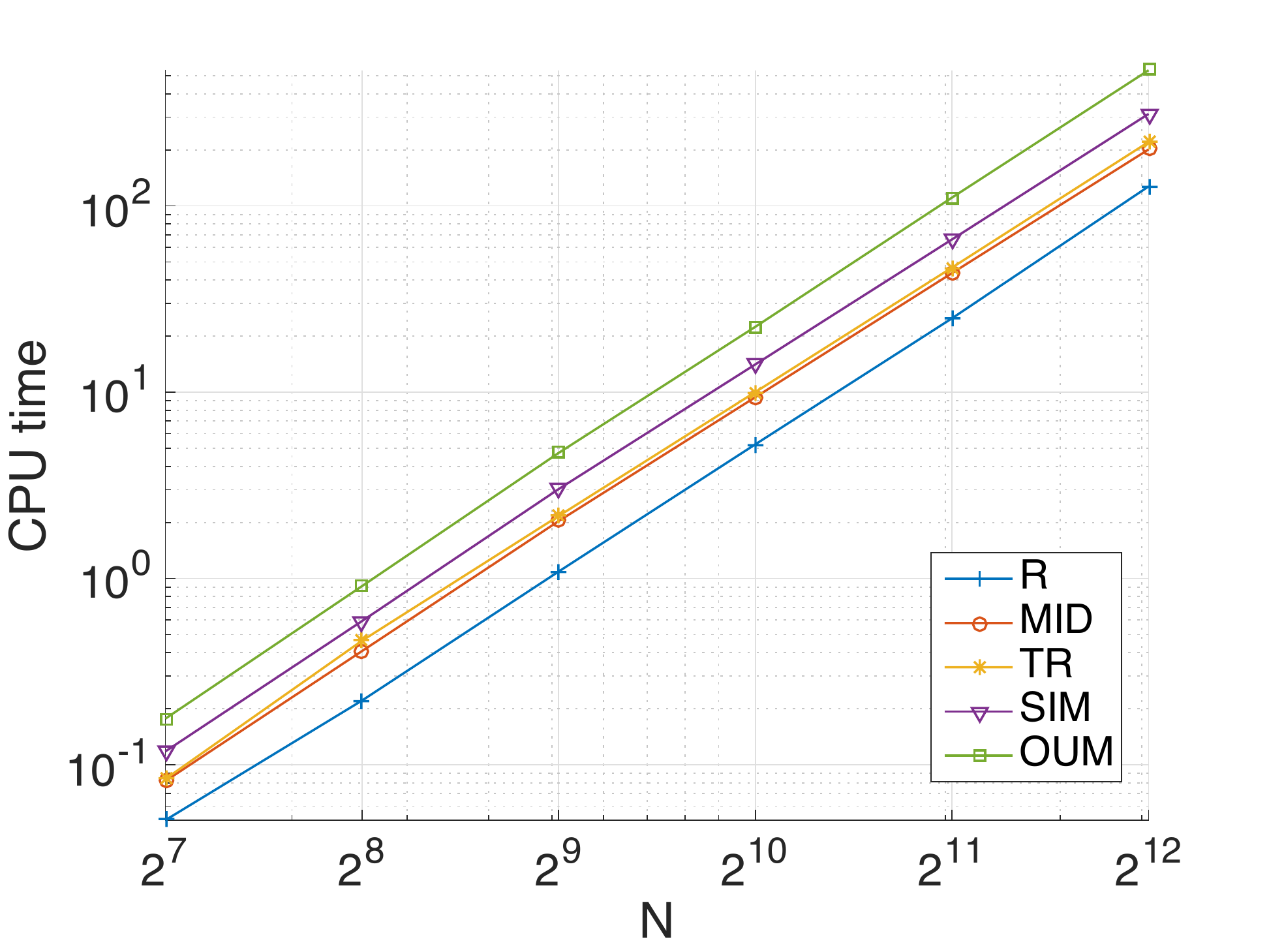}
(d)\includegraphics[width = 0.4\textwidth]{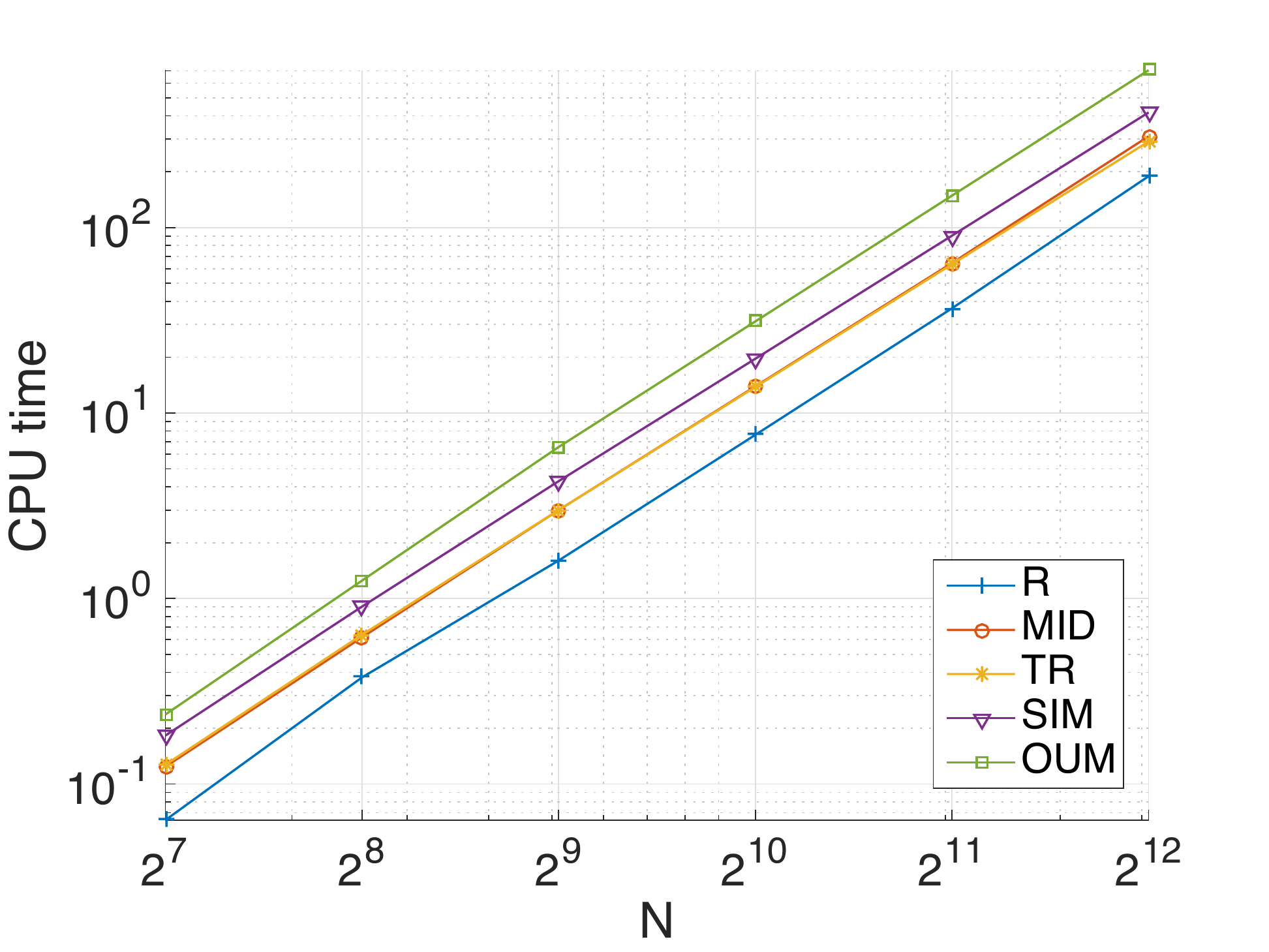}

(e)\includegraphics[width = 0.4\textwidth]{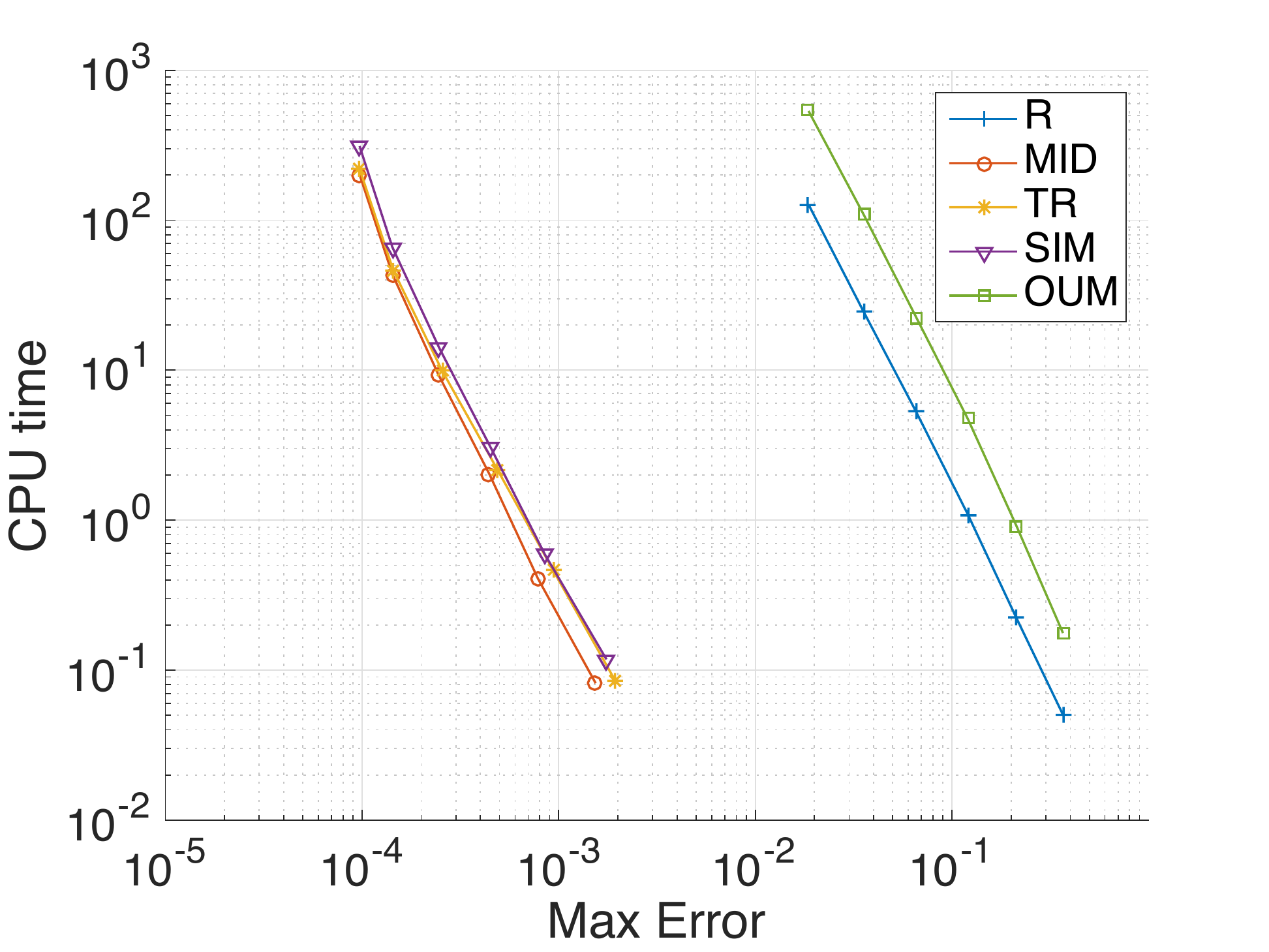}
(f)\includegraphics[width = 0.4\textwidth]{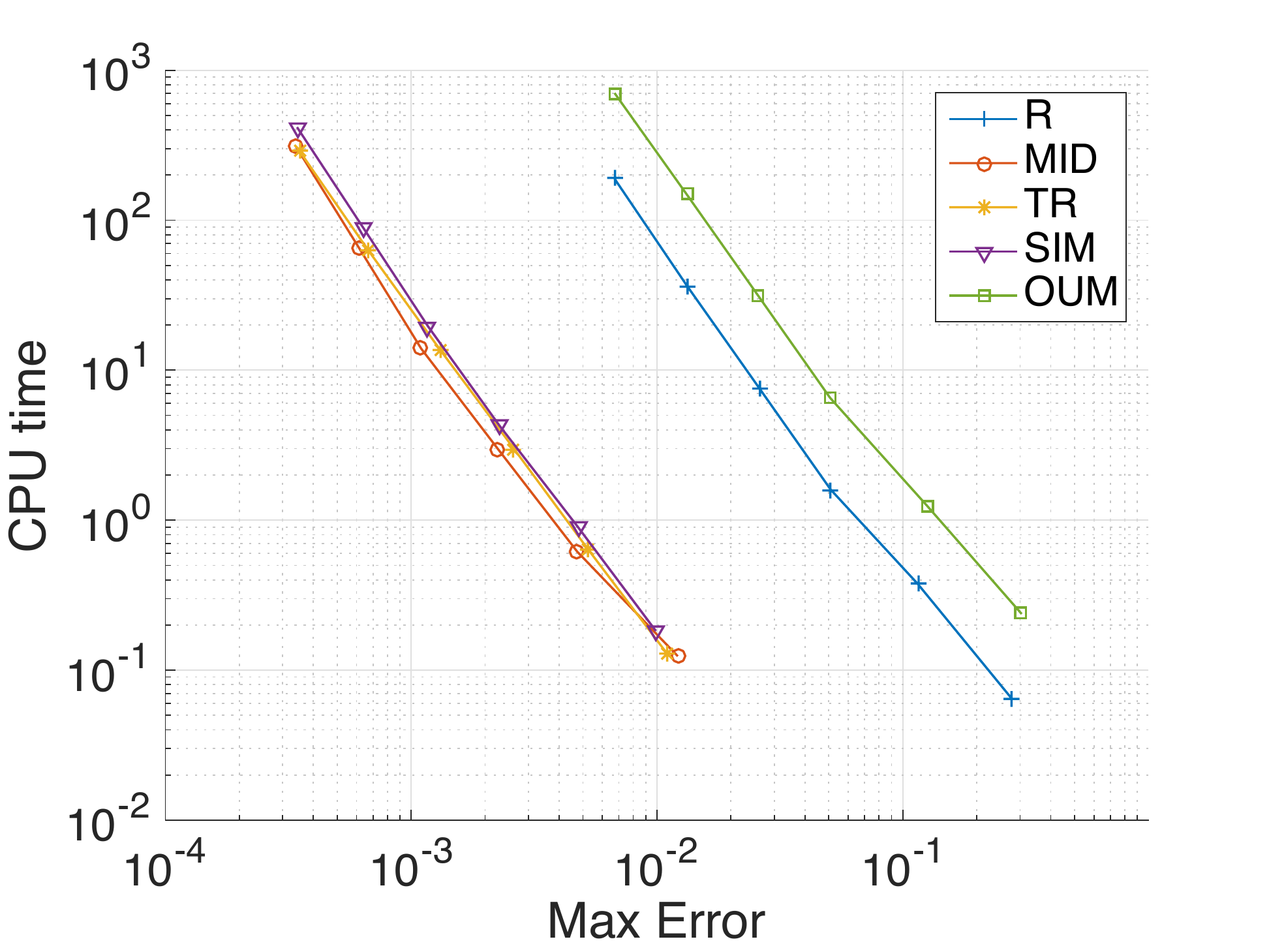}

\caption{Comparison of OLIM-MID, OLIM-TR, and OLIM-SIM to OLIM-R and the OUM.
The computational domain is $N\times N$, $N = 2^p$, $7\le p\le 12$.
(a,b): The maximum absolute error versus $N$.
(c,d): The CPU time versus $N$.
(e,f): The CPU time versus the maximal error.
Left column (a,c,e): SDE \eqref{linSDE}. Right column (b,d,f): SDE \eqref{cirSDE}.
}
\end{center}
\label{fig:comparison1}
\end{figure*}

\begin{figure*}
\begin{center}

(a)\includegraphics[width = 0.4\textwidth]{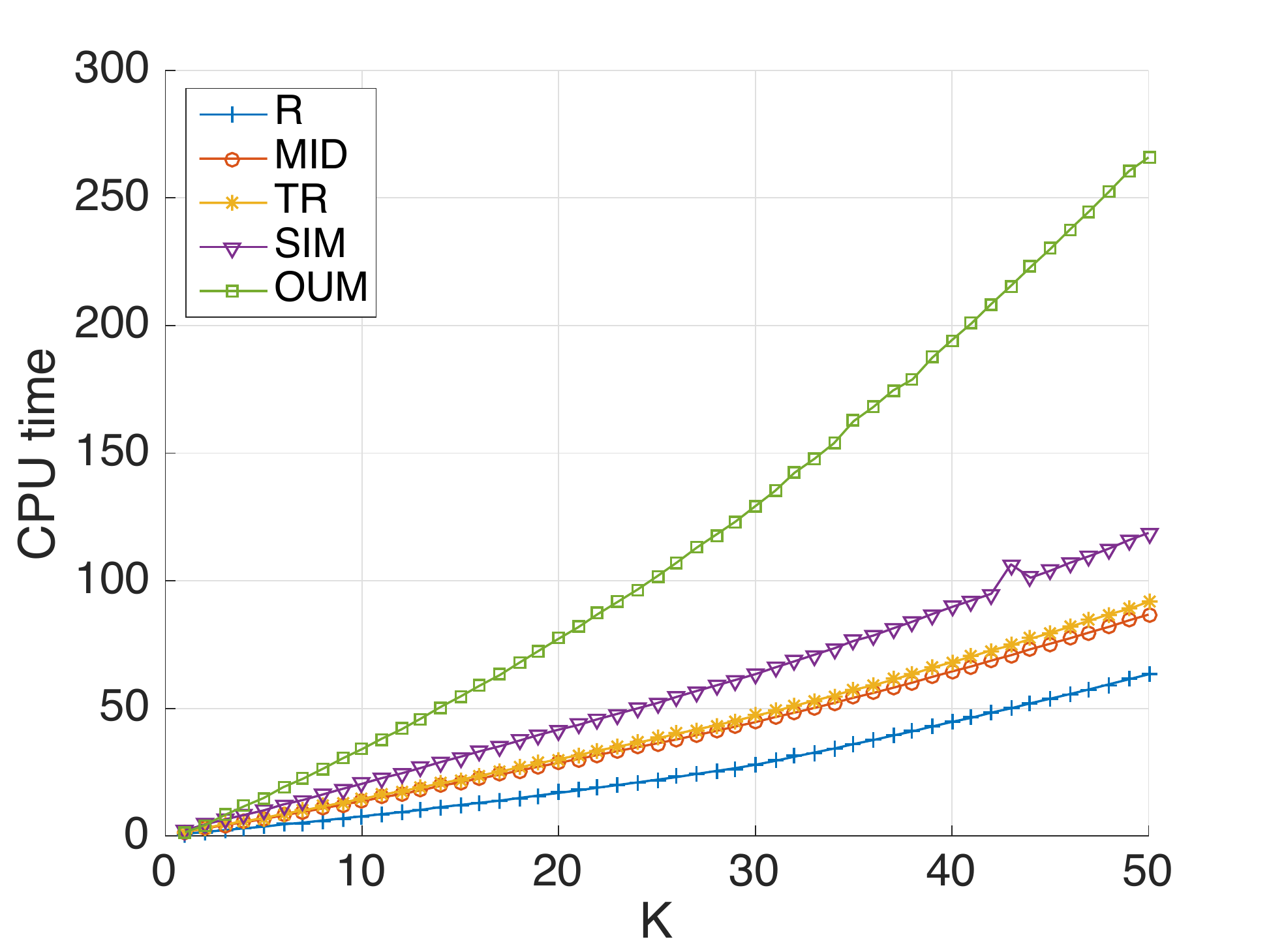}
(b)\includegraphics[width = 0.4\textwidth]{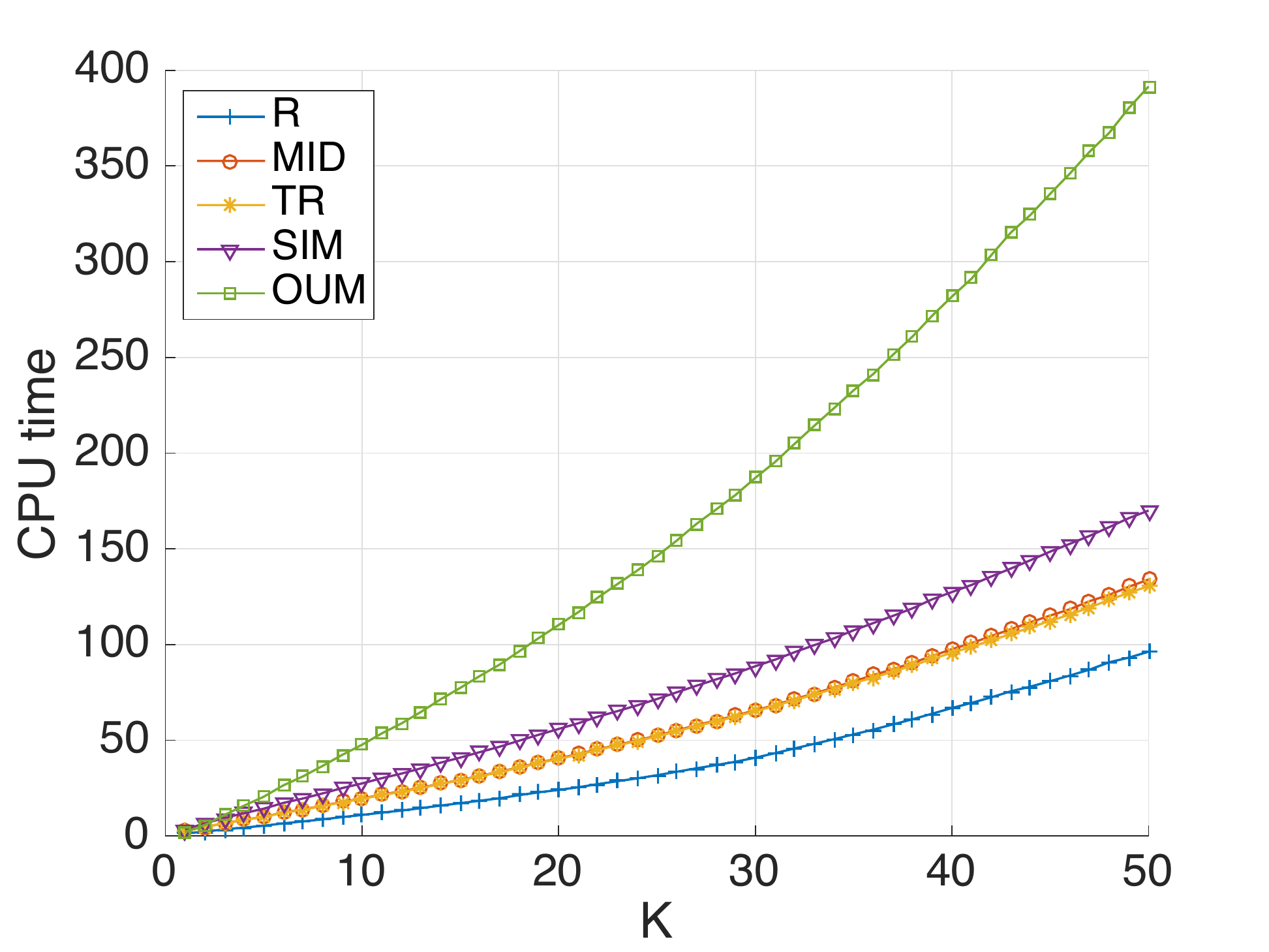}

\caption{
The CPU time (in seconds) plotted versus $K$ for $N=1024$ for SDEs \eqref{linSDE} (a) and \eqref{cirSDE} (b).
}
\label{fig:CPUvsK}
\end{center}
\end{figure*}

The graphs of the maximum absolute error and the CPU time as functions of $N$ respectively are shown in Fig. \ref{fig:comparison1} (a,b,c,d).
The least squares fits to the formulas $E = C N^{-q}$ and $T = CN^{q}$ for 
the maximum absolute errors and the CPU times  respectively are given in Table \ref{Table:LSfit}.
Fig. \ref{fig:comparison1} (a,b) and Table  \ref{Table:LSfit} (Columns 2 and 4) show that  OLIM-MID, OLIM-TR, 
and OLIM-SIM { are} 100 to 1000 times more accurate  then
OLIM-R and the OUM for values of  $K$ optimized for the latter methods. 
Fig.  \ref{fig:comparison1} (c,d) and Table  \ref{Table:LSfit} (Columns 3 and 5) show that  OLIM-MID, OLIM-TR, 
and OLIM-SIM are faster than the OUM. OLIM-MID is faster than the OUM at least by the factor of 1.5. OLIM-R is faster than the OUM at least by the factor of 3.
The CPU time (in seconds) versus the update factor $K$ for $N=1024$ is plotted 
for the OLIMs and the OUM in Fig. \ref{fig:CPUvsK}. 
These plots illustrate the advantage of our time-saving update strategy.

The CPU time as the function of the maximum absolute error is plotted in Fig. \ref{fig:comparison1}. It is clear that 
OLIM-MID, OLIM-TR, and OLIM-SIM are significantly better methods than  
OLIM-R and the OUM in terms of the balance between the accuracy and the CPU time.

\begin{table}
\caption{The least squares fits to the formulas $E = C N^{-q}$ and $T = CN^{q}$ for 
the maximum absolute errors and the CPU times as functions of $N$  ($N = 2^p$) respectively 
for SDEs \eqref{linSDE} (Columns 2 and 3) and \eqref{cirSDE} (Columns 4 and 5).
The update factor $K(p) =p - 3$  is chosen according to the Rule-of-Thumb for OUM and OLIM-R.}

\label{Table:LSfit}      % Give a unique label
\begin{tabular}{lllll}
\hline\noalign{\smallskip}
Method & Max Error & CPU time & Max Error & CPU time  \\
\noalign{\smallskip}\hline\noalign{\smallskip}
OUM & $24.8 \cdot N^{-0.860}$ &  $2.48\cdot 10^{-6}\cdot N^{2.31}$ & $19.1\cdot N^{-0.954}$ &  $3.52\cdot 10^{-6}\cdot N^{2.30}$\\
OLIM-R & $24.8 \cdot N^{-0.860}$ & $0.817 \cdot 10^{-6} \cdot N^{2.26}$ & $19.5\cdot N^{-0.957}$ & $1.09\cdot 10^{-6} \cdot N^{2.27}$\\
OLIM-MID & $47.5\cdot N^{-1.56}$ & $1.55\cdot 10^{-6}\cdot N^{2.25}$ & $0.843\cdot N^{-0.944}$ & $2.32\cdot 10^{-6}\cdot N^{2.25}$\\
OLIM-TR & $50.3\cdot N^{-1.57}$ & $1.60\cdot 10^{-6}\cdot N^{2.26}$ & $1.67\cdot N^{-1.02}$ & $2.71\cdot 10^{-6}\cdot N^{2.23}$\\
OLIM-SIM & $48.3\cdot N^{-1.56}$ & $2.01\cdot 10^{-6}\cdot N^{2.27}$ & $0.923\cdot N^{-0.951}$ & $3.87\cdot 10^{-6}\cdot N^{2.23}$\\
\noalign{\smallskip}\hline
\end{tabular}
\end{table}

%%%
{
\subsection{Effects of the hierarchical update strategy on numerical errors in OLIM-R}
\label{sec:HUS}
The one-point update and the triangle update in OLIM-R and the OUM are equivalent.  Therefore,
if the CPU-saving hierarchical update strategy would not be implemented in  OLIM-R, 
the numerical errors produced by these methods would coincide (in the exact arithmetics).
However, the hierarchical update combined with numerical errors might lead to non-identical numerical solutions by the OUM and  OLIM-R.
Figs. \ref{fig:comparison1}(a) and (b) show that the maximal errors 
in these two methods are very close (the curves visually coincide).
More detailed data for the maximal and RMS errors and CPU times in the OUM and OLIM-R are displayed in Tables  \ref{Table:HUS-l} and \ref{Table:HUS-c}
for SDEs \eqref{linSDE} and \eqref{cirSDE} respectively.
They indicate that the hierarchical update, in some cases, might increase the numerical error, but this increase is negligible.
On the other hand, the CPU times in OLIM-R are  approximately 4 times smaller.
}
\begin{table}
\caption{
Comparison of Maximal Errors, RMS errors and CPU times of the OUM and  OLIM-R 
applied to SDE \eqref{linSDE}. 
%For each mesh size $N$, the Maximal and the RMS errors coincide
%up to the recorded precision
%starting from the maximal update factor $K$ listed. 
}
\label{Table:HUS-l}      % Give a unique label
\begin{tabular}{lllll}
\hline\noalign{\smallskip}
Method, $N$, $K$ & Max Error & RMS error & CPU time, seconds  \\
\noalign{\smallskip}\hline\noalign{\smallskip}
OUM, $N = 512$ & & &\\
$K = 3$ &	1.7669e-01&	1.0440e-01 &	2.05 \\
$K = 5$ &	1.2133e-01&	7.9878e-02 &	3.69\\
$K =7$ &	1.2058e-01&	7.9659e-02 &	5.60\\
\hline
OLIM-R, $N = 512$ & & &\\
$K = 3$ &	1.8368e-01 &	1.0706e-01	&0.55\\
$K = 5$ &	1.2133e-01&	7.9878e-02&	0.88\\
$K = 7$ &	1.2058e-01&	7.9659e-02&	1.26\\
\hline\hline
OUM, $N = 1024$ & & &\\
$K = 4$ &	8.3905e-02&	5.2161e-02	&11.63\\
$K = 6$ &	6.6225e-02&	4.4102e-02	&19.15\\
$K = 8$ &	6.5912e-02&	4.3992e-02	&26.29\\
\hline
OLIM-R, $N = 1024$ & & &\\
$K = 4$ &	8.4836e-02	&5.2346e-02&	3.02\\
$K = 6$ &	6.6225e-02	&4.4102e-02&	4.56\\
$K = 8$ &	6.5912e-02	&4.3992e-02&	6.06\\
\hline\hline
OUM, $N = 2048$ & & &\\
$K = 5$ &	4.1289e-02&	2.6584e-02	&62.80\\
$K = 7$ &	3.5510e-02	&2.3803e-02	&94.16\\
$K = 9$ &	3.5350e-02&	2.3743e-02	&127.43\\
\hline
OLIM-R, $N = 2048$ & & &\\
$K = 5$ &	4.1529e-02 &	2.6609e-02 &	15.04\\
$K = 7$ &	3.5510e-02 &	2.3803e-02 &	21.55\\
$K = 9$ &	3.5350e-02 &	2.3743e-02	&28.60\\
\hline\hline
OUM, $N = 4096$ & & &\\
$K = 5$	& 2.6959e-02	& 1.6629e-02	& 263.71\\
$K = 7$	& 1.9089e-02	& 1.2776e-02	& 397.53\\
$K =  9$	 & 1.8656e-02	& 1.2599e-02	& 535.13\\
\hline
OLIM-R, $N = 4096$ & & &\\
$K=5$ &	2.7204e-02 &	1.6702e-02	 & 66.76\\
$K=7$ &	1.9051e-02 &	1.2769e-02	 & 96.03\\
$K= 9$ &	1.8656e-02	 & 1.2599e-02	 & 127.22\\
\noalign{\smallskip}\hline
\end{tabular}
\end{table}

\begin{table}
\caption{
Comparison of Maximal Errors, RMS errors and CPU times of the OUM and OLIM-R 
applied to SDE \eqref{cirSDE}. 
}
\label{Table:HUS-c}      % Give a unique label
\begin{tabular}{lllll}
\hline\noalign{\smallskip}
Method, $N$, $K$ & Max Error & RMS error & CPU time, seconds  \\
\noalign{\smallskip}\hline\noalign{\smallskip}
OUM, $N = 512$ & & &\\
$K = 3$ &		5.1079e-02&		2.1963e-02&		2.71\\
$K = 5$ &		5.0563e-02&		2.1874e-02&		5.05\\
$K = 7$ &		5.0563e-02&		2.2332e-02&		7.71\\
\hline
OLIM-R, $N = 512$ & & &\\
$K = 3$ &		5.1925e-02&		2.2371e-02&		0.81\\
$K = 5$ &		5.0850e-02&		2.1940e-02&		1.29\\
$K = 7$ &		5.0845e-02&		2.2371e-02&		1.83\\
\hline\hline
OUM, $N = 1024$ & & &\\
$K = 5$ &		2.5930e-02&		1.1245e-02&		20.50\\
$K = 8$ &		2.5890e-02&		1.1447e-02&		36.65\\
$K = 11$ &		2.5890e-02&		1.1696e-02&		53.57\\
\hline
OLIM-R, $N = 1024$ & & &\\
$K = 5$ &		2.6182e-02&		1.1283e-02&		5.36\\
$K = 8$ &		2.6031e-02&		1.1454e-02&		8.71\\
$K = 11$ &		2.6025e-02&		1.1692e-02&		12.21\\
\hline\hline
OUM, $N = 2048$ & & &\\
$K = 5$ &		1.3478e-02&		5.8271e-03&		83.24\\
$K = 10$ &		1.3204e-02&		5.8695e-03&		193.00\\
$K = 15$ &		1.3189e-02&		6.0301e-03&		311.84\\
$K = 20$ &		1.3190e-02&		6.1811e-03&		443.35\\
\hline
OLIM-R, $N = 2048$ & & &\\
$K = 5$ &		1.3656e-02&		5.8558e-03&		22.78\\
$K = 10$ &		1.3215e-02&		5.8607e-03&		46.23\\
$K = 15$ &		1.3200e-02&		6.0178e-03&		72.81\\
$K = 20$ &		1.3201e-02&		6.1618e-03&		103.35\\
\hline\hline
OUM, $N = 4096$ & & &\\
$K = 5$ &		7.2210e-03&		3.0923e-03&		341.70\\
$K = 10$ &		6.6921e-03&		2.9666e-03&		791.00\\
$K = 15$ &		6.6527e-03&		3.0332e-03&		1276.04\\
$K = 20$ &		6.6506e-03&		3.0975e-03&		1811.45\\
\hline
OLIM-R, $N = 4096$ & & &\\
$K = 5$ &		7.3480e-03&		3.1267e-03&		101.85\\
$K = 10$ &		6.7011e-03&		2.9587e-03&		210.50\\
$K = 15$ &		6.6555e-03&		3.0230e-03&		326.48\\
$K = 20$ &		6.6536e-03&		3.0859e-03&		463.64\\
\noalign{\smallskip}\hline
\end{tabular}
\end{table}

%%%%

\subsection{Comparison of OLIM-MID, OLIM-TR, and OLIM-SIM to each other}
\label{sec:OLIM2}
\begin{figure*}
\begin{center}

(a)\includegraphics[width = 0.4\textwidth]{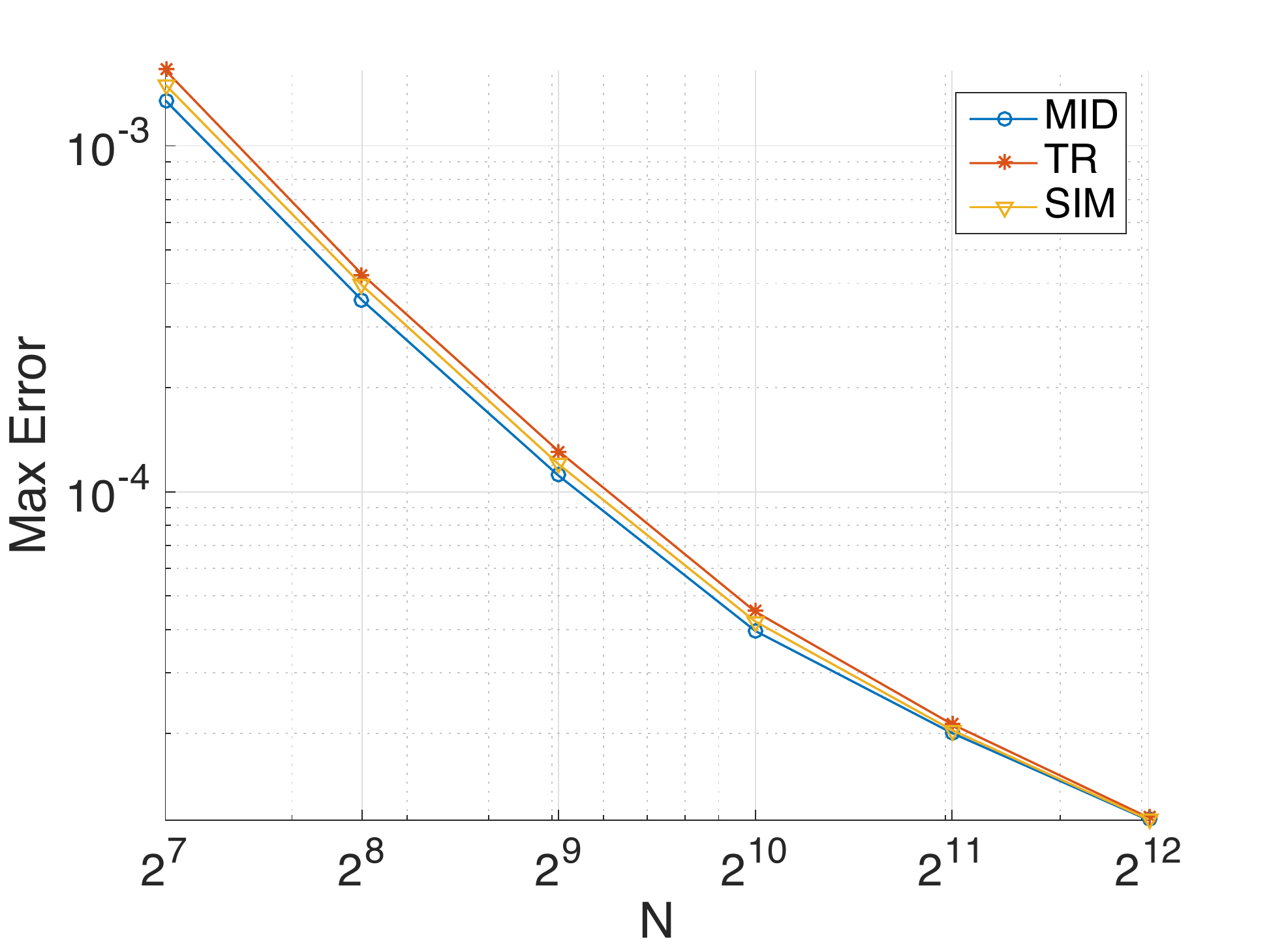}
(b)\includegraphics[width = 0.4\textwidth]{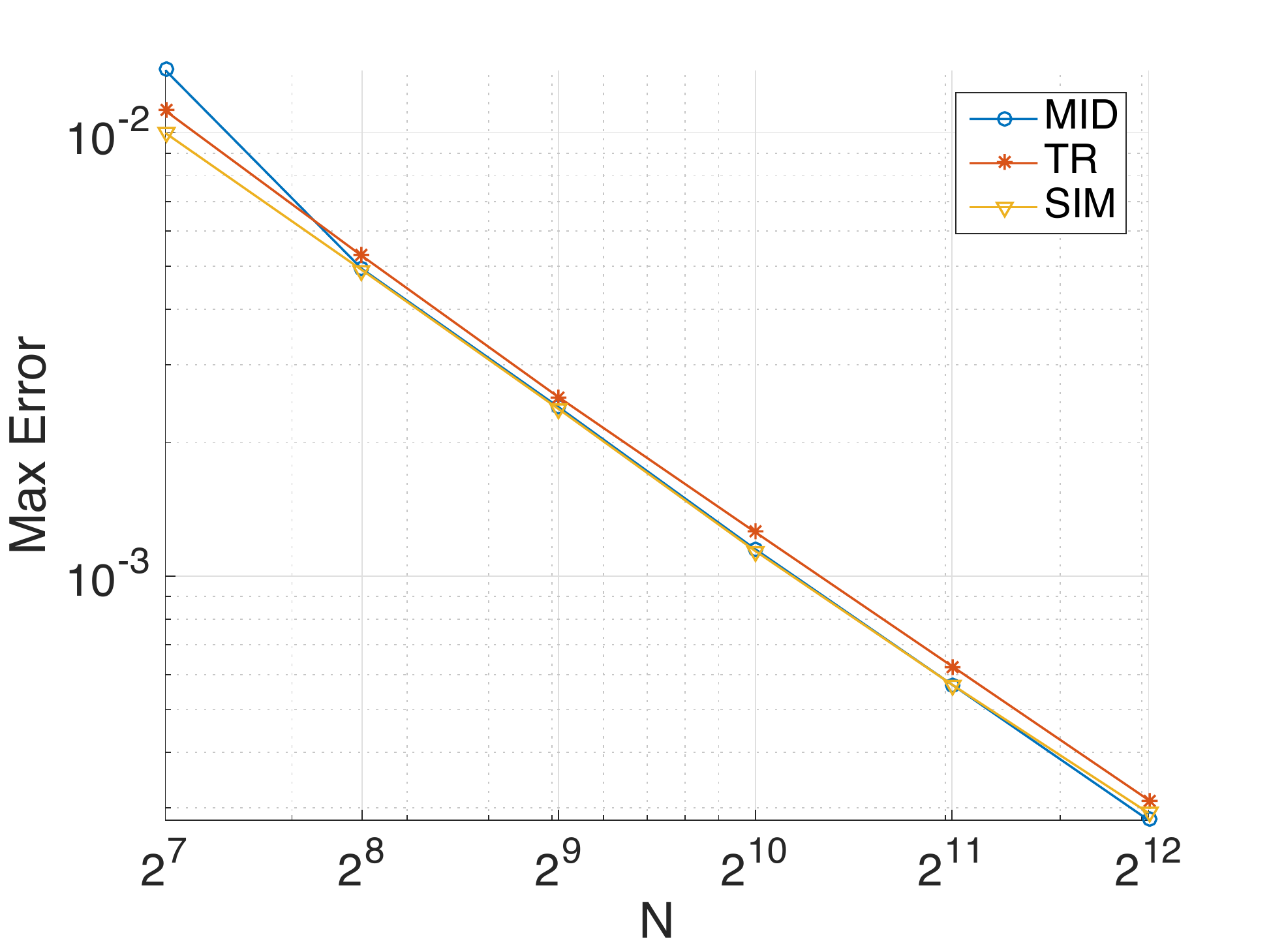}

(c)\includegraphics[width = 0.4\textwidth]{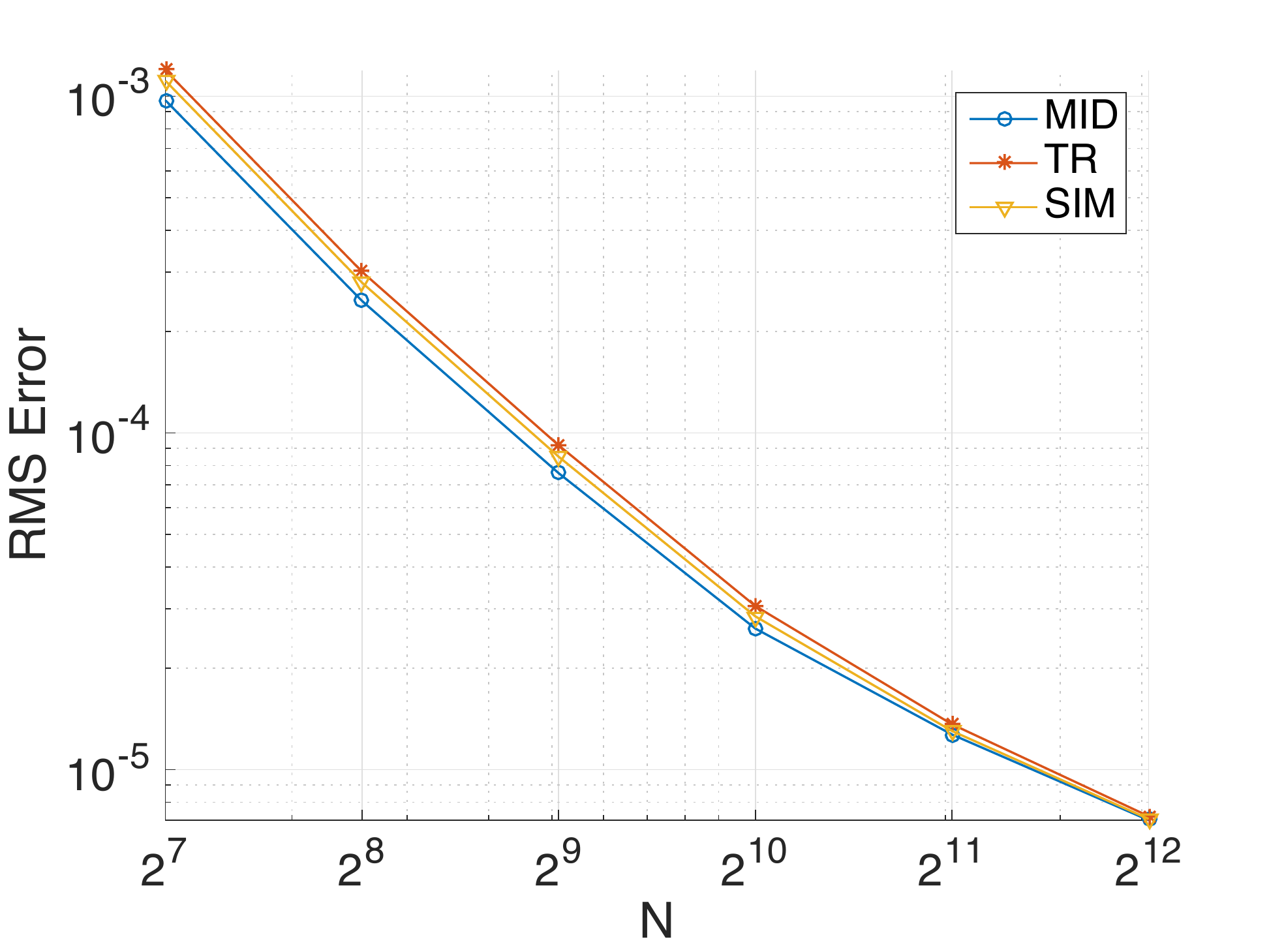}
(d)\includegraphics[width = 0.4\textwidth]{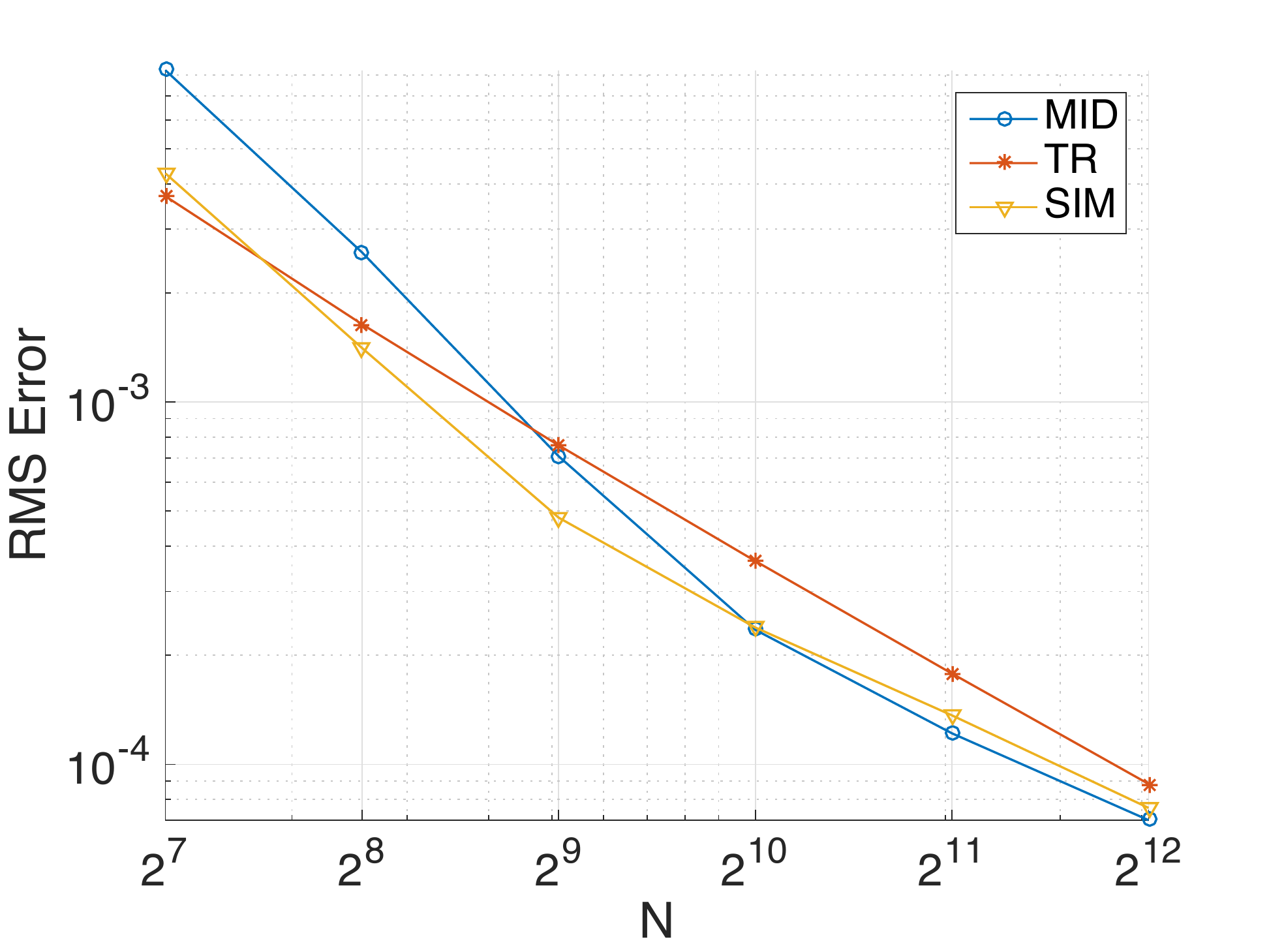}

(e)\includegraphics[width = 0.4\textwidth]{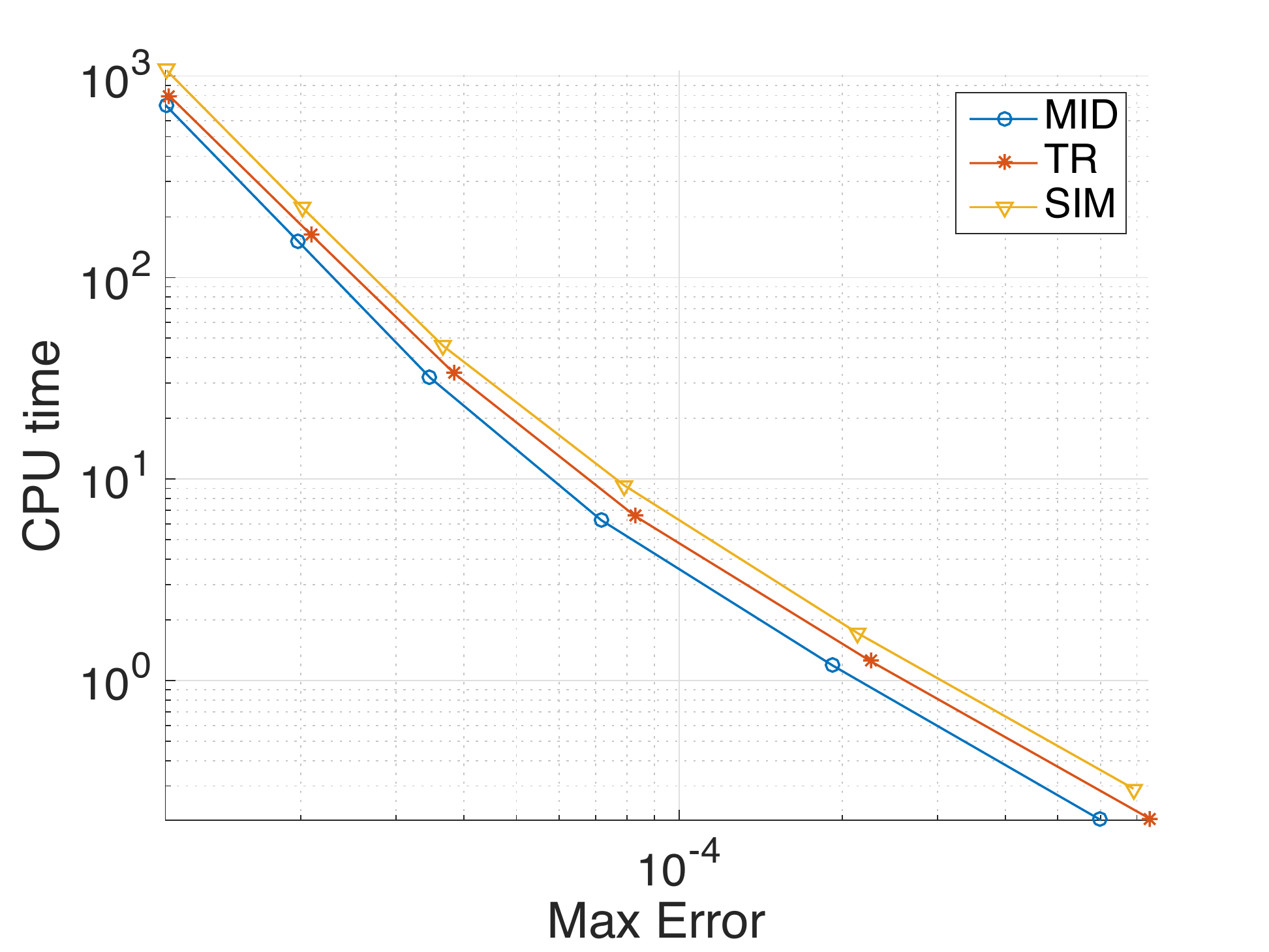}
(f)\includegraphics[width = 0.4\textwidth]{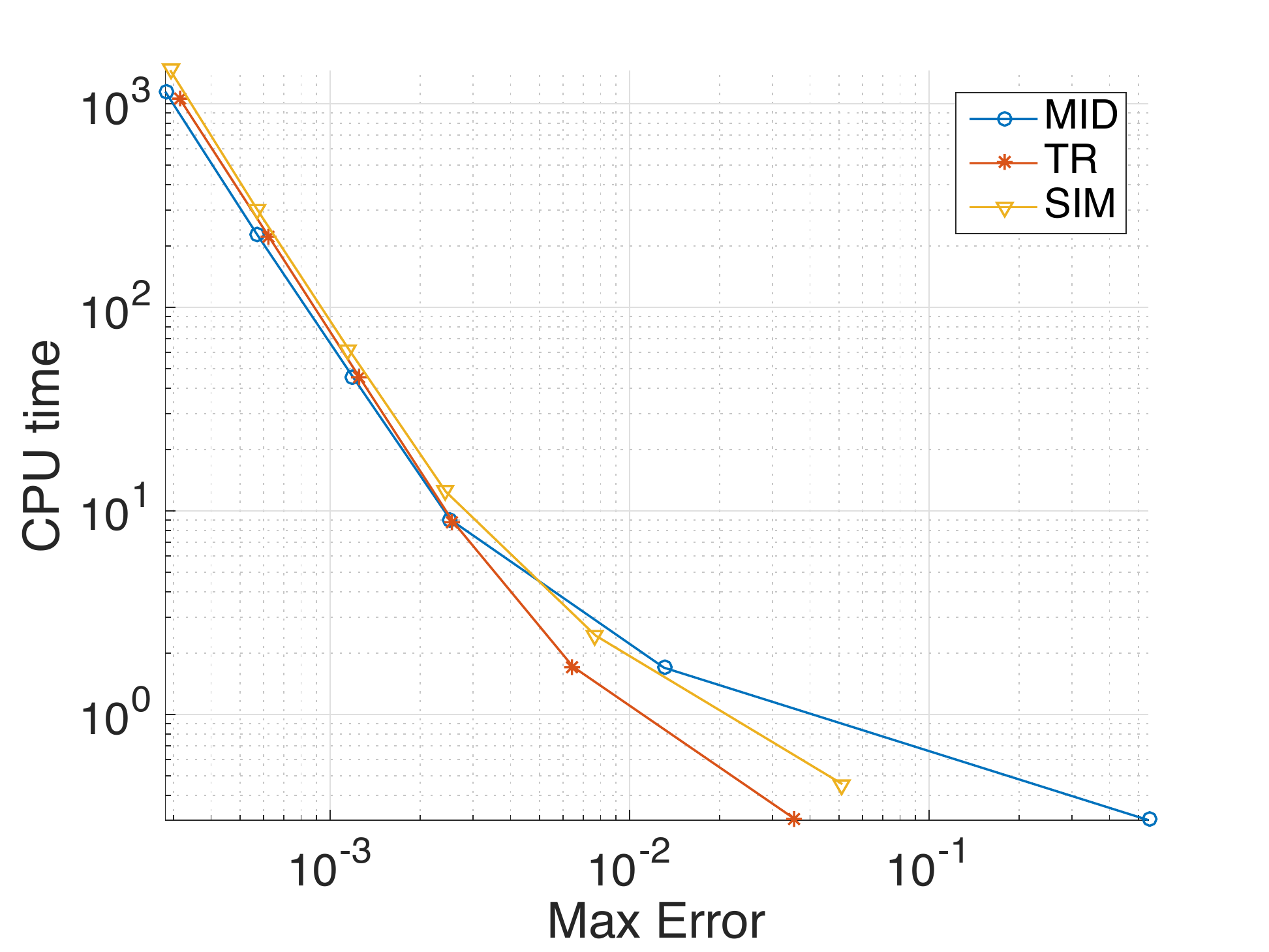}

\caption{Comparison of OLIM-MID, OLIM-TR, and OLIM-SIM to each other.
The computational domain is $N\times N$, $N = 2^p$, $7\le p\le 12$.
(a,b): The maximum absolute error versus $N$.
(c,d): The RMS error versus $N$.
(e,f): The CPU time versus the maximal error.
Left column (a,c,e): SDE \eqref{linSDE}. Right column (b,d,f): SDE \eqref{cirSDE}.
}
\end{center}
\label{fig:comparison2}
\end{figure*}

The comparison of OLIM-MID, OLIM-TR, and OLIM-SIM is conducted using $K(p) = 10 + 4(p - 7)$, i.e., 
according to the Rule-of-Thumb for OLIM-MID, OLIM-TR, and OLIM-SIM. 
%The computations are done on the  $N\times N$ mesh where $N=2^p$, $7\le p\le 12$.
Fig. \ref{fig:comparison2} shows that all these three methods are quite close in accuracy. Fig. \ref{fig:CPUvsK}
indicates that OLIM-SIM is somewhat slower than OLIM-MID and OLIM-TR.  
The least squares fits to the formulas $E = C N^{-q}$
for the maximum absolute errors and the RSM errors as functions of $N$  ($N = 2^p$) are given in Table \ref{Table:LSfit2}.
OLIM-MID gives the best results on SDE \eqref{linSDE} 
while OLIM-TR and OLIM-SIM challenge it on SDE \eqref{cirSDE} for rough meshes.

We favor OLIM-MID method. For fine meshes, it has the best balance between the 
accuracy and the CPU time among all methods
considered in this work. 

\begin{table}
\caption{
The least squares fits to the formulas $E = C N^{-q}$
for the maximum absolute errors and the RSM errors as functions of $N$  ($N = 2^p$) 
for SDEs \eqref{linSDE} (Columns 2 and 3) and \eqref{cirSDE} (Columns  4 and 5).
The update factor $K(p) = 10 + 4(p - 7)$  is chosen according to the { Rule-of-Thumb} for OLIM-MID, OLIM-TR, and OLIM-SIM.
}
\label{Table:LSfit2}      % Give a unique label
\begin{tabular}{lllll}
\hline\noalign{\smallskip}
Method & Max Error & RMS Error   & Max Error & RMS Error  \\
\noalign{\smallskip}\hline\noalign{\smallskip}
OLIM-MID & $0.817\cdot N^{-1.39}$  & $0.705\cdot N^{-1.43}$  & $2.47\cdot N^{-1.10}$  & $5.85\cdot N^{-1.41}$ \\
%\hline
OLIM-TR & $1.31\cdot N^{-1.44}$  & $1.16\cdot N^{-1.48}$  & $1.61\cdot N^{-1.03}$  & $0.646\cdot N^{-1.07}$ \\
%\hline
OLIM-SIM & $1.07\cdot N^{-1.42}$  & $0.99\cdot N^{-1.46}$  & $1.42\cdot N^{-1.02}$  & $0.846\cdot N^{-1.15}$ \\
\noalign{\smallskip}\hline
\end{tabular}
\end{table}

The error plots for each of the methods for SDEs \eqref{linSDE} and \eqref{cirSDE} for the $256\times 256$ meshes and $K=11$ 
(corresponds to the Rule-of-Thumb)
are shown in Fig. \ref{fig:errors}. 
We will return to the discussion on the error distributions in Section \ref{sec:errors} below. 
Now we just note that these error plots are consistent with the results of Section \ref{sec:errors}.
As we will show, the midpoint quadrature rule tends to underestimate the line integral, the trapezoid rule 
tends to exaggerate it, and the Simpson rule can make errors of either sign. 

\begin{figure*}
\begin{center}

(a)\includegraphics[width = 0.4\textwidth]{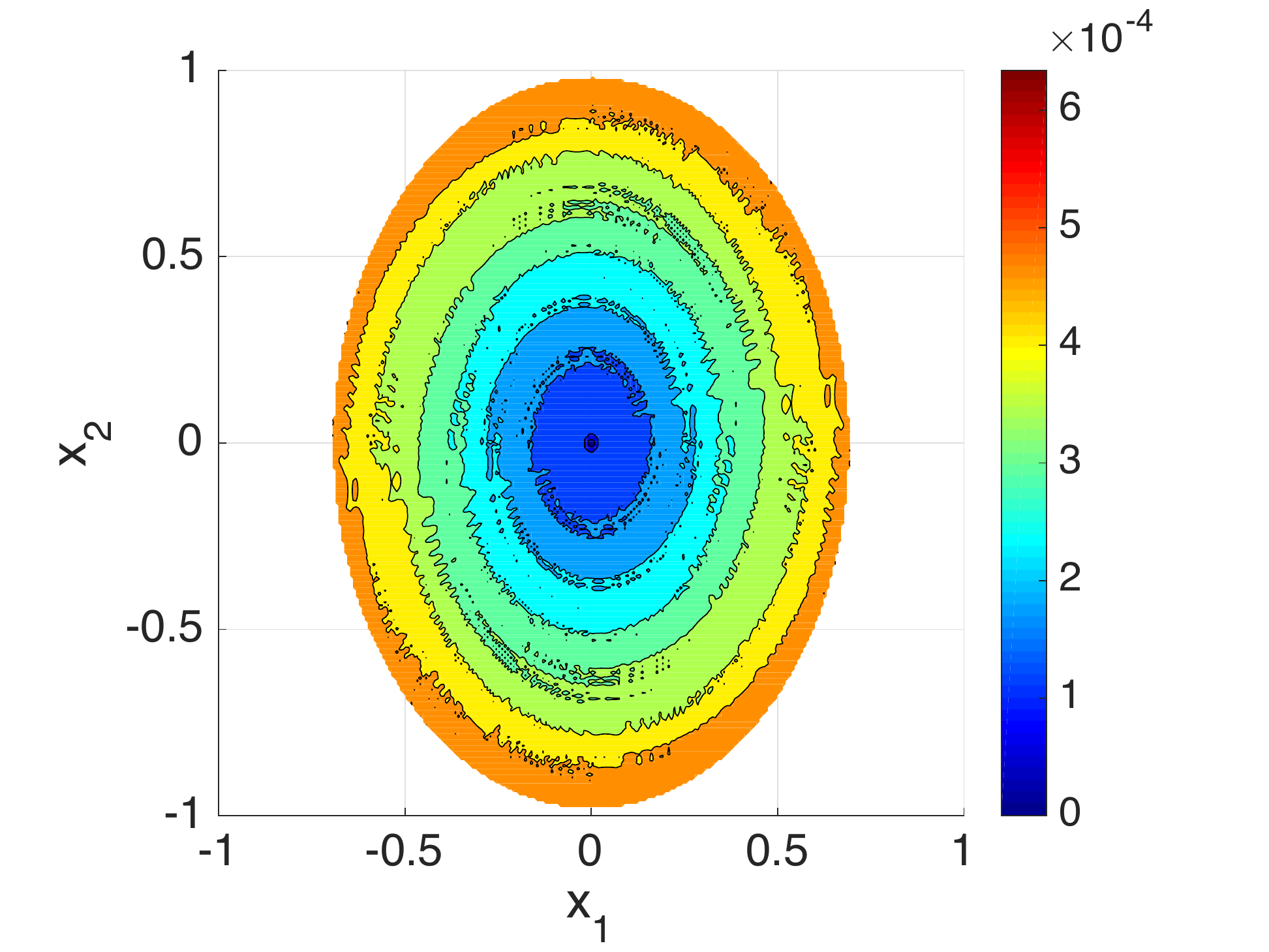}
(b)\includegraphics[width = 0.4\textwidth]{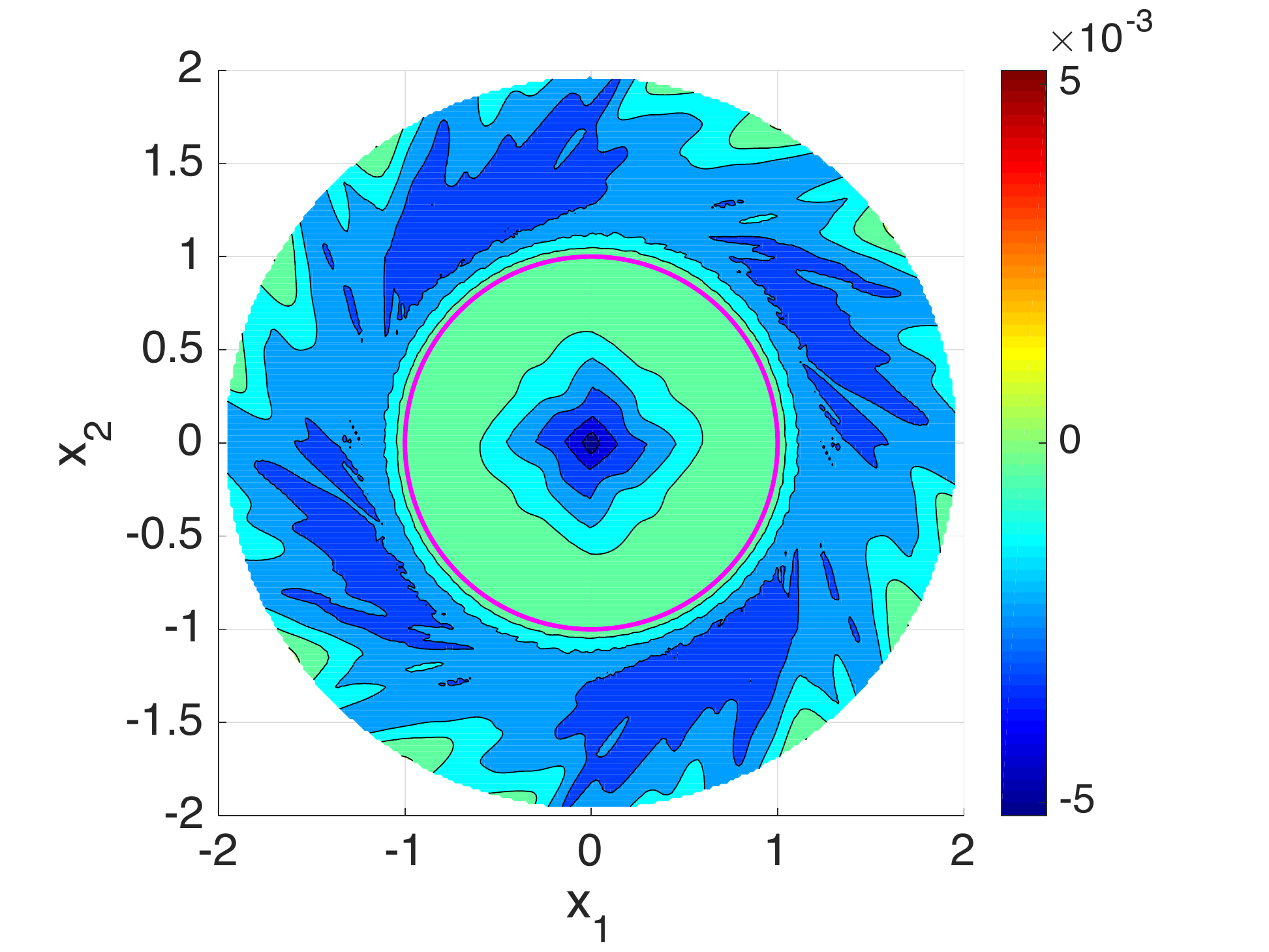}

(c)\includegraphics[width = 0.4\textwidth]{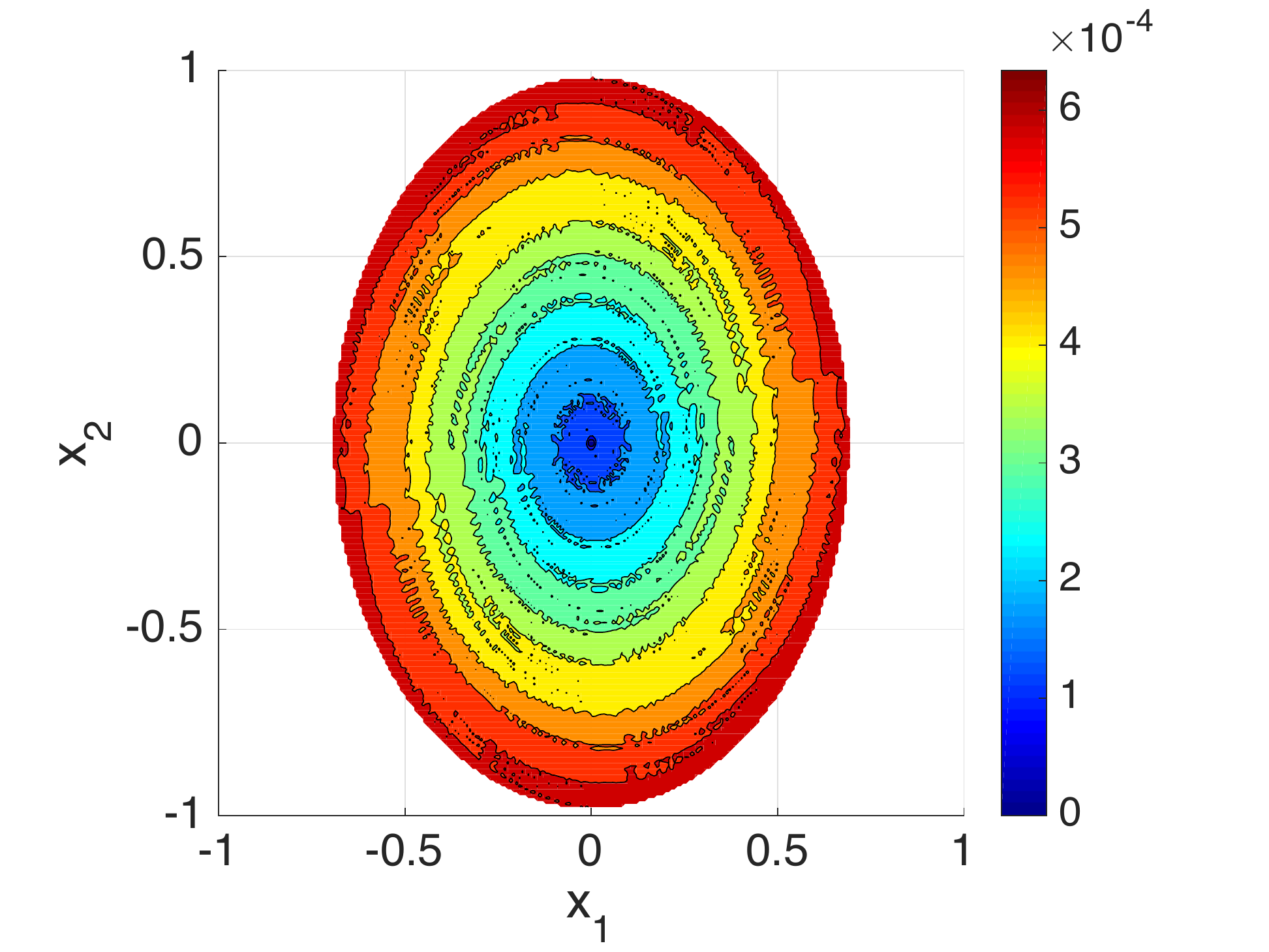}
(d)\includegraphics[width = 0.4\textwidth]{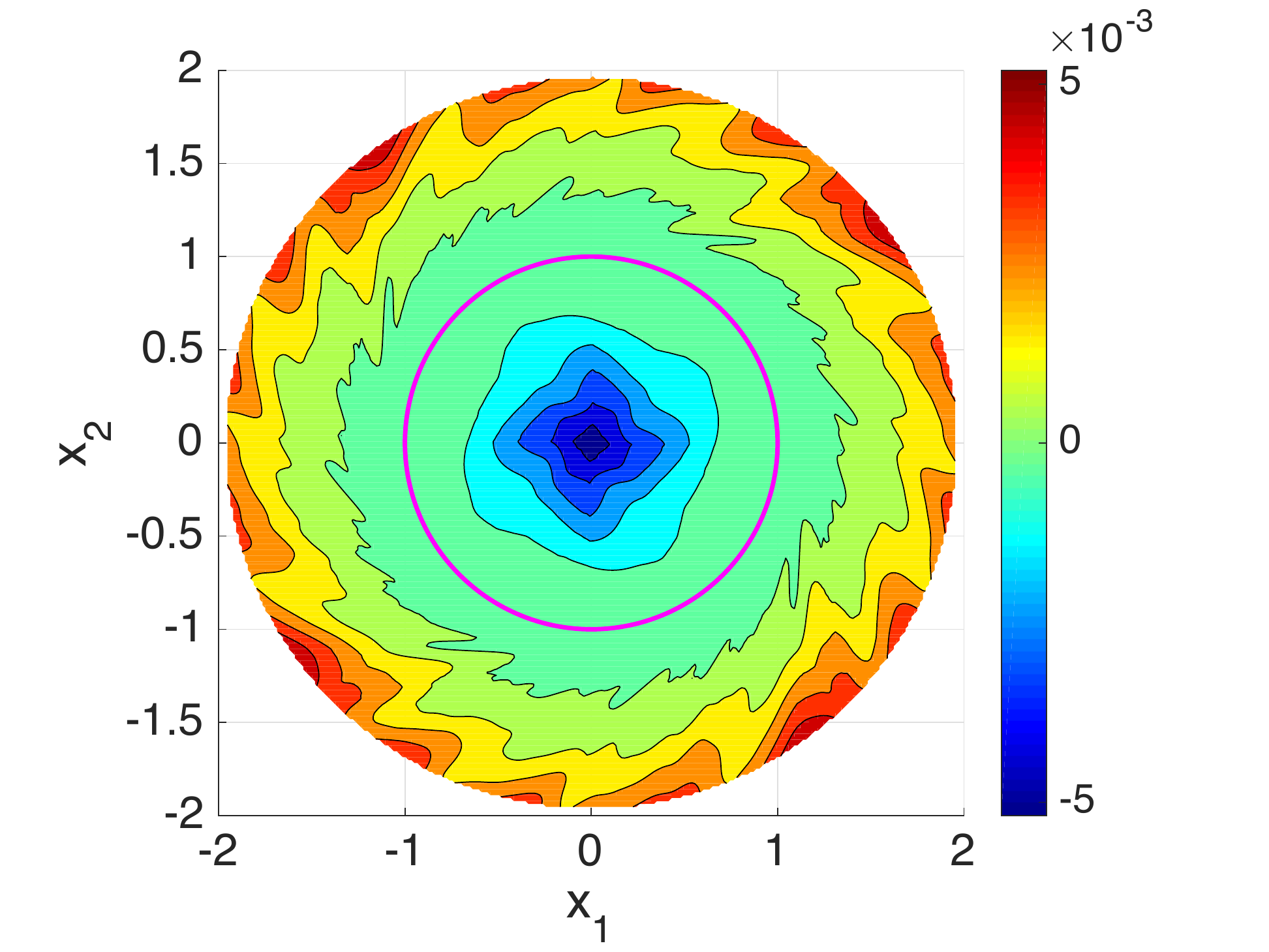}

(e)\includegraphics[width = 0.4\textwidth]{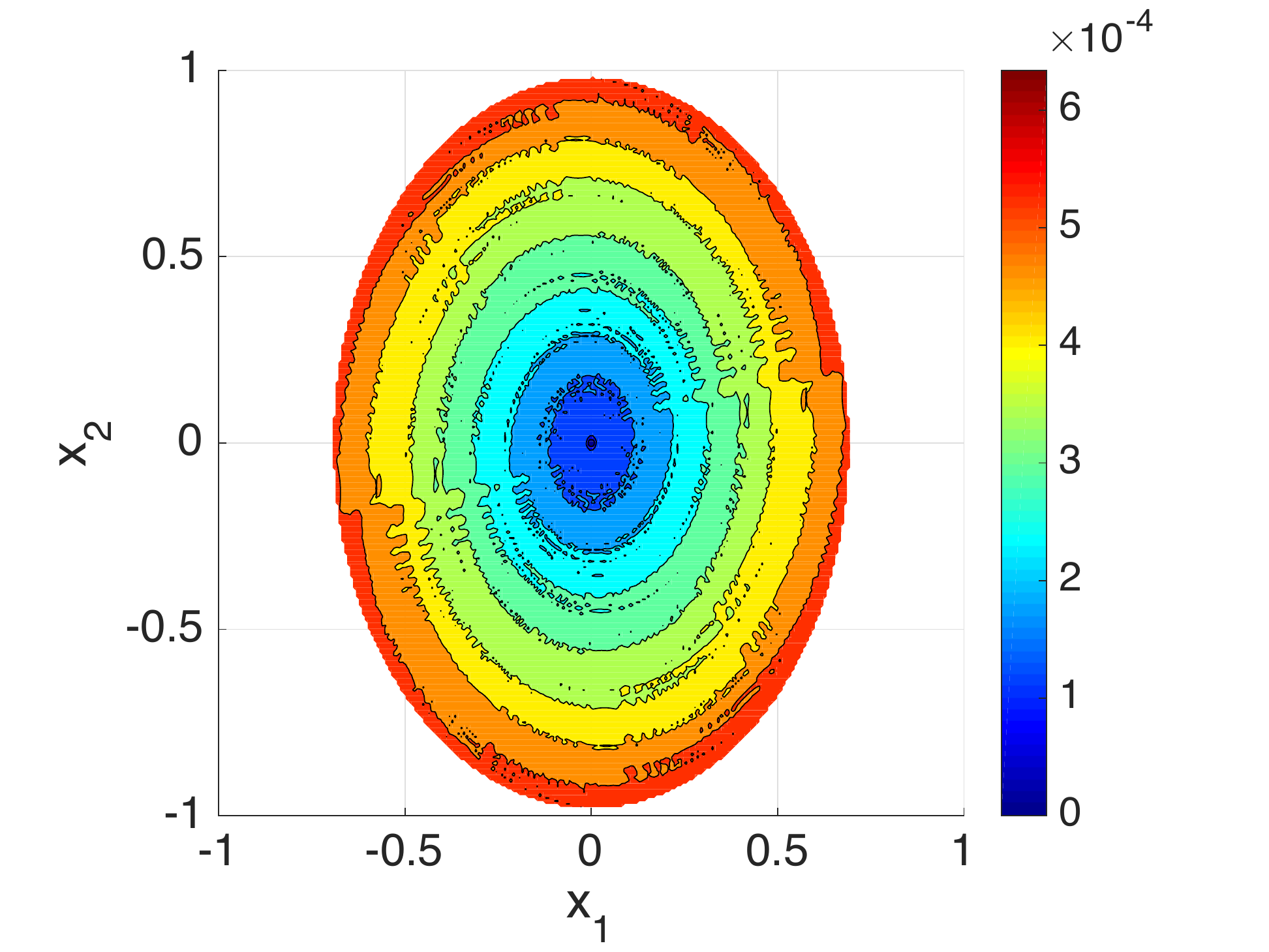}
(f)\includegraphics[width = 0.4\textwidth]{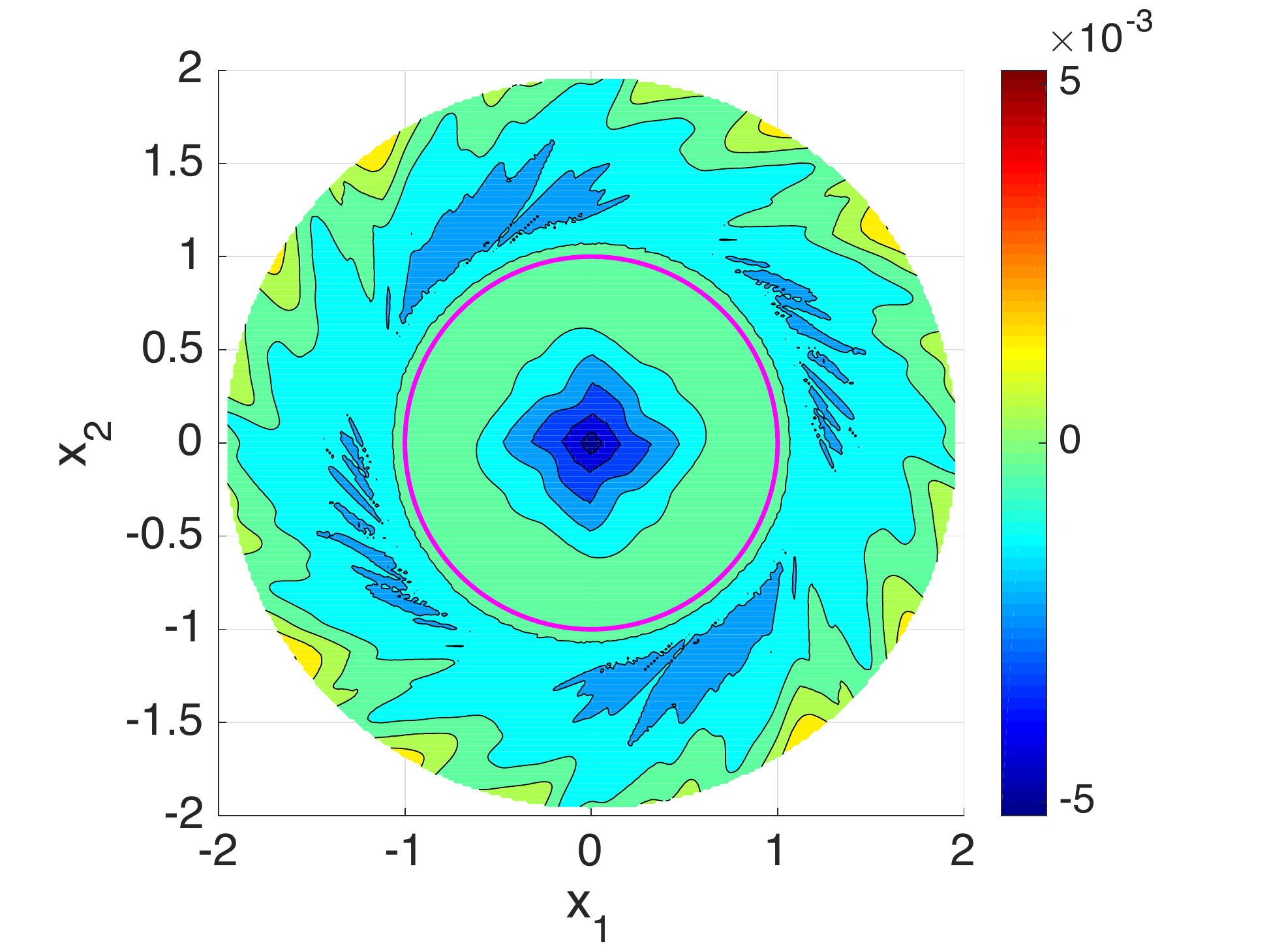}

\caption{Error plots $U_{computed} - U_{exact}$ produced by OLIM-MID (a,b), OLIM-TR(c,d), and OLIM-SIM (e,f) 
on $256\times 256$ meshes with $K=11$.
Left column (a,c,e): SDE \eqref{linSDE}. Right column (b,d,f): SDE \eqref{cirSDE}.
}
\end{center}
\label{fig:errors}
\end{figure*}

%%%%%%%%%%%%%%%%%%%%%%%%%%%%%%%%%%%

%%%%%%%%%%%%%%%%%%%%%%%%%%%%%%%%%%%

\section{Origin of errors in the OLIMs}
\label{sec:errors}
In this section, we discuss various factors which contribute to the error of the numerical solution by the OLIMs. 

%%%
\subsection{Errors of quadrature rules}
\label{seq:quad}
In this Section, we discuss quadrature rule errors applied to the integral along a 
straight line segment connecting the points $\mx_0$ and $\mx_1$
\begin{equation}
\label{line_int}
I(\mx_0,\mx_1):=\int_0^{l}\left( \|\mb(\mx(t))\| - \mb(\mx(t))\cdot \mv\right)dt,
\end{equation}
where 
$$
l:=\|\mx_1-\mx_0\|,~~\mx(t) = \mx_0 + \mv t,~~ \mv:=\frac{\mx_1 - \mx_0}{\|\mx_1 - \mx_0\|}, ~~{\rm i.e.}~~\|\mv\| = 1.
$$
Let $\mb_0$ be the vector field at point $\mx_0$. 
Assuming that $l$ is sufficiently small,
we approximate $\mb(\mx(t))$ by its Taylor expansion around $\mx_0$:
$$
 \mb(\mx(t))=\mb_0+t J \mv  + \tfrac{1}{2}t^2 H + O(t^3),
 $$ 
where 
$$
\mb(\mx): = \left[\begin{array}{c}b_1(\mx)\\b_2(\mx)\end{array}\right],\quad J:= \left[\begin{array}{c}\nabla b_1^T(\mx_0)\\\nabla b_2^T(\mx_0)\end{array}\right],
\quad H = \left[\begin{array}{c}\mv^T\nabla\nabla b_1(\mx_0) \mv\\ \mv^T\nabla\nabla b_2(\mx_0) \mv\end{array}\right].
$$
Then the integral in Eq. \eqref{line_int} becomes
\begin{align}
I(\mx_0,\mx_1)&=\int_0^{l}f(t)dt,~~{\rm where}\notag \\ 
f(t):&= \|\mb(\mx_0+\mv t)\| - \mb(\mx_0+\mv t)\cdot \mv \notag \\
&= \|\mb_0+J \mv t + \tfrac{1}{2}t^2H \|-(\mb_0+J \mv t + \tfrac{1}{2}t^2H)^{T} \mv + O(t^3), \label{line_int1}
\end{align}
Approximating it using the quadrature rule $\mathcal{Q}$ with the error $E_{\mathcal{Q}}$, we obtain
\begin{equation}
\label{line_int2}
I(\mx_0,\mx_1) = \mathcal{Q}(f) + E_{\mathcal{Q}}(f).
\end{equation}
The errors of the right-hand, midpoint and  trapezoid quadrature rules involve $f'$, $f''$, and $f''$ (see e.g. \cite{conte}).
For $f(t)$ given by Eq. \eqref{line_int1}, one can easily evaluate these derivatives at $t=0$. 
The resulting error estimates for the line integrals are
\begin{align}
I(\mx_0,\mx_1) &= \mathcal{R}(f) + \frac{l^2}{2}  \left. \frac{df}{dt}\right\vert_{t=0}  + O(l^3),\label{ER}\\
I(\mx_0,\mx_1) &= \mathcal{M}(f) + \frac{l^3}{24} \left. \frac{d^2f}{dt^2}\right\vert_{t=0} + O(l^4),\label{EM}\\
I(\mx_0,\mx_1) &= \mathcal{T}(f) - \frac{l^3}{12}  \left. \frac{d^2f}{dt^2}\right\vert_{t=0}  + O(l^4),\label{ET}
%I(\mx_0,\mx_1) &= \mathcal{S}(f) - \frac{l^5}{2880} \left. \frac{d^4 f }{dt^4} \right\vert_{t=0} + O(l^3), \label{ES}
\end{align}
where $\mathcal{R}$, $\mathcal{M}$, and  $\mathcal{T}$ denote the right-hand, the midpoint, 
and the trapezoid basic quadrature rules respectively, and 
\begin{align}
 \left. \frac{df}{dt}\right\vert_{t=0} & = \left( \frac{\mb_{0}}{\|\mb_{0}\|}-\mv \right)^{T} J \mv,\label{d1}\\
 \left. \frac{d^2f}{dt^2}\right\vert_{t=0}&=  \frac{\| J \mv \|^{2}  \| \mb_{0} \|^{2} - (\mb_{0}^{T} J \mv)^2}{\|\mb_{0}\|^{3}}
 + \left( \frac{\mb_{0}}{\|\mb_{0}\|}-\mv \right)^{T} H. \label{d2}
%  \left. \frac{d^4 f }{dt^4}  \right\vert_{t=0} & = -3\frac{\|J\mv\|^4 \|\mb_0\|^4-6 \|J\mv\|^2 (\mb_0^T J \mv)^2 \|\mb_0\|^2 + 5 (\mb_0^T J \mv)^4 }{\|\mb_0\|^7} \label{d4}.
\end{align} 
Eqs. \eqref{ER} and \eqref{d1} show that the integration error in OLIM-R is $O(l^2)$ and can be of either sign.
Eqs. \eqref{EM} and \eqref{ET} indicate that  the integration errors in OLIM-MID and OLIM-TR are $O(l^3)$.
The first term of $f''(0)$   (Eq. \eqref{d2}) is due to the linear part of $\mb$. It is nonnegative due to the Schwarz inequality. 
The second term of $f''(0)$ is due to the nonlinearity of $\mb$. It can have an arbitrary sign.
Therefore, the contribution to the error due to the linear part of 
$\mb$ is nonpositive for the midpoint rule, and nonnegative for the trapezoid rule. 
This is consistent with the error plots in Fig. \ref{fig:errors}(a) and (c) for the linear SDE.

Simpson's quadrature rule has an error $O(l^5)$:
\begin{equation}
\label{ES}
I(\mx_0,\mx_1) = \mathcal{S}(f) - \frac{l^5}{2880} \left. \frac{d^4f }{dt^4} \right\vert_{t=0} + O(l^6),
\end{equation}
It is easy to calculate the contribution to $f^{(4)}(0)$ due to the linear part of $\mb$:
\begin{equation}
\label{d4}
\left. \frac{d^4 f }{dt^4}  \right\vert_{t=0}  = -3\frac{\left( \|J\mv\|^2 \|\mb_0\|^2 -  (\mb_0^T J \mv)^2\right)
\left( \|J\mv\|^2 \|\mb_0\|^2 -   5(\mb_0^T J \mv)^2\right) }{\|\mb_0\|^7} +\ldots.
\end{equation}
It can have an arbitrary sign.
It is clear that the integration error in OLIM-SIM is small, and the total error of the numerical solution by OLIM-SIM, 
which is comparable to those by OLIM-MID and OLIM-TR (see Sec. \ref{sec:OLIM2}), is due to the other factors discussed 
in Sections \ref{sec:Kerr} and \ref{sec:curve}.

\begin{figure*}
\begin{center}

(a)\includegraphics[width = 0.4\textwidth]{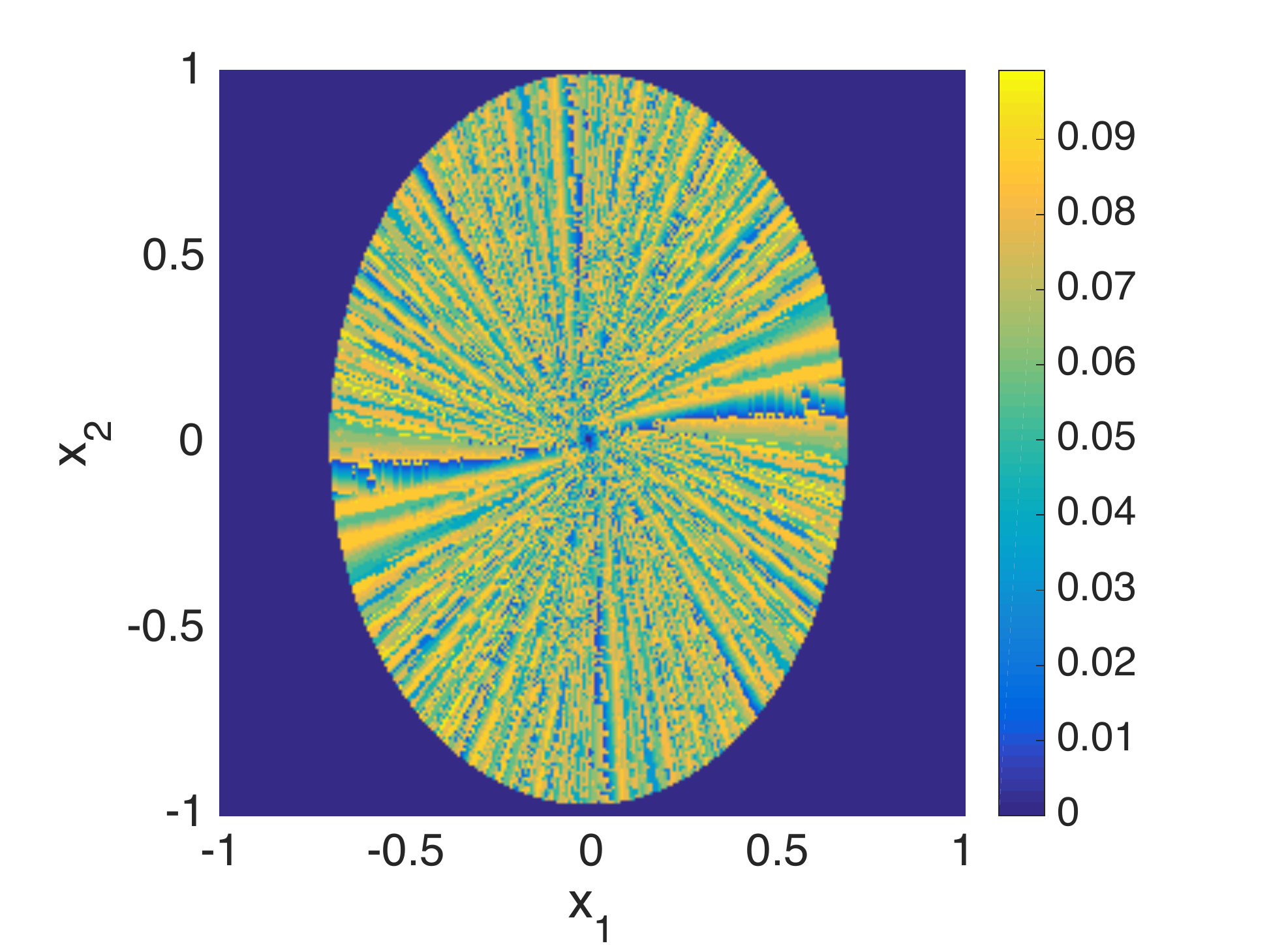}
(b)\includegraphics[width = 0.4\textwidth]{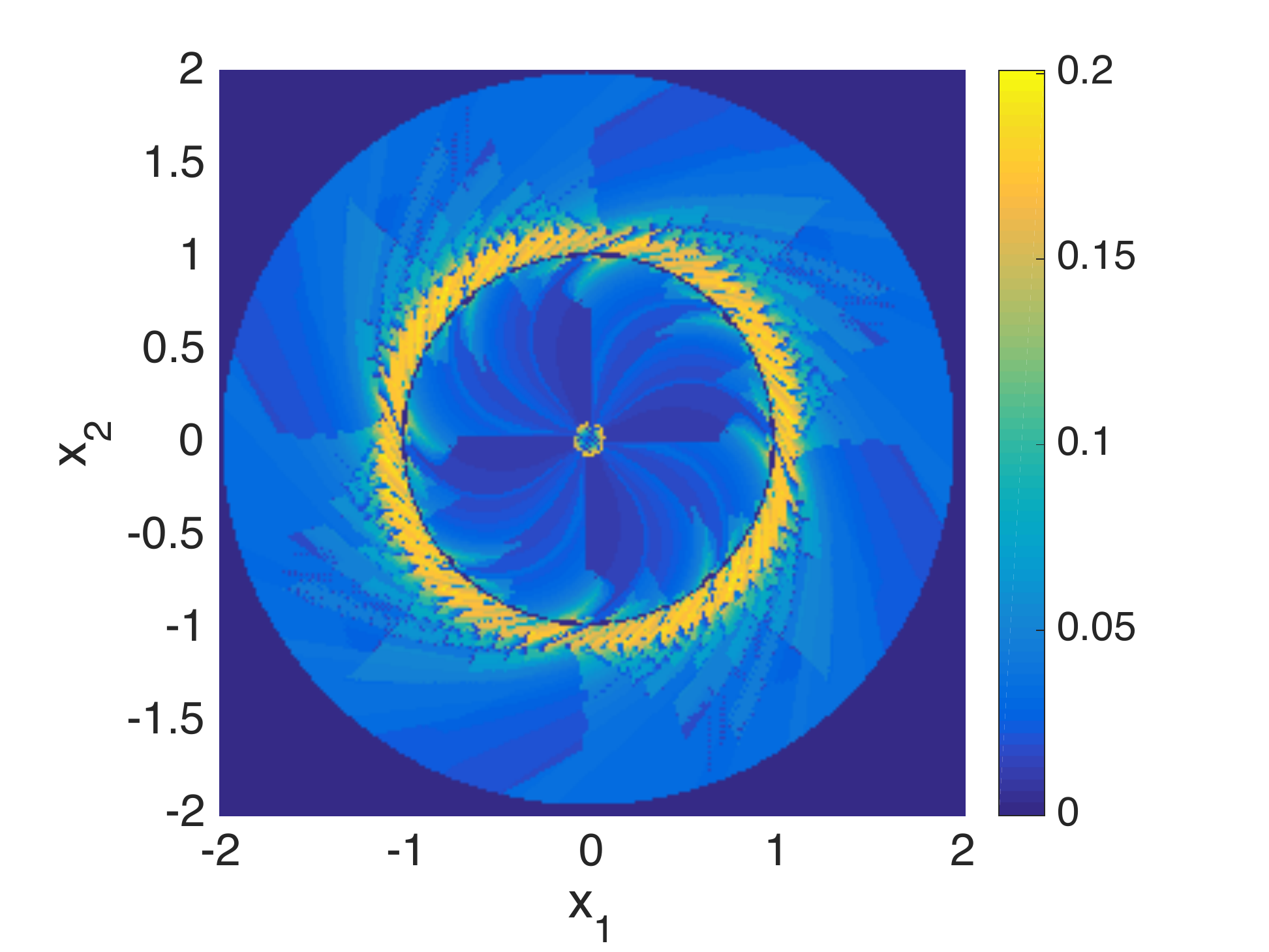}

\caption{
The update length $l$ in the numerical solution by OLIM-MID on $256\times 256$ mesh with $K=11$
for SDE \eqref{linSDE} (a) and SDE \eqref{cirSDE} (b).
}
\end{center}
\label{fig:ulength}
\end{figure*}
{ The lengths $l$ in Eqs. \eqref{ER}-\eqref{ET} and \eqref{ES} can be recorded while the numerical solution is computed.
For each mesh point $\mx$ we define the \emph{update length} as the distance
$$
l(\mx): = \|\mx - \mx_0\|
$$ 
for the one-point update, 
and 
$$l(\mx):=\|\mx - [s^{\ast}\mx_0 + (1-s^{\ast})\mx_1]\|$$ 
for the triangle update, 
where $s^{\ast}$ is the solution of the corresponding minimization problem in Eq. \eqref{2ptOLIM}. 
We update $l(\mx)$ every time when the value of the quasi-potential at $\mx$ is updated.}
The update lengths $l$ in the numerical solutions by OLIM-MID on $256\times 256$ mesh with $K=11$
for SDEs \eqref{linSDE} and \eqref{cirSDE} are shown in Fig. \ref{fig:ulength} (a) and (b) respectively.
The OLIMs are designed so that the update length $l$ can be at most $(K+\sqrt{2})h$ which is $0.097$ 
and $0.194$ for the cases in Figs. \ref{fig:ulength} (a) and (b)
respectively. 
For SDE \eqref{linSDE}, the update length is close to its maximal value at a significant 
fraction of mesh points. For SDE \eqref{cirSDE}, it is close to its maximum in the outer neighborhood of the 
unit circle from which the computation starts. For the other OLIMs, the update lengths are similar.
Therefore, the reduction of the integration error over rather long line segments 
by the use of second and higher order quadrature rules significantly improves the accuracy.
This is exactly what we observe in Figs. \ref{fig:K} and \ref{fig:comparison1}.

 %%%%%%%%%%%%%%%%%%%%%%%%%%%%%%%%%
 
\subsection{Errors due to the curvatures of the level sets of the quasi-potential and MAPs}
\label{sec:curve}

\begin{figure}
\begin{center}

\includegraphics[width = 0.6\textwidth]{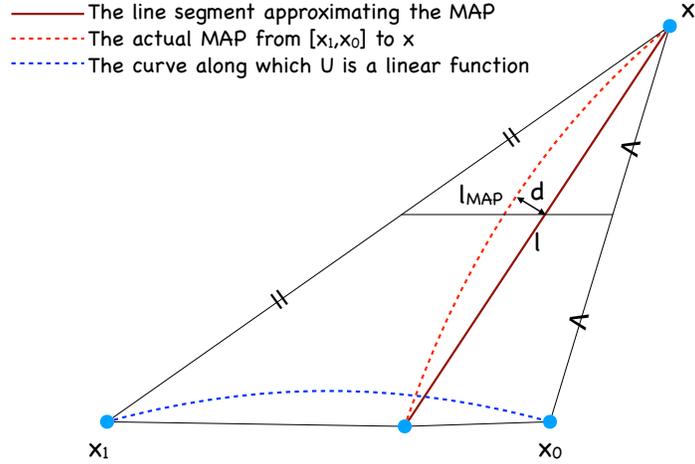}

\caption{ 
Errors due to the curvature of the MAPs and the level sets of the quasi-potential.
}
\end{center}
\label{fig:Ecurv}
\end{figure}

\begin{figure*}
\begin{center}
(a)\includegraphics[width = 0.4\textwidth]{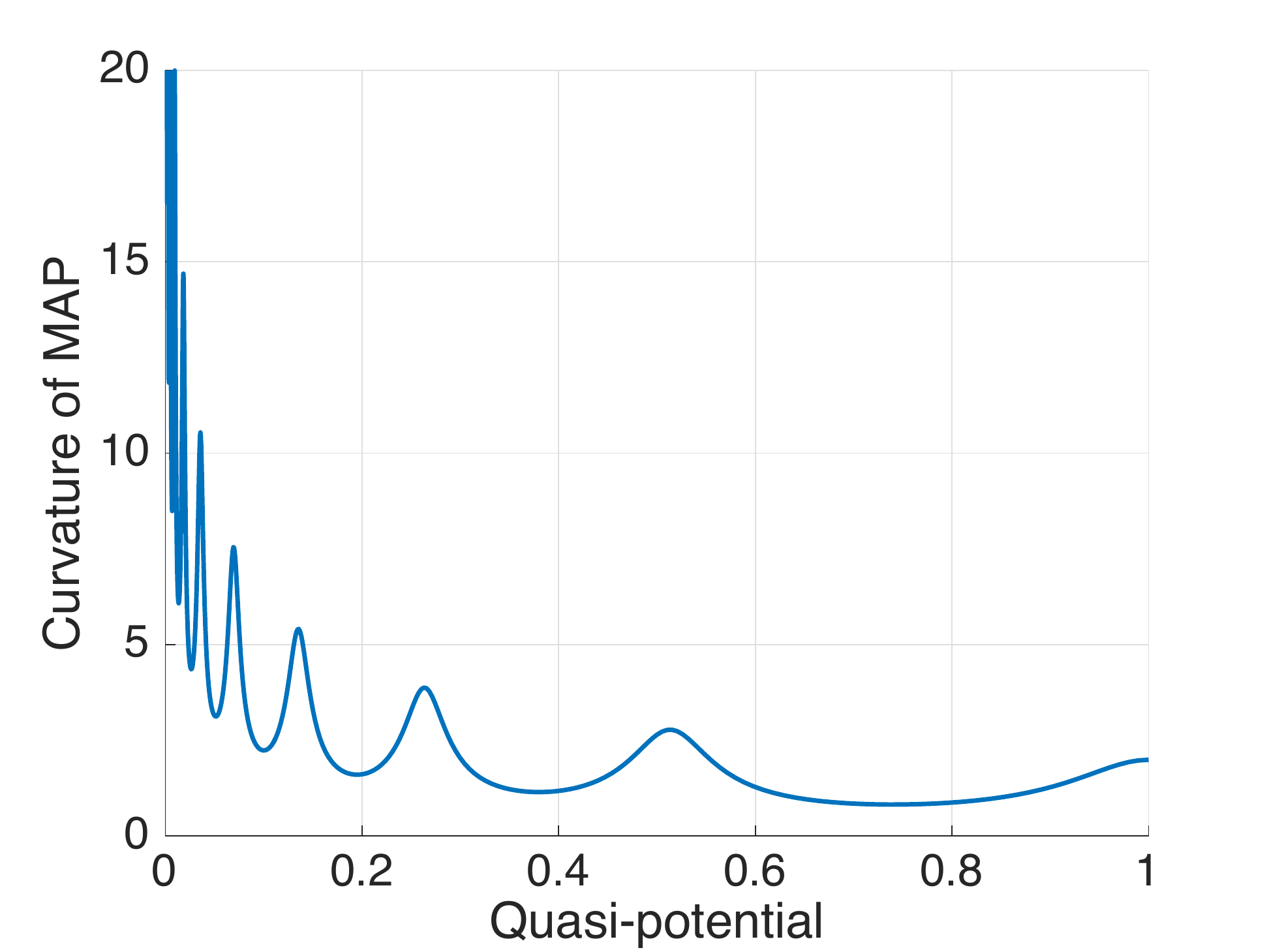}
(b)\includegraphics[width = 0.4\textwidth]{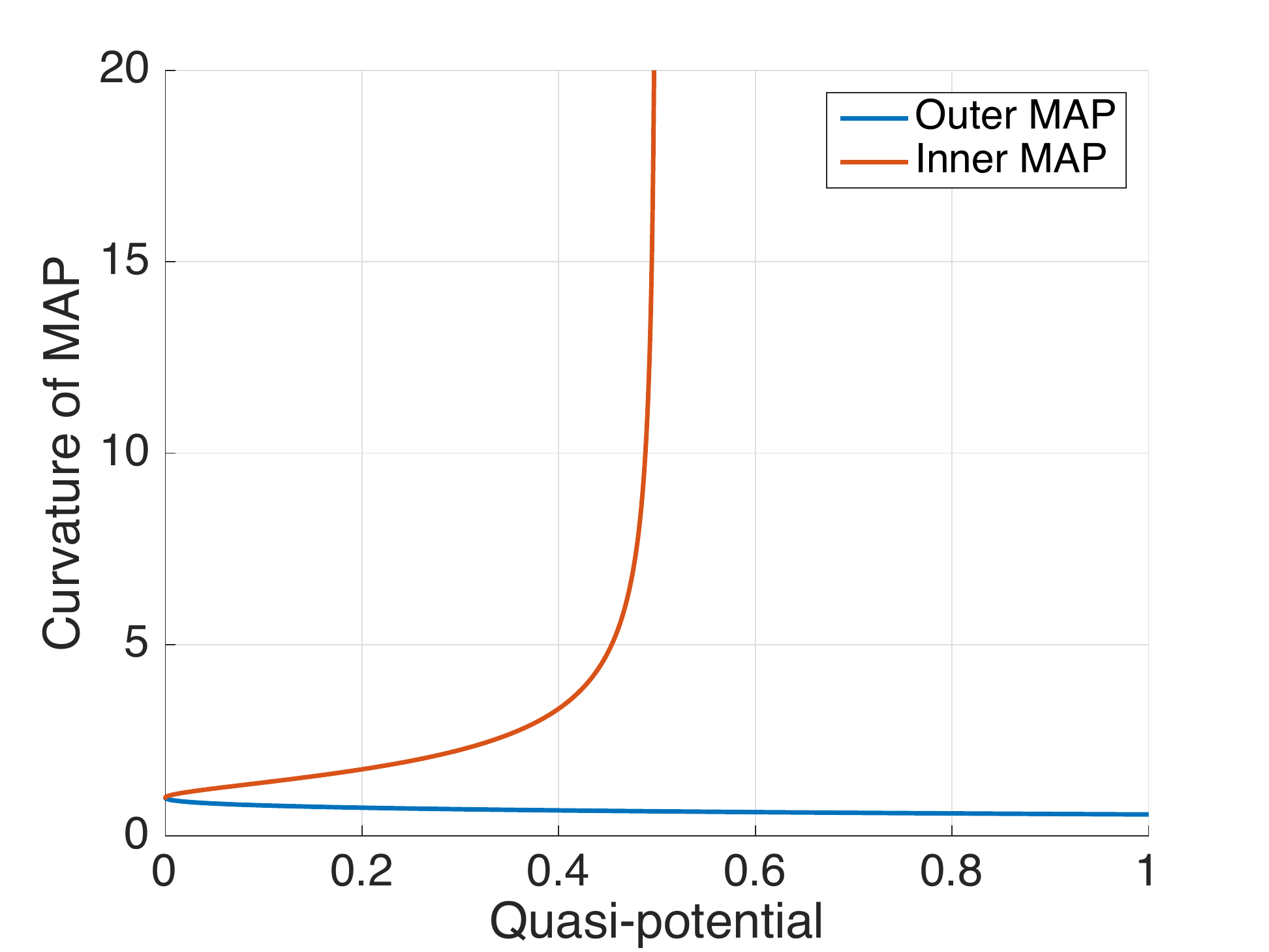}

\caption{
(a): Curvature along a MAP for SDE \eqref{linSDE} plotted versus the quasi-potential along the MAP.
(b): Curvatures along inner and outer MAPs for SDE \eqref{cirSDE} plotted versus the quasi-potential along them.
}
\end{center}
\label{fig:curve}
\end{figure*}

The OLIMs involve two contributions to their numerical errors due to the 
approximations of curved segments with straight line segments (Fig. \ref{fig:Ecurv}).
First, in both, the {one-point } and the {triangle updates}, curved segments of MAPs of lengths at most $Kh(1+o(1))$ 
and  $(Kh + \sqrt{h_1^2 + h_2^2})(1 + o(1))$ respectively
are approximated by straight line segments. Second, 
in the {triangle update}, curved segments of length at most $\sqrt{h_1^2 + h_2^2}(1+o(1))$ 
along which the quasi-potential changes linearly, are approximated by straight line segments. 

Let $l$ be the length of the  line segment in the {one-point update} or in the {triangle update}
approximating the MAP segment. Suppose that there is a MAP connecting the endpoints of this line segment, and its
curvature is $\kappa > 0$. Assume that $l$ is small and $\kappa$ is constant. Then the length of the MAP segment is
\begin{equation}
\label{lenMAP}
l_{MAP}  = \frac{2}{\kappa}\arcsin\left(\frac{\kappa l}{2}\right)= l\left(1 + \frac{(\kappa l)^2}{24}+ O((\kappa l)^4)\right).
\end{equation}
OLIM-MID and OLIM-SIM involve the evaluation of $\mb$ at the midpoint of the line segment. 
The midpoint of the MAP segment is at distance $d$ from it:
\begin{equation}
\label{bmid}
d = \frac{1}{\kappa}\left(1-\sqrt{1-\frac{(\kappa l)^2}{4}}\right) = l\left(\frac{\kappa l}{8} + O((\kappa l)^3)\right).
\end{equation}
The norm of the difference of $\mb$ evaluated at these two midpoints is proportional to $d$. 
The curvature of MAPs can blow up near the equilibria of $\mb$. As an example, the graphs of the curvature of the MAPs
versus the quasi-potential along them for SDEs \eqref{linSDE} and \eqref{cirSDE} are plotted in Fig. \ref{fig:curve}.

The bases of the triangles $[\mx_1,\mx_0]$ in the {triangle update} are at most of length $h\sqrt{2}$. In the {triangle update},
the quasi-potential along the line segment $[\mx_1,\mx_0]$ is assumed to change linearly, while the actual curve along which the quasi-potential 
changes linearly, is not  a straight line in general. For example, if $U(\mx_1)=U(\mx_0)$, such a curve is the corresponding level set. 
For SDEs \eqref{linSDE} and \eqref{cirSDE}, such level sets are ellipses and circles respectively.  Let $\mu>0$ be the curvature of the 
curve connecting $\mx_1$ and $\mx_0$ along which the quasi-potential is a linear function. 
Then the length of the segment of MAP connecting the point $\mx$ to be updated  
with this curve is shorter (as in Fig. \ref{fig:Ecurv}) or longer than $l_{MAP}$ in Eq. \eqref{lenMAP}
by at most 
\begin{equation}
\label{dlen}
\mu h^2 + hO((\mu h)^3).
\end{equation}

In summary, the local error of the {triangle update} due to approximating curved segments with line segments is bounded by
\begin{equation}
\label{ecurve}
E_{curve}\lesssim \frac{\left(Kh+\sqrt{h_1^2 + h_2^2}\right)^3\kappa^2}{24} + \|J\|\frac{\left(Kh+\sqrt{h_1^2 + h_2^2}\right)^2\kappa }{8} + \mu h^2,
\end{equation}
where $J$ is the Jacobian matrix of $\mb$. 
Eq. \eqref{ecurve} explains relatively large numerical errors near the origin by OLIM-MID, OLIM-TR, and OLIM-SIM applied to SDE \eqref{cirSDE}
(see Fig. \ref{fig:errors} (b), (d), and (f)). 

The comparison of $E_{curve}$ and the integration errors of the OLIMs (Eqs. \eqref{ER}, \eqref{EM}, \eqref{ET}, and \eqref{ES}),
combined with the results of our numerical tests in Section \ref{sec:tests}, show that 
the OLIMs are at most first order accurate regardless of the quadrature rule used. However, for a wide range of reasonable mesh sizes,
OLIM-MID, OLIM-TR, and OLIM-SIM might appear to have superlinear convergence (see the least squares fits in Table \ref{Table:LSfit2})
due to the change of the relative magnitudes 
of error terms of different orders. 
On the contrary, the OLIM-R exhibits sublinear convergence (see the least squares fits in Table \ref{Table:LSfit}). 

%%%
\subsection{Finite update radius}
\label{sec:Kerr}
The direct application of the OUM developed in \cite{oum1,oum2} to the Hamilton-Jacobi PDE 
\eqref{HJ} for the quasi-potential would require an infinite update radius.
The OUM adjusted for computing the quasi-potential \cite{quasi} as 
well as the OLIMs set the update radius to $Kh$ where the update factor $K$ is a finite positive integer, 
and $h=\max\{h_1,h_2\}$ is the mesh step. 
The extra error due to the finite update radius may appear in the case  
where the angle between the normal to the level set of $U$ approximating the {\sf Accepted Front} and the vector field $\mb$ at the point $\mx$ to be
updated is close to $\pi/2$ (the angle $\alpha$ in Fig. \ref{fig:rat} (a, Top)).
 This extra error was quantified in \cite{quasi} and shown not to exceed
$0.5\|\mb(\mx)\|h(K^{-1} + O(K^{-3}))$ for sufficiently small $h$. 

{ To minimize this extra error, one might be tempted to use large values of the update factor $K$.
However, there is a trade-off. 
The quadrature errors are approximately proportional to 
$h^pK^p$ where $p = 2,3,$ or  5 depending on the quadrature rule used (see Section \ref{seq:quad}), and the errors due to 
the curvature of MAPs and level sets of the quasi-potential are approximately $C_1h^2K^2 + C_2h^3K^3$.
Hence, the optimal $K$ is the minimizer of the  local error function of the form
$$
E_h(K) = A\frac{h}{K} + Bh^2K^2 + Ch^3K^3 + Dh^2,
$$
where the coefficients $A$, $B$, $C$ and $D$ depend on the vector field $\mb$ and its derivatives, on the direction and the
curvature of the MAPs, and  
on the curvature of the
 level sets of the quasi-potential. They change from one mesh point to another, and can have different orders of magnitude.
Differentiating $E_h$ with respect to $K$ we get the following equation for the optimal  $K$:
$$
\frac{dE_h(K)}{dK} = -\frac{Ah}{K^2} + 2Bh^2K+3Ch^3K^2 = \frac{h}{K^2}\left[-A + 2BhK^3 + 3Ch^2K^4\right] = 0.
$$
In two simple cases, where $B\ll Ch$ or $Ch \ll B$, the optimal $K$ is $O(h^{-1/3})$ or $O(h^{-1/2})$ respectively, i.e.,
$O(N^{1/3})$ or $O(N^{1/2})$. In the general case, the dependence  of the optimal $K$ on $N$ is more complicated,
however, it is clear that $K$ must grow with $N$ but not faster that $N^{1/2}$. 
Our numerical experiments suggest simple Rules-of-Thumb for choosing $K$ (see Section \ref{sec:Ktest}).
It turns out that $K$ does not need to be very large.}

%%%%%%%%

%%%%%%%%%%%%%%%%%%%%%%%%%%%%%
\section{Conclusion}
\label{sec:conclusion}
We have introduced the family of Ordered Line Integral Methods (OLIMs) for computing the quasi-potential 
on a rectangular mesh. The update rules in OLIM-R, employing the right-hand rectangle quadrature rule, are equivalent to those in the OUM.
OLIM-MID, OLIM-TR, and OLIM-SIM (employing the midpoint, the trapezoid, and Simpson's quadrature rules) converge faster and 
admit errors  two to three orders of magnitude smaller than  OLIM-R  and the OUM (See Fig. \ref{fig:K}). 
Nevertheless, asymptotically, the OLIMs, like the OUM, are at most
first order accurate due to the use of linear interpolation. 
While the use of second or higher order quadrature rules requires the use of a nonlinear solver and hence makes the {triangle update}
more expensive { than} the one in the OUM, the proposed hierarchical update strategy for reducing the number { of} calls for the {triangle update} makes the OLIMs
faster than the OUM. Due to this strategy,  OLIM-R is about four  times faster than the OUM.
OLIM-MID and OLIM-TR are about 1.5 times faster than the OUM for the optimal $K$'s and about three times faster for $K=50$
(see Fig. \ref{fig:CPUvsK}). Our C codes
{\tt OLIM\_righthand.c}, {\tt OLIM\_midpoint.c}, {\tt OLIM\_trapezoid.c} and {\tt OLIM\_simpson.c}
implementing the corresponding OLIMs are available on M. Cameron's { website} \cite{mariakc}.

We have investigated the dependence of the error upon the update factor $K$ and 
concluded that, while large $K$ nearly eliminates 
the error in the direction of the MAP due to insufficient update radius in the case of a slowly changing vector field $\mb$, 
it increases the integration error and the error
due to the curvature of the MAP. 
Based on our study of the relationship between the numerical error and the update factor $K$,
we have proposed Rules-of-Thumb for choosing $K$. 

Our comparison of the OLIMs with second or higher order quadrature rules shows that the best one 
of them in terms of the balance between the accuracy and the CPU time
is achieved by OLIM-MID.

The present work is focused on the 2D case with isotropic diffusion. 
{ The OLIMs can be extended to the case with anisotropic diffusion and to 3D. 
We will report these developments in the future.}

%%%%%%%%%%%%%%%%%%%%%%%%%%%%%%%%%

\section*{Acknowledgements}
We thank Professor A. Vladimirsky for a valuable discussion.
This work was supported in part by the NSF grant DMS1554907.

%%%%%%%%%%%%%%%%%%%%%%%%%%%%%
%
%%%%%%%%%   A P P E N D I X   A
%
%%%%%%%%%%%%%%%%%%%%%%%%%%%%%

 \appendix
\setcounter{equation}{0}
\renewcommand{\theequation}{A-\arabic{equation}}

{
 \section*{Appendix A. The Freidlin-Wentzell action vs the Geometric action}
 \label{sec:AppA}
The Freidlin-Wentzell action functional for SDE \eqref{sde1} is defined on the set of absolutely continuous paths $\phi(t)$ by \cite{FW}
\begin{equation}
\label{FWA}
S_T(\phi) = \frac{1}{2}\int_0^T\|\dot{\phi} - \mb(\phi)\|^2dt.
\end{equation}
The original definition of the quasi-potential \cite{FW} with respect to a compact set $A$ (an attractor of $\dot{\mx}=\mb(\mx)$) at a point $\mx$ is 
\begin{equation}
\label{odef}
U_A(\mx) = \inf_{T,\phi}\left\{S_T(\phi)~|~\phi(0)\in A,~\phi(T)=\mx,~\phi~\text{is absolutely continuous}\right\}.
\end{equation}
The minimization with respect to the travel-time $T$ can be performed analytically \cite{FW,hey1,hey2} resulting at the geometric action $S(\psi)$.
Let  $\phi(t)$ be  a fixed absolutely continuous path $\phi(t)$. Expanding $\|\cdot\|^2$ in Eq. \eqref{FWA} and using the inequality
$y^2 + z^2 \ge 2yz$ for all nonnegative real numbers $y$ and $z$, we get:
\begin{align}
S_T(\phi) &= \frac{1}{2}\int_0^T\|\dot{\phi} - \mb(\phi)\|^2dt = \frac{1}{2}\int_0^T\left(\|\dot{\phi}\|^2 - 2\dot{\phi}\cdot\mb(\phi) + \|\mb(\phi)\|^2\right)dt \notag \\
&\ge  \frac{1}{2}\int_0^T\left(2\|\dot{\phi}\| \|\mb(\phi)\|- 2\dot{\phi}\cdot\mb(\phi) \right)dt \label{derg} \\
& = \int_0^T\left(\|\dot{\phi}\| \|\mb(\phi)\|- \dot{\phi}\cdot\mb(\phi) \right)dt.\notag
\end{align}
The inequality in Eq. \eqref{derg} becomes an equality if and only if $\|\dot{\phi}\| = \|\mb(\phi)\|$.
Let $\chi$ be the path obtained from $\phi$ by a reparametrization such that $\|\dot{\chi}\| = \|\mb(\chi)\|$. Then
\begin{equation}
\label{derg2}
S_T(\phi)\ge S_{T_{\chi}}(\chi) =  \int_0^{T_{\chi}}\left(\|\dot{\chi}\| \|\mb(\chi)\|- \dot{\chi}\cdot\mb(\chi) \right)dt.
\end{equation}
Note that $T_{\chi}$ can be infinite. The integral in right-hand side of Eq. \eqref{derg2} is invariant with respect to the parametrization of the path $\chi$.
Hence, we can pick the most convenient one, for example, the arclength parametrization, and denote the reparametrized path by $\psi$. Hence,
\begin{equation}
\label{derg2}
S_{T_{\chi}}(\chi) =  \int_0^L\left(\|\psi_s(s)\| \|\mb(\psi(s))\|- \psi_s(s)\cdot\mb(\psi(s)) \right)ds =:S(\psi),
\end{equation}
where $L$ is the length of the paths $\chi$ and  $\psi$ (corresponding to the same curve).
For computation of the quasi-potential, it is more convenient to deal with the geometric action $S(\psi)$ than with the Freidlin-Wentzell action $S_T(\phi)$.
}

%%%%%%%%%%%%%%%%%%%%%%%%%%%%%
%
%%%%%%%%%   A P P E N D I X   B
%
%%%%%%%%%%%%%%%%%%%%%%%%%%%%%

 \appendix
\setcounter{equation}{0}
\renewcommand{\theequation}{B-\arabic{equation}}

 \section*{Appendix B. The {triangle updates} for the OLIMs}
 \label{sec:AppB}

\textbf{OLIM-R}\\
OLIM-R performs the {triangle update} by solving the following minimization problem
\begin{align}
&u = \min_{s\in[0,1]}\left[su_0 + (1-s)u_1 + \|\mb\|\|\mx - \mx_s\| - \mb\cdot(\mx - \mx_s)\right],
\label{r_olim1}\\
&\text{where}~~\mb\equiv\mb(\mx),~~\mx_s = s\mx_0 + (1-s)\mx_1,~~u_0 \equiv U(\mx_0),~~u_1 \equiv U(\mx_1) .\notag
\end{align}
Taking the derivative of the function to be minimized
$$
f(s): =  s u_0 + (1-s) u_1 + \|\mb\|  \|\mx - s\mx_0 - (1-s) \mx_1\| - \mb \cdot (\mx - s\mx_0 - (1-s) \mx_1)
$$
with respect to $s$ and setting it to zero, we obtain the following equation for $s$:
\begin{equation}
\label{OR0}
u_0-u_1 + \|\mb(\mx)\| \frac{ (\mx-\mx_s)\cdot(\mx_1-\mx_0) }{ \|\mx - \mx_s\| } - \mb(\mx) \cdot (\mx _1- \mx_0 )=0.
\end{equation}
Regrouping terms and taking squares we obtain the following quadratic equation for $s$:
\begin{align}
&As^2 + 2Bs + C = 0,~~{\rm where} \label{OR}\\
A& = \|\mx _1- \mx_0 \|^2 \left([\mb(\mx) \cdot (\mx _1- \mx_0 ) - (u_0-u_1)]^2 -  \|\mb(\mx)\|^2 \|\mx _1- \mx_0 \|^2\right),\\
B& =   \left([\mb(\mx) \cdot (\mx _1- \mx_0 ) - (u_0-u_1)]^2 -  \|\mb(\mx)\|^2 \|\mx _1- \mx_0 \|^2\right)\left[(\mx-\mx_1)\cdot(\mx_1-\mx_0)\right],\\
C& =  [\mb(\mx) \cdot (\mx _1- \mx_0 ) - (u_0-u_1)]^2\|(\mx-\mx_1)\|^2 -  \|\mb(\mx)\|^2 \left((\mx-\mx_1)\cdot(\mx_1-\mx_0)\right).
\end{align}
We solve Eq. \eqref{OR}, select its root $s^{\ast}$, if any, on the interval $[0,1]$, and verify that it is also the root of Eq. \eqref{OR0}.
In the case of success, the {triangle update} returns
$$
{\sf Q}_{\Delta}(\mx_1,\mx_0,\mx) = s^{\ast} u_0 + (1-s^{\ast}) u_1 + \|\mb\|  \|\mx - s^{\ast}\mx_0 - (1-s^{\ast}) \mx_1\| - \mb \cdot (\mx - s^{\ast}\mx_0 - (1-s^{\ast}) \mx_1).
$$
Otherwise, it returns ${\sf Q}_{\Delta}(\mx_1,\mx_0,\mx) =+\infty$.

%%%%%%%%%%
\textbf{OLIM-MID}\\
OLIM-MID performs the {triangle update} by solving the following minimization problem
\begin{align}
u &= \min_{s\in[0,1]}\left[su_0 + (1-s)u_1 + \|\mb_{ms}\|\|\mx - \mx_s\| - \mb_{ms}\cdot(\mx - \mx_s)\right],~~\text{where}
\label{mid_olim}\\
\mx_s& = s\mx_0 + (1-s)\mx_1,~~\mb\equiv \mb(\mx)\notag\\
\mb_{ms}& = s\mb_{m0} + (1-s)\mb_{m1},~~ { \mb_{m0}}\equiv  \mb\left(\frac{\mx_0+\mx}{2}\right),~~\mb_{m1}\equiv \mb\left(\frac{\mx_1+\mx}{2}\right) .\notag
\end{align}
Taking the derivative of 
$$
f(s): =  su_0 + (1-s)u_1 + \|\mb_{ms}\|\|\mx - \mx_s\| - \mb_{ms}\cdot(\mx - \mx_s)
$$
with respect to $s$ and setting it to zero, we obtain the following equation for $s$:
\begin{align}
& u_0-u_1 + \|\mb_{ms}\| \frac{ (\mx-\mx_s)\cdot(\mx_1-\mx_0) }{ \|\mx - \mx_s\| }  +  \|\mx - \mx_s\|\frac{\mb_{ms}\cdot (\mb_{m0}-\mb_{m1}) }{\|\mb_{ms}\|} \notag \\
& - \mb_{ms} \cdot (\mx _1- \mx_0 ) - (\mx-\mx_s)\cdot (\mb_{m0}-\mb_{m1}) =0. \label{OM0}
\end{align}
%where $\mb_0\equiv \mb(\mx_0)$ and $\mb_1\equiv\mb(\mx_1)$.
The hybrid nonlinear solver \cite{wilkinson,stewart} is used for finding a root $s^{\ast}$ of Eq. \eqref{OM0} in the interval $[0,1]$.
In the case of success, the {triangle update} returns
$$
{\sf Q}_{\Delta}(\mx_1,\mx_0,\mx) = s^{\ast} u_0 + (1-s^{\ast}) u_1 + 
\|\mb_{ms^{\ast}}\|  \|\mx - s^{\ast}\mx_0 - (1-s^{\ast}) \mx_1\| - \mb_{ms^{\ast}} \cdot (\mx - s^{\ast}\mx_0 - (1-s^{\ast}) \mx_1).
$$
Otherwise, it returns ${\sf Q}_{\Delta}(\mx_1,\mx_0,\mx) =+\infty$.

%%%%%%%%%%
\textbf{OLIM-TR}\\
OLIM-TR performs the {triangle update} by solving the following minimization problem
\begin{align}
u &= \min_{s\in[0,1]}\left[su_0 + (1-s)u_1 +\frac{1}{2}\left\{ (\|\mb_s\| +\|\mb\|)\| \mx - \mx_s\| - (\mb_s +\mb)\cdot(\mx - \mx_s)\right\}\right], \label{tr_olim}\\
&\text{where}\notag \\
\mx_s& = s\mx_0 + (1-s)\mx_1\notag\\
\mb_s& = s\mb_0+ (1-s)\mb_1,\quad \mb_0\equiv \mb(\mx_0),~~\mb_1\equiv\mb(\mx_1),~~\mb\equiv \mb(\mx).\notag
\end{align}
Taking the derivative of 
$$
f(s): =  su_0 + (1-s)u_1 +\frac{1}{2}\left\{ (\|\mb_{s}\| +\|\mb\|)\| \mx - \mx_s\| - (\mb_{s} +\mb)\cdot(\mx - \mx_s)\right\}
$$
with respect to $s$ and setting it to zero we obtain the following equation for $s$:
\begin{align}
& u_0-u_1 + \frac{1}{2}\{ (\|\mb_s\| +\|\mb\|) \frac{ (\mx-\mx_s)\cdot(\mx_1-\mx_0) }{ \|\mx - \mx_s\| }  +  \|\mx - \mx_s\| \frac{\mb_s\cdot (\mb_0-\mb_1) }{\|\mb_s\|} \notag \\
& - (\mb_s +\mb) \cdot (\mx _1- \mx_0 ) - (\mx-\mx_s)\cdot (\mb_0-\mb_1)  \}  =0. \label{OT0}
\end{align}
The hybrid nonlinear solver \cite{wilkinson,stewart} is used for finding a root $s^{\ast}$ of Eq. \eqref{OT0} in the interval $[0,1]$.
In the case of success, the {triangle update} returns
$$
{\sf Q}_{\Delta}(\mx_1,\mx_0,\mx) = s^{\ast} u_0 + (1-s^{\ast}) u_1 + \frac{1}{2}\left\{ (\|\mb_{s^{\ast}} \| +\|\mb\|)\| \mx - \mx_{s^{\ast}}\| - (\mb_{s^{\ast}} +\mb)\cdot(\mx - \mx_{s^{\ast}})\right\}.
$$
Otherwise, it returns ${\sf Q}_{\Delta}(\mx_1,\mx_0,\mx) =+\infty$.

%%%%%%%%%%
\textbf{OLIM-SIM}\\
OLIM-SIM performs the {triangle update} by solving the following minimization problem
\begin{align}
u &= \min_{s\in[0,1]} [su_0 + (1-s)u_1 +\frac{1}{6} \{ (\|\mb_s\| +4\|\mb_{ms}\|+\|\mb\|)\| \mx - \mx_s\| \notag \\
& - (\mb_s +4\mb_{ms}+\mb)\cdot(\mx - \mx_s) \} ], \label{sim_olim}\\
&\text{where}\notag \\
\mx_s& = s\mx_0 + (1-s)\mx_1,~~\mb\equiv\mb(\mx)  \notag\\
\mb_s& = s\mb_0+ (1-s)\mb_1,\quad \mb_0\equiv \mb(\mx_0),~~\mb_1\equiv\mb(\mx_1).\notag\\
\mb_{ms}& = s\mb_{m0} + (1-s)\mb_{m1},~~ {\mb_{m0}}\equiv  \mb\left(\frac{\mx_0+\mx}{2}\right),~~\mb_{m1}\equiv \mb\left(\frac{\mx_1+\mx}{2}\right) .\notag
\end{align}
Taking the derivative of 
$$
f(s): =  su_0 + (1-s)u_1 +\frac{1}{6} \{ (\|\mb_s\| +4\|\mb_{ms}\|+\|\mb\|)\| \mx - \mx_s\| \notag \\
 - (\mb_s +4\mb_{ms}+\mb)\cdot(\mx - \mx_s) \} 
$$
with respect to $s$ and setting it to zero, we obtain the following equation for $s$:
\begin{align}
& u_0-u_1 + \frac{1}{6}\{ (\|\mb_s\| +4\|\mb_{ms}\|+\|\mb\|) \frac{ (\mx-\mx_s)\cdot(\mx_1-\mx_0) }{ \|\mx - \mx_s\| }  +  \notag \\
& \|\mx - \mx_s\|\left[ 4 \frac{\mb_{ms}\cdot (\mb_{m0}-\mb_{m1}) }{\|\mb_{ms}\|} +  \frac{\mb_{s}\cdot (\mb_{0}-\mb_{1}) }{\|\mb_{s}\|}\right]\label{OS0} \\
& - (\mb_s +4\mb_{ms}+\mb) \cdot (\mx _1- \mx_0 ) - (\mx-\mx_s)\cdot (4(\mb_{m0}-\mb_{m1}) + (\mb_0-\mb_1))  \}  =0. \notag
\end{align}
The hybrid nonlinear solver \cite{wilkinson,stewart} is used for finding a root $s^{\ast}$ of Eq. \eqref{OS0} in the interval $[0,1]$.
In the case of success, the {triangle update} returns
\begin{align*}
{\sf Q}_{\Delta}(\mx_1,\mx_0,\mx) & = s^{\ast} u_0 + (1-s^{\ast}) u_1  +\frac{1}{6} \{ (\|\mb_{s^{\ast}}\| +4\|\mb_{ms^{\ast}}\|+\|\mb\|)\| \mx - \mx_{s^{\ast}}\| \notag \\
& - (\mb_{s^{\ast}} +4\mb_{ms^{\ast}}+\mb)\cdot(\mx - \mx_{s^{\ast}}) \}.
\end{align*}
Otherwise, it returns ${\sf Q}_{\Delta}(\mx_1,\mx_0,\mx) =+\infty$.

%%%%%%%%%%%%%%%%%%%%%%%%%%%%%
%
%%%%%%%%%   A P P E N D I X   C
%
%%%%%%%%%%%%%%%%%%%%%%%%%%%%%

\appendix
\setcounter{equation}{0}
\renewcommand{\theequation}{C-\arabic{equation}}
\section*{Appendix C. Proof of Theorem \ref{thm_equiv}}

\begin{proof}
Without the loss of generality we assume that $\mx_1$ is the origin.

\begin{figure*}
\begin{center}

  \includegraphics[width=0.7\textwidth]{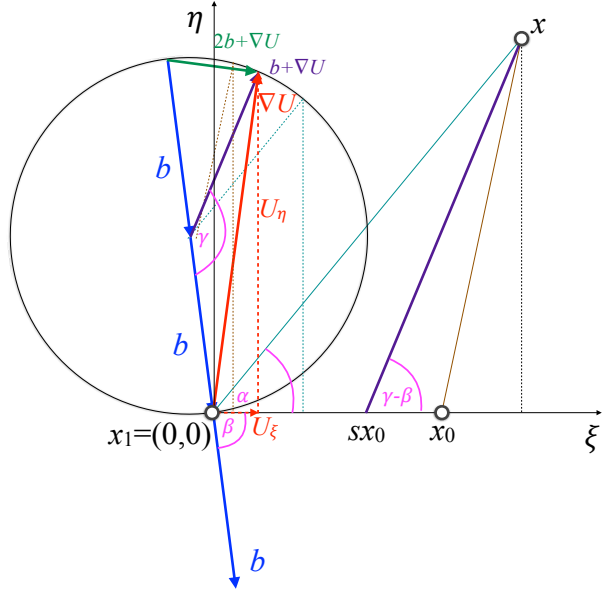}
\caption{An illustration for Section \ref{sec:equiv} and Appendix C.
A geometrical interpretation of the solution of finite difference Eq. \eqref{dfeq} and minimization problem \eqref{r_olim}. }
\end{center}
\label{fig:equiv}       % Give a unique label
\end{figure*}

{\bf Step 1.} Show that $u$ is a solution of Eq. \eqref{dfeq} if and only if $u - u_1= \|\mx\|(U_{\xi}\cos(\alpha) + U_{\eta}\sin(\alpha))$ where (see Fig. \ref{fig:equiv})
$\alpha$ ($0<\alpha <\pi$) is the angle between the vectors $\mx_0$ and $\mx$, $U_{\xi} = \|\mx_0\|^{-1}(u_0-u_1)$, and
$U_{\eta}$ is a solution of 
\begin{equation}
\label{dfeq2}
U_{\xi}^2 + U_{\eta}^2 +2(b_{\xi}U_{\xi} + b_{\eta}U_{\eta}) = 0,\quad
\mb = \left[\begin{array}{c}b_{\xi}\\b_{\eta}\end{array}\right] \equiv \left[\begin{array}{r}\|\mb\|\cos(\beta) \\ -\|\mb\|\sin(\beta)\end{array}\right],
\end{equation}
 which is Eq. \eqref{HJ} written in the $(\xi,\eta)$-coordinates at the point $\mx$.

First observe that both Eqs. \eqref{dfeq} and \eqref{r_olim} are invariant with respect to translations. Therefore,
we shift $\mx_1$ to the origin as shown in Fig. \ref{fig:equiv} without changing their solutions. 

Second, Eq. \eqref{dfeq} is invariant with respect to orthogonal transformations. Indeed, the multiplication of $\mx$ and $\mx_0$ by an orthogonal matrix $O$
converts Eq. \eqref{eq1} to 
\begin{equation}
\label{eq3}
\left[\begin{array}{c}u - u_0\\u-u_1\end{array}\right] = \left[\begin{array}{c}(\mx - \mx_0)^T\\(\mx-\mx_1)^T\end{array}\right]O^{T}\nabla U = PO^T\nabla U.
\end{equation}
Hence the matrix $P$ in Eq. \eqref{dfeq} changes to $PO^T$ and $\mb$ becomes $O\mb$ leading to the equation
\begin{equation}
\label{dfeq1}
\left[u - u_0, u-u_1\right] P^{-T}O^TOP^{-1}\left[\begin{array}{c}u - u_0\\u-u_1\end{array}\right] + 
2\mb^TO^TOP^{-1}\left[\begin{array}{c}u - u_0\\u-u_1\end{array}\right] = 0,
\end{equation}
which is equivalent to Eq. \eqref{dfeq}.
Hence, we apply an orthogonal transformation to map the original coordinate system onto the $(\xi,\eta)$ system in which 
$\mx_0 $ lies on the positive $\xi$-semiaxis  and the $\eta$-coordinate of $\mx$ is positive:
$$
\mx_0 = \left[\begin{array}{c}\|\mx_0\|\\0\end{array}\right],\quad \mx = \left[\begin{array}{c}\|\mx\|\cos(\alpha)\\\|\mx\|\sin(\alpha)\end{array}\right],
$$
where $\alpha$ ($0<\alpha <\pi$) is the angle between vectors $\mx_0$ and $\mx$ as shown in Fig. \ref{fig:equiv}.

Finally, if $u$ is a solution of Eq. \eqref{dfeq} then 
\begin{align*}
\nabla u &= \left[\begin{array}{cc}
 \|\mx\|\cos(\alpha) -\|\mx_0\|~~~ &~~~ \|\mx\|\sin(\alpha) \\ \|\mx\|\cos(\alpha)~~ &~~ \|\mx\|\sin(\alpha) \end{array}\right]^{-1}
\left[\begin{array}{c}u - u_0\\u-u_1\end{array}\right] \\
&= \frac{1}{\|\mx\|\|\mx_0\|\sin(\alpha)}
\left[\begin{array}{cc}
- \|\mx\|\sin(\alpha) ~~~&~~~ \|\mx\|\sin(\alpha) \\
\|\mx\|\cos(\alpha)  ~~~&~~~ -\|\mx\|\cos(\alpha) +\|\mx_0\|
\end{array}\right]
\left[\begin{array}{c}u - u_0\\u-u_1\end{array}\right] \\
&= \left[\begin{array}{c} \frac{u_0-u_1}{ \|\mx_0\|} \\ \frac{(u_1 - u_0)\cos(\alpha)}{\|\mx_0\|\sin(\alpha)} + \frac{u-u_1}{\|\mx\|\sin(\alpha)}
\end{array}\right]\equiv
\left[\begin{array}{c}U_{\xi}\\U_{\eta}\end{array}\right].
\end{align*}
Hence, if $u$ is the solution of Eq. \eqref{dfeq}, then $U_{\xi}$ is exactly $(u_0-u_1)/\|\mx_0\|$ which shows that it is independent of $u$.
Hence Eq. \eqref{dfeq} can be rewritten as an equation Eq. \eqref{dfeq2} for $U_{\eta}$. 

{\bf Step 2.} 
Find geometric conditions guaranteeing the existence of  solution(s) of Eq. \eqref{dfeq2} satisfying the consistency check and  determine the selection rule if it has two solutions.

 Eq. \eqref{HJ} implies that $\nabla U$ is orthogonal to $2\mb + \nabla U$. 
 Therefore, the locus of the vectors $\nabla U$ satisfying Eq. \eqref{HJ} is the circle \cite{quasi} shown in Fig. \ref{fig:equiv}.
 This circle passes through the origin and has center at the end of the vector $-\mb$ originating from the origin. 
 Since $\|\nabla U\|^2 = U_{\xi}^2 + U_{\eta}^2$,  Eq. \eqref{HJ} has a solution if and only 
 if the line normal to the $\xi$-axis and passing through the point $(U_{\xi},0)$ (the red dashed line in Fig. \ref{fig:equiv})
intersects the circle. The MAP is collinear to the vector $\mb+\nabla U$ \cite{quasi}. 
The consistency condition requires that 
the MAP passing through the point $\mx$ crosses the interval $[\mx_1,\mx_0]$. 
This means that the angle between the vector $\mb+\nabla U$ and the positive $\xi$-semiaxis 
should be not less { than} the angle $\alpha$ between the vector $\mx-\mx_1\equiv \mx$ and the positive $\xi$-semiaxis, 
and not greater than the angle between the vector $\mx-\mx_0$ and the positive $\xi$-semiaxis.
Drawing rays parallel to  $\mx$ and $\mx-\mx_0$ from the center of the circle and then dropping normals from their intersections with the circle 
to the $\xi$-axis as shown in Fig. \ref{fig:equiv},  we obtain the interval on the $\xi$-axis where $U_{\xi}$ should belong in order to make the solution $U_{\eta}$ 
of Eq. \eqref{dfeq2} satisfy the
consistency condition. This interval bounded by the endpoints of the thin brown and green-blue dashed lines in Fig. \ref{fig:equiv}.
Note that the consistency condition can be satisfied only by the larger root of Eq. \eqref{dfeq2}, i.e., we should select the root
\begin{align}
U_{\eta} &= -b_{\eta} +\sqrt{b_{\eta}^2 - 2b_{\xi}U_{\xi}-U_{\xi}^2}\notag \\
&\equiv \|\mb\|\sin(\beta)+\sqrt{\|\mb\|^2 \sin^2(\beta)- 2\|\mb\|\cos(\beta)U_{\xi}-U_{\xi}^2}.\label{ueta}
\end{align}

{\bf Step 3.} 
Find the solution of the minimization problem \eqref{r_olim} and show that, if the minimizer $s^{\ast}\in(0,1)$ then it coincides with 
$u = \|\mx\|(U_{\xi}\cos(\alpha) + U_{\eta}\sin(\alpha))$, where $U_{\xi} = (u_0-u_1)/\|\mx\|$ and $U_{\eta}$ is given by Eq. \eqref{ueta}. 

Consider the function to be minimized in Eq. \eqref{r_olim} rewritten for $\mx_1$  shifted to the origin:
\begin{align}
f(s)&:= u_1 + s(u_0-u_1) +\|\mb\|\|\mx - s\mx_0\| - \mb\cdot(\mx-s\mx_0) \notag \\
& \equiv u_1 + U_{\xi}s\|\mx_0\| + \|\mb\|\|\mx - s\mx_0\|(1 - \cos(\gamma)), \label{fun1}
\end{align}
where $\gamma$ is the angle between the vectors $\mb$ and $\mx-s\mx_0$. The point $s\mx_0$, and hence the value of $s$, 
is uniquely determined by the angle $\gamma$ (Fig. \ref{fig:equiv}):
\begin{equation}
\label{s}
\|\mx-s\mx_0\| = \frac{\|\mx\|\sin(\alpha)}{\sin(\gamma-\beta)},\quad s\|\mx_0\| = \|\mx\|\left(\cos(\alpha) - \sin(\alpha)\cot(\gamma - \beta)\right).
\end{equation}
Moreover, since $\cot(\gamma -\beta)$ is a monotone function on the interval $0<\gamma - \beta < \pi$, 
there is a one-to-one correspondence between $-\infty < s<\infty$
and $\beta < \gamma < \beta + \pi$.
Therefore, the function $f(s) =: F(\gamma(s))$, $\beta < \gamma < \beta + \pi$, where
\begin{align}
F(\gamma) &= u_1 + U_{\xi} \|\mx\|\left(\cos(\alpha) - \sin(\alpha)\cot(\gamma - \beta)\right) + 
\frac{\|\mb\|\|\mx\|\sin(\alpha)}{\sin(\gamma-\beta)}(1 - \cos(\gamma)) \notag \\
& = \|\mx\|\left(U_{\xi}\cos(\alpha) +\left[\frac{\|\mb\|(1-\cos(\gamma))}{\sin(\gamma-\beta)} 
-\frac{U_{\xi}\cos(\gamma -\beta)}{\sin(\gamma-\beta)}\right]\sin(\alpha)\right).
\label{Fgamma}
\end{align}
If $(s^{\ast},f(s^{\ast}))$ is a minimum of $f(s)$, then there is a unique minimum $(\gamma^{\ast},F(\gamma^{\ast})=f(s^{\ast}))$ of $F(\gamma)$.

Let us minimize $F(\gamma)$. Its derivative is given by:
\begin{align*}
\frac{dF}{d\gamma}&  = \|\mx\|\sin(\alpha) \frac{[\|\mb\|\sin(\gamma) +U_{\xi}\sin(\gamma - \beta)]}{\sin(\gamma - \beta)}  \notag \\
& - \frac{[\|\mb\|(1-\cos(\gamma)) - U_{\xi}\cos(\gamma -\beta)]\cos(\gamma - \beta)}{\sin^2(\gamma - \beta)}.
\end{align*}
Setting it to zero, cancelling the positive constant $\|\mx\|\sin(\alpha)/\sin^2(\gamma - \beta)$, 
regrouping the terms, and applying trigonometric formulas, we obtain the following equation for $\gamma$:
\begin{equation}
\label{geq}
U_{\xi} + \|\mb\|\cos(\beta) -\|\mb\|\cos(\gamma - \beta) = 0.
\end{equation}
Hence, the optimal angle $\gamma$ satisfies:
\begin{equation}
\label{cosg-b}
\cos(\gamma - \beta) = \frac{U_{\xi} +\|\mb\|\cos(\beta)}{\|\mb\|} = \frac{U_{\xi} +b_{\xi}}{\|\mb\|}.
\end{equation}
Let us denote  by $\gamma^{\ast}$ the solution of Eq. \eqref{cosg-b} lying in the interval $(\beta,\beta + \pi)$.
To check whether $\gamma^{\ast}$ is a maximizer or a minimizer, we evaluate the second derivative of $F(\gamma)$ at $\gamma^{\ast}$ and find:
\begin{equation}
\label{der2F}
\frac{d^2F(\gamma^{\ast})}{d\gamma^2} = \|\mx\|\sin(\alpha)\frac{\|\mb\|}{\sin(\gamma -\beta)} >0,
\end{equation}
as the angle $\gamma - \beta\in (0,\pi)$ by construction. Hence the optimal $\gamma$ is the minimizer of $F$.
{ Next, we recall Eq. \eqref{HJ}: $\|\nabla U\| + 2\mb\cdot\nabla U=0$. 
Adding $\|\mb\|^2$ to both sides, we obtain $\|\nabla U + \mb\|^2 = \|\mb\|^2$. 
Then Eq. \eqref{cosg-b} and the equality $\|\nabla U + \mb\|=\|\mb\|$ imply}
\begin{equation}
\label{sing-b}
\sin(\gamma - \beta) = \frac{U_{\eta} +b_{\eta}}{\|\mb\|}= \frac{U_{\eta} -\|\mb\|\sin(\beta)}{\|\mb\|} .
\end{equation}
Therefore, 
\begin{equation}
\label{ueta2}
U_{\eta} = -b_{\eta}+ \|\mb\|\sin(\gamma - \beta).
\end{equation}
On the other hand, from Eq. \eqref{cosg-b} we obtain:
\begin{equation}
\label{sing-b2}
\sin(\gamma - \beta)= \frac{\sqrt{\|\mb\|\sin^2(\beta) - 2\|\mb\|U_{\xi}\cos(\beta) -U_{\xi}^2}}{\|\mb\|}.
\end{equation}
Plugging Eq. \eqref{sing-b2}  into Eq. \eqref{ueta2} we get
\begin{equation}
\label{ueta3}
U_{\eta} = \|\mb\|\sin(\beta)+ \sqrt{\|\mb\|\sin^2(\beta) - 2\|\mb\|U_{\xi}\cos(\beta) -U_{\xi}^2},
\end{equation}
which coincides with Eq. \eqref{ueta}.

Finally, the solution of the minimization problem \eqref{r_olim}
$$
u =  { \min_{s\in[0,1]}f(s) }
$$
is achieved either at $s^{\ast}$ if $0\le s^{\ast}\le 1$, or at the endpoints $s=0$ or $s=1$. 
Hence, if $0 < s^{\ast} < 1$, then the solution of the minimization problem \eqref{r_olim} coincides with the one of the finite difference scheme \eqref{dfeq}, and the 
latter meets the consistency conditions.
Conversely, the solution of the finite difference scheme \eqref{dfeq} satisfying  the 
consistency conditions coincides with the one of the minimization problem \eqref{r_olim}, and the corresponding minimizer $s^{\ast}\in[0,1]$.

\end{proof}



\end{document}